\documentclass[mnsc,nonblindrev]{informs3_hide}
\OneAndAHalfSpacedXI % default for MS

\usepackage[english]{babel}
\usepackage[autostyle, english = american]{csquotes}
\MakeOuterQuote{"}
\usepackage[colorlinks,linkcolor=blue, citecolor=blue]{hyperref}
\usepackage[normalem]{ulem}

\usepackage{cleveref}
\crefname{subsection}{section}{subsections}
\usepackage{eqnarray}
\usepackage{algorithm,algpseudocode}
\usepackage{amsmath, bbm, enumitem, xparse, bm, lipsum}
\usepackage{amssymb}
\usepackage{comment}

\usepackage{listings}
\usepackage{xcolor} % Required for colors
\usepackage{tcolorbox}

\DeclareRobustCommand{\mybox}[2][gray!10]{%
\begin{tcolorbox}[   %% Adjust the following parameters at will.
        left=0pt,
        right=0pt,
        top=0pt,
        bottom=0pt,
        colback=#1,
        colframe=#1,
        enlarge left by=0mm,
        boxsep=15pt,
        arc=5pt,outer arc=5pt,
        ]
        #2
\end{tcolorbox}
}

\newcommand{\eps}{\varepsilon}
\newcommand{\bI}{\mathbbm{1}}
\newcommand{\bE}{\mathbb{E}}

\newcommand{\Pois}{\mathrm{Pois}}
\newcommand{\ALG}{\mathsf{ALG}}

\newcommand{\by}{\bm{y}}

\newcommand{\bD}{\bm{D}}

\newcommand{\bA}{\bm{A}}

\newcommand{\bx}{\bm{x}}

\newcommand{\PoisOPTRe}{\mathrm{PoisOPTRe}}

\newcommand{\Dual}{\mathrm{Dual}}

\usepackage[dvipsnames]{xcolor}

\usepackage{xparse}

\NewDocumentEnvironment{myproof}{o}
{\IfNoValueTF{#1}{\paragraph{{Proof.} }} {\paragraph{{#1.} }} }
{\hfill$\Halmos$}

\usepackage{natbib,bbm,multirow,multicol}
 \bibpunct[, ]{(}{)}{,}{a}{}{,}%
 \def\bibfont{\small}%
\TheoremsNumberedThrough     % Preferred (Theorem 1, Lemma 1, Theorem 2)
\ECRepeatTheorems

\EquationsNumberedThrough    % Default: (1), (2), ...
\allowdisplaybreaks

\begin{document}
%%%%%%%%%%%%%%%%

%\setcounter{page}{1}
%\DOI{10.1287/xxxx.0000.0000}%

\RUNAUTHOR{Jiang}

\RUNTITLE{Competitive Analysis of Stock-based Pricing}

\TITLE{  Competitive Analysis of Stock-based Thresholds via Prophet Inequalities in Continuous Time
}

\ARTICLEAUTHORS{%
\AUTHOR{Jiashuo Jiang}

\AFF{\  \\
Department of Industrial Engineering and Decision Analytics, Hong Kong University of Science and Technology
}
}

\ABSTRACT{
We study a continuous-time $K$-unit online resource allocation problem with nonhomogeneous Poisson arrivals and time-varying valuation distributions. While the optimal dynamic policy generally depends on both the remaining inventory and the time left in the horizon, we focus on a simpler and practically appealing class of stock-based threshold policies, whose limited number of thresholds depend only on the number of units remaining. We evaluate these policies against the multi-unit prophet benchmark, which selects the best $K$ realized values in hindsight. Our main contribution is a new competitive-analysis framework for stock-based thresholds in continuous time. We first reformulate the problem through a type-covering dual.
The central challenge for analyzing the dual is that the dual is both \textit{infinite-dimensional} and \textit{non-convex}: the adversary can choose time-varying arrival and valuation processes, while the policy performance depends nonlinearly on the stochastic inventory trajectory. 

We overcome these challenges by reducing the continuous-time adversarial problem to a Poisson optimization $\PoisOPTRe_K$, and then proving sharp structural properties of its worst-case solutions. In particular, adversarial arrivals admit cutoff and late-filling structures, which yield an exact four-parameter formulation for two thresholds and a finite nested-interval representation for general thresholds.
These reductions make the guarantees directly computable. For example, for two thresholds, we obtain a competitive ratio $0.6269$ for $K=2$; with three thresholds, we obtain the ratio $0.6816$ for $K=3$. In this way, we show that simple stock-based thresholds achieve strong prophet-inequality guarantees despite ignoring calendar time.
}



\maketitle
%%%%%%%%%%%%%%%%%%%%%%%%%%%%%%%%%%%%%%%%%%%%%%%%%%%%%%%%%%%%%%%%%%%%%%

\section{Introduction}
The allocation of scarce resources to sequentially arriving agents is a fundamental problem in operations research, revenue management, and online mechanism design. A seller begins with $K$ units of inventory and faces customers who arrive over a finite horizon with uncertain valuations. Upon each arrival, the seller must decide immediately whether to allocate one unit and collect the realized reward, or reject the customer and preserve inventory for future, potentially higher-value demand. This trade-off appears in many canonical applications, including hotel room booking (e.g. \cite{bitran1995application}), seasonal retail and fashion inventory management (e.g. \cite{fisher1996reducing, gallego1994optimal}), cloud-computing and reusable-resource services (e.g. \cite{besbes2022static, elmachtoub2025power}), and online advertising allocation \citep{mehta2007adwords}. We study this problem in a continuous-time model with Poisson arrivals, a classical framework in revenue management and dynamic pricing \citep{gallego1994optimal, gallego1997multiproduct}, while allowing arrival rates and valuation distributions to vary over time.

Although dynamic programming characterizes optimal policies when the model is fully specified, the resulting policies typically depend on both the remaining inventory and the precise time left in the horizon. Such policies can be difficult to implement and sensitive to misspecified demand forecasts, especially when demand is non-stationary or exhibits unpredictable spikes and lulls. Motivated by a growing literature in operations research and operations management on the effectiveness of simple policies for online problems (e.g. \cite{ma2021dynamic, elmachtoub2025power}), we focus on stock-based threshold policies: the decision maker posts thresholds that depend only on the number of units remaining, not on calendar time. These policies are simple but not naive. They use inventory as a real-time summary of realized demand, automatically becoming more selective when inventory is depleted quickly and less restrictive when sales are slow. They are also operationally attractive because they resemble booking limits, protection levels, and scarcity-based rules that are easy to communicate and implement. Moreover, by tying decisions to scarcity rather than to arbitrary calendar-time fluctuations, stock-based policies can become more robust to forecast errors and more acceptable to customers than finely time-dependent pricing rules.

We analyze stock-based threshold policies through the lens of prophet inequalities \citep{krengel1978semiamarts}. The offline benchmark is a ``prophet'' who observes all future arrivals and valuations and selects the best $K$ customers in hindsight. An online policy, in contrast, must make irrevocable allocation decisions as customers arrive. The performance of a policy is measured by its competitive ratio: the worst-case ratio between the policy's expected reward and the prophet's expected reward. This benchmark provides a robust way to evaluate online resource allocation policies and is closely connected to posted pricing and online mechanism design (e.g. \cite{hajiaghayi2007automated, correa2019pricing}). Within this framework, our central question is: 

\mybox{How much performance is lost by ignoring time, and what is the best guarantee achievable by stock-based threshold policies?}

\subsection{Our Main Results and Contributions}

We develop a new framework for analyzing stock-based threshold policies in the $K$-unit online resource allocation problem with Poisson arrivals. Our results should be interpreted in two layers. First, we reduce the original minimax problem to a Poisson optimization problem, denoted by $\PoisOPTRe_K(\bm{s})$, where $\bm{s}$ is a vector specifying the number of remaining units that we will switch to a new threshold. The value of $\PoisOPTRe_K(\bm{s})$ provides a valid lower bound on the competitive ratio of stock-based threshold policies. Second, and more importantly, we give an exact structural characterization of the optimal solution to this reduced problem. This structural characterization turns an infinite-dimensional adversarial optimization problem into a finite-dimensional one, which enables us to compute strong guarantees for policies with limited number of stock-based thresholds.

Our first contribution is a reduction from the original competitive analysis problem to $\PoisOPTRe_K(\bm{s})$. The original problem is difficult because the adversary may choose time-varying arrival rates and valuation distributions over a continuous horizon. After dualizing the policy evaluation problem, the resulting formulation contains infinitely many type-covering constraints, and the interaction between the online policy's stochastic inventory path and the prophet's hindsight benchmark creates \emph{non-linear} and \emph{non-convex} structure. We show that, for any switching vector $\bm{s}$, it is enough to consider a reduced Poisson optimization problem with finitely many covering constraints; see \eqref{lp:PoissonDualRe}. The reduced problem captures the central tension faced by stock-based policies: setting thresholds too low may cause the policy to stock out before high-value customers arrive, whereas setting thresholds too high may lead the policy to reject too much demand.

Our second contribution is an exact structural characterization of the optimal solution to the reduced problem $\PoisOPTRe_K(\bm{s})$ for arbitrary $m$, where $m$ denotes the number of thresholds. Although $\PoisOPTRe_K(\bm{s})$ has only finitely many constraints, it still contains continuous-time control variables describing the adversary's arrival pattern, and it remains \emph{non-convex}. We prove that there exists an optimal solution with an extreme, cutoff structure: at every time, the arrival-rate vector across threshold levels is an extreme point, meaning that lower-index threshold levels are active and higher-index levels are inactive. Moreover, the intermediate arrival rates can be chosen in a late-filled form between adjacent threshold levels. After fixing the total masses of these intermediate controls, we further show that the optimal controls admit a finite nested-interval representation: the first nontrivial level is active on a single interval, and each deeper level has only a bounded number of active ``islands'' inside the active time of the previous level. Consequently, the infinite-dimensional optimization over arrival-rate functions reduces to an optimization over $O(m^2)$ number of scalar variables, independent of the inventory level $K$. %We propose a computational approach for solving the general $m$ problem in \Cref{sec:general-switching-computation}.

This structural result leads to an easy way to compute guarantees. For $m=2$, the reduced problem admits an especially compact finite-dimensional representation, involving only $4$ main parameters after normalization. For $m=3$, the nested-interval characterization yields a finite-dimensional program with $9$ timing and mass variables. These formulations can be solved for any fixed inventory level $K$, and the resulting lower bounds are reported in \Cref{tab:Ratios} for $K=1,\ldots,8$. The guarantees improve substantially over the best known ratios for a single static threshold ($m=1$), scaling as $1-\Theta\left(\sqrt{\frac{\log(K)}{K}}\right)$, which was previously developed in \cite{chawla2024static} and shown to be tight in \cite{jiang2025tightness}. Moreover, for small inventory levels, the lower bounds obtained by stock-based policies with only two or three thresholds already exceed the best known guarantees for fully adaptive policies from \cite{jiang2025tight}\footnote{Though the ratios in \cite{jiang2025tight} are developed for the more general Bernoulli arrivals, the worst-case happens when the Bernoulli arrival collapses to the Poisson arrival}. This comparison should be interpreted carefully: the fully adaptive guarantees in \cite{jiang2025tight} are benchmarked through an LP relaxation, whereas our framework works directly with the prophet benchmark.

\begin{table}[H]
    \centering
    \begin{tabular}{|c|c|c|c|c|c|c|c|c|}
    \hline
     Value of $K$    & 1 & 2 & 3 & 4 & 5 & 6 & 7 & 8 \\
     \hline
     Existing Ratios ($m=1$) & 0.5 & 0.5859 & 0.6309 & 0.6605 & 0.6821 & 0.6989 & 0.7125 & 0.7240\\
     \hline
     Our Ratios ($m=2$)    & N.A. & 0.6269 & 0.6732 & 0.7014 & 0.7218 & 0.7376 & 0.7503 & 0.7607\\
     \hline
     Our Ratios ($m=3$) & N.A. & N.A. & 0.6816 & 0.7119 & 0.7324 & 0.7479 & 0.7604 & 0.7707 \\
     \hline
     Existing Ratios (fully adaptive) & 0.5 & 0.6148 & 0.6741 & 0.7120 & 0.7389 & 0.7593 & 0.7754 & 0.7887 \\   
     \hline
    \end{tabular}
    \caption{The comparison between our ratios and existing ratios for $K=1$ up to $K=8$. The existing ratios for static threshold with $m=1$ are developed in \cite{chawla2024static} and have been proven to be tight in \cite{jiang2025tightness}. The existing ratios for fully adaptive policies are from \cite{jiang2025tight}.}
    \label{tab:Ratios}
\end{table}

Overall, our results provide strong evidence for the power of simple stock-based policies. Even though these policies ignore calendar time, a small number of inventory contingent thresholds already achieves better guarantees for stock-based threshold policies, and in some small-inventory regimes exceeds the best known lower bounds for fully adaptive policies, though under a different benchmark. Thus, the price of using simple, inventory-based rules appears to be surprisingly small. Beyond the numerical guarantees, our structural characterization of $\PoisOPTRe_K(s)$ provides a general framework for analyzing and optimizing stock-based policies in continuous-time online resource allocation.

\subsection{Our Techniques}

Our analysis combines ideas from continuous-time duality, prophet inequalities, and optimal control. The main difficulty is that the original problem is both infinite-dimensional and non-convex. The adversary may choose time-varying arrival rates and valuation distributions over a continuous horizon, while the online policy's performance depends on the stochastic evolution of its remaining inventory. Rather than solving this non-convex problem directly, we identify structural properties that any worst-case instance for a reduced problem must satisfy, and use these properties to finally reduce the problem to a finite-dimensional optimization that can be easily optimized over.

\textbf{A continuous-time type-covering formulation.}
We first transform arbitrary valuation distributions into a common type space. This allows us to separate two roles of the adversary: choosing how often different types arrive, and choosing how types are mapped to values. For any fixed arrival structure and stock-based threshold policy, we formulate the competitive ratio problem as an infinite-dimensional linear program over the valuation mapping. Its dual gives a continuous time analogue of the type-covering formulation used in prophet inequalities \citep{jiang2025tightness}. Intuitively, for every type level, the online policy must accept enough higher-type customers to cover a fixed fraction of what the prophet would select in hindsight. This reformulation turns the comparison with the prophet into a family of covering constraints indexed by customer types.

\textbf{Identifying the relevant binding constraints.}
The type-covering dual still contains infinitely many constraints. We reduce this continuum by exploiting the Poisson structure of arrivals. For any type level, the number of arrivals above that level is Poisson distributed, and the marginal value of increasing this arrival mass is governed by a Poisson tail probability: the probability that the prophet has not yet exhausted its $K$ units. On the online side, the corresponding marginal value is governed by the probability that the stock-based policy still has enough inventory to accept that type. Comparing these two marginal quantities at a binding type yields necessary conditions for the binding constraints. These conditions allow us to construct a finite-constraint relaxation/certificate that lower-bounds the original minimax value, leading to the reduced Poisson optimization problem $\PoisOPTRe_K(\bm{s})$, given a fixed switching vector $\bm{s}$.

\textbf{Using optimal-control structure to handle non-convexity.}
The reduced problem $\PoisOPTRe_K(\bm{s})$ has finitely many covering constraints, but it remains a non-convex continuous-time control problem. The non-convexity comes from the interaction between the adversary's arrival-rate controls and the online policy's inventory-state probabilities. We exploit the special state dynamics induced by stock-based threshold policies. Because the threshold depends only on the remaining inventory, the state process has a structured form, and the adversary's controls affect the system through a small number of acceptance rate levels. This structure allows us to analyze the associated Hamiltonian and the geometry of the costate variables. We show that an optimal adversarial control can be chosen to be extremal: at each time, the adversary activates a cutoff set of threshold levels rather than using fractional arrival rates across many levels. Moreover, the intermediate controls can be chosen in a late-filled form between adjacent threshold levels.

To further control the possible switching patterns, we reparametrize time by the active arrival mass of the first threshold level. In this active-time scale, the adversary's control becomes a nested sequence of binary activity indicators. The Hamiltonian switching functions and the geometric shape of the costates imply that the first nontrivial level can be chosen to be active on a single interval, while each deeper level has only a bounded number of active ``islands'' inside the active time of the previous level. This shows that the worst-case arrival pattern is not arbitrary: despite the non-convexity, an optimal adversary admits a bang-bang, nested-interval structure.

\textbf{Reducing to a finite-dimensional computation.}
The structural characterization above turns the infinite-dimensional optimization over arrival-rate functions into an optimization over finitely many scalar variables, independent of the inventory level $K$. For two thresholds, the optimal adversarial pattern reduces to a simple three-block structure, leading to a compact finite-dimensional program. For three thresholds, the nested-interval characterization gives a finite-dimensional program with 9 variables, for each $K$. Solving these reduced programs produces the lower bounds reported in \Cref{tab:Ratios}. Importantly, these computations are not based on a direct time discretization of the original infinite-dimensional problem. Instead, they follow from an exact structural characterization of the reduced Poisson optimization problem.

\subsection{Other Related Literature}
Our work contributes to three stream of literature that i) studies the effectiveness of simple policies, ii) the performances of prophet inequalities, static thresholds, and posted pricing, and iii) dynamic pricing and resource allocation in continuous time, which we now briefly review as follows.

\textbf{Effectiveness of simple policies}. 
A recurring theme in decision making is that carefully chosen low-dimensional policies can retain much of the value of fully dynamic optimization while requiring less information and fewer operational adjustments. In robust revenue management, \citet{ball2009toward} establish tight competitive guarantees for booking-limit policies that become more selective as inventory is depleted. \citet{ma2021dynamic} show that precommitted price and assortment calendars obtain $(1-1/e)$ and $1/2$ guarantees under stationary and nonstationary demand. In service systems, \citet{besbes2022static} and \citet{elmachtoub2025power} establish universal constant-factor guarantees for static pricing even though the optimal policies are state dependent. More closely related to the policy class considered in this paper, \citet{balseiro2026dynamic} show that a stock-dependent policy using only two prices can substantially improve the asymptotic performance loss of static pricing. Together, these results suggest that a small amount of carefully chosen state dependence can capture much of the value of a fully dynamic policy. More broadly, the effectiveness of simple policies has been widely studied in various problems, including stochastic depletion problems \citep{chan2009stochastic}, submodular optimization problems \citep{asadpour2016maximizing, kapralov2013online}, inventory theory \citep{goldberg2016asymptotic, xin2016optimality}, information design \citep{agrawal2026simple}, mechanism design \citep{chawla2010multi, wang2025power}, online matching/advertising \citep{karp1990optimal, mehta2007adwords}.

\textbf{Prophet inequalities, static thresholds, and posted pricing.}
The prophet inequality framework was introduced by \citet{krengel1978semiamarts} and has since become a central tool for analyzing online allocation problems; see, for example, the survey by \citet{correa2019recent}. The literature has developed guarantees under many feasibility constraints, including matroids \citep{kleinberg2012matroid}, knapsacks \citep{feldman2021online}, and multi-unit constraints \citep{hajiaghayi2007automated}. In the multi-unit setting, \citet{alaei2014bayesian} obtained a guarantee of $1-1/\sqrt{k+3}$ using dynamic thresholds, and \citet{jiang2025tight} later improved the guarantees for every $k>1$ and established tightness with respect to the ex-ante benchmark, while \cite{amil2025multi} extends to order fulfillment problems. Static thresholds have also received substantial attention. \citet{hajiaghayi2007automated} showed that static thresholds achieve asymptotically strong guarantees, and \citet{arnosti2021tight} studied tight guarantees in the prophet-secretary setting. More closely related to our work, \citet{chawla2024static} and \citet{jiang2025tightness} analyze static thresholds for multi-unit prophet inequalities, while \citet{chawla2025multi} studies stock-based policies in multi-unit combinatorial prophet inequalities and obtains a competitive ratio that converges to $1-1/e$ as $K\rightarrow\infty$. \citet{perez2026iid} studies threshold policies with limited switches for prophet inequalities under the IID setting.
Prophet inequalities are also closely connected to online pricing and posted-price mechanisms \citep{hajiaghayi2007automated,chawla2010multi,alaei2014bayesian,lucier2017economic, correa2019pricing}. Our results can therefore be interpreted as guarantees for simple posted-price policies with a small number of inventory-contingent prices. These policies are more flexible than a single static threshold but much simpler than fully adaptive policies that depend on both inventory and time.

\textbf{Resource allocation with continuous-time arrival process}. 
Continuous-time stochastic arrival models, particularly Poisson and nonhomogeneous Poisson processes, are foundational in revenue management and pricing. In the seminal work of \citet{gallego1994optimal}, a seller dynamically controls the demand intensity through price while allocating a fixed inventory over a finite horizon. \citet{gallego1997multiproduct} extends this intensity-control framework to multiple products sharing common resources and develops applications to network yield management. Closely related continuous-time models arise in the dynamic and stochastic knapsack problem of \citet{kleywegt1998dynamic}. Subsequent work allows increasingly rich forms of nonstationarity and correlation, for example in \citet{zhao2000optimal, feng2000perishable, bai2023fluid, jiang2025constant, li2025revenue, yuan2026dynamic}.
Continuous-time Poisson models have also become prominent in online matching and service-system resource allocation, where both the availability of resources and the arrival of allocation opportunities evolve stochastically \citep{blanchet2022asymptotically, collina2020dynamic, amanihamedani2024improved}, which leads to the stationary prophet-inequality model of \citet{kessel2022stationary, patel2024combinatorial, aminian2026stationary}.

\section{Problem Formulation}

We consider the resource allocation problem for a single product. There are initially $K$ units of the product. We adopt the standard formulation that customers arrive sequentially in a fixed time interval denoted by $[0,1]$, requesting one unit of the product. We assume that the customer arrives according to a Poisson process over $[0,1]$ and at time $t$, the rate of arrival is given by $\lambda_t\in[0, M]$ for some positive constant $M$ (Here $M$ can be chosen arbitrarily for the worst-case). We denote by $\bm{\lambda}=\{\lambda_t, t\in[0,1]\}$ the sequence of all arrival rates. Each customer arrives with a valuation for the product.
If the customer does not arrive at time $t$, then we simply say that the customer at time $t$ has a valuation $0$. In contrast, if the customer arrives at time $t$,
it is also assumed that the valuation, denoted by $v_t$,  is drawn independently from a given distribution $G_t$, which is positively supported. After the customer at time $t$ reveals its valuation, the decision maker has to decide whether to accept the customer or not. If accepted, then one unit of the product will be consumed, and a reward that equals the valuation $v_t$ will be collected. Otherwise, no unit will be consumed, and no reward will be collected. The goal of the decision maker is to maximize the total collected reward without consuming more than $K$ units of the product.

%If $v_t\geq p_t$, then the customer will buy the product if there is still a remaining unit by paying $p_t$. Otherwise, if $v_t<p_t$, then the customer will simply leave without buying. The price $p_t$ posted by the decision maker is determined based on the problem instance $\bm{G}=\left\{G_t, \forall t\in[0,1]\right\}$ and the historical valuation realization $\left\{v_{\tau}, 0 \leq\tau<t\right\}$, without the awareness of the valuation realization of the customers from time $t$ to the end of the horizon. In this work, we consider the problem of social welfare maximization, i.e., the decision maker sets the prices $\{p_t, t\in[0,1]\}$ to maximize the total valuations of the customers who has purchased the product. 

%\noindent\textbf{Multi-unit Prophet Inequalities}. The dynamic pricing problem under the social welfare maximization objective can be interpreted as a multi-unit prophet inequalities problem in continuous time. Note that in the multi-unit prophet inequalities, the decision maker is allowed to select $K$ values from the valuation set $\left\{v_t, t\in[0,1]\right\}$ in an online manner, and the goal is to maximize the summation of the selected values. This can be done by setting a threshold at each period and selecting the value only if the value exceeds the threshold. Clearly, any pricing policy corresponds to a threshold policy in the multi-unit prophet inequalities problem. %Approximating the optimal dynamic pricing policy is equivalent to approximating the optimal threshold policy for the multi-unit prophet inequalities problem. 

\noindent\textbf{Stock-based Threshold Policies}. Among all possible policies, we restrict to the class of stock-based threshold policies with limited switches. To be specific, after observing the problem instance $\bm{G}=\left\{(G_t, \lambda_t), \forall t\in[0,1]\right\}$, the decision maker selects a threshold vector $\bm{p}=(p_1,\dots, p_m)$ and a switching vector $\bm{s}=(s_1,\dots,s_m)$ such that $K=s_1 >s_2>\dots s_m>s_{m+1}=0$. If we have $k$ units remain and $s_{i+1}<k\leq s_{i}$, the policy accepts an arriving customer if and only if its value is at least $p_i$, subject to the prescribed tie-breaking rule. Note that the threshold vector is fixed before operations begin and is not revised while the inventory level remains unchanged. Thus, the policy reacts to realized demand only through the remaining inventory.

\noindent\textbf{Performance Metrics}. 
We denote by $\ALG_K(\pi^{\bm{p}, \bm{s}}, \bm{G})$ the expected total reward collected by the stock-based threshold policies with the threshold options $\bm{p}$ and switch options $\bm{s}$, under the problem instance $\bm{G}$. We compare the performance of our algorithm to that of the offline optimum, i.e., a prophet who is able to see the value realizations in the future and makes the optimal decision in hindsight. Clearly, for the offline optimum, the optimal decision is made to select the $K$ largest values from the realization $\bm{v}=\{v_t, t\in[0,1]\}$. We denote by $\ALG_K(\pi^{\mathrm{off}}, \bm{G})$ the expected total reward collected by the prophet where $\pi^{\mathrm{off}}$ denotes the offline optimal policy, i.e., 
\begin{equation}\label{eqn:Prophet}
\ALG_K(\pi^{\mathrm{off}}, \bm{G}) = \mathbb{E}_{\bm{v}\sim\bm{G}}\left[ \sum_{k=1}^K v_{(k)} \right]
\end{equation}
where $v_{(1)}\geq v_{(2)}\geq\dots\geq v_{(T)}$ denotes the order statistics of  $\bm{v}=\{v_t, t\in[0,1]\}$.
Then, we aim to derive the optimal stock-based policy represented by $\pi^{\bm{p}^*, \bm{s}^*}$ (here $\bm{p}^*$ and $\bm{s}^*$ are allowed to depend on the instance $\bm{G}$) and characterize the competitive ratio, i.e., the ratio defined by
\begin{equation}\label{def:CompetitiveRatio}
    \gamma^*_{K, m} = \inf_{\bm{G}}\frac{\ALG_K(\pi^{\bm{p}^*, \bm{s}^*}, \bm{G})}{\ALG_K(\pi^{\mathrm{off}}, \bm{G})}.
\end{equation}
In the above definition, we aim to derive a robust guarantee for the performance of our policy. On one hand, we derive the optimal stock-based policy to maximize the ratio of $\frac{\ALG_K(\pi^{\bm{p}^*, \bm{s}^*}, \bm{G})}{\ALG_K(\pi^{\mathrm{off}}, \bm{G})}$, based on the given instance $\bm{G}$. On the other hand, an adversary selects the problem instance $\bm{G}$ to minimize the performance of our policy. The guarantee $\gamma^*_{K, m}$ holds for all problem instances. In this paper, we construct computable lower bounds on $\gamma^*_{K, m}$.

\begin{remark}
It is easy to see that the competitive ratio defined in \eqref{def:CompetitiveRatio} also gives a performance guarantee for the optimal stock-based policies with respect to the optimal policy, since the total expected values collected by the optimal policy is upper bounded by that of the offline optimum. Our results can also be applied to the posted pricing problem, where the decision maker posts prices to sell units to the customers. Following a virtual valuation procedure as detailed in \cite{correa2019pricing}, we can show that the optimal stock-based threshold policy, applied to the transferred virtual valuations, collects at least $\gamma^*_{K, m}$ fraction of the optimal total revenue.
\end{remark}

\section{Stock-based Policies and Problem Reformulation}\label{sec:reformulation}
We now characterize the competitive ratio for the stock-based policies. We derive an optimization problem that could explicitly lower bound the competitive ratios for the stock-based policies. Note that the key difficulty to approximate the competitive ratios is to characterize the worst-case distributions selected by the adversary to minimize the ratio of $\frac{\ALG_K(\pi^{\bm{p}^*, \bm{s}^*}, \bm{G})}{\ALG_K(\pi^{\mathrm{off}}, \bm{G})}$. On one hand, the adversary can select different support sets for the distributions $\bm{G}$. On the other hand, the adversary can select arbitrary probability masses on the support size. To characterize the freedom of the adversary for selecting the support sets and the probability masses, we introduce the two mappings in the following way. 
\begin{definition}[valuation mapping]\label{def:ValueMap}
Let $F:[0,1]\rightarrow\mathcal{G}$ be a strictly monotone mapping from an index set $[0,1]$ to a valuation set $\mathcal{G}$ that covers the support sets of all the distributions in $\bm{G}$. Then, any valuation $v_t$ drawn from the distribution $G_t$, for any $t\in[0, 1]$, can be represented as $v_t=F(q_t)$, for some index $q_t\in[0,1]$.  
\end{definition}

Here, an index $q$ can be understood as the type of a customer, and the valuation $v$ is the corresponding reward for accepting this type of customer.
Note that following the valuation mapping specified in \Cref{def:ValueMap}, a threshold on the valuation can be interpreted as a threshold on the index set $[0,1]$ since $F$ is a one-to-one mapping from the index $[0,1]$ to the valuation set $\mathcal{G}$. Also, a valuation picked by the offline optimum can also be interpreted as an index among the set $[0,1]$. Therefore, we can reduce our problem of computing the ratio in \eqref{def:CompetitiveRatio} for arbitrary distributions to a simplified problem for distributions with support set restricted to $[0,1]$, while the adversary enjoys the freedom to select the valuation mapping $F$. To characterize the distributions over the index set $[0,1]$, we introduce the density mapping as follows. 

\begin{definition}[density mapping]\label{def:DensityMap}
For any problem instance $\bm{G}$ and a valuation mapping $F$, let $\bm{G}'=\{G_t', t\in[0,1]\}$ be distributions on the index set $[0,1]$ such that $G_t(F(q))=G'_t(q)$, for any index $q\in[0,1]$ and $t\in[0, 1]$. 
\end{definition}
For notational simplicity, throughout the paper, we normalize the bottom type
so that $F(0)=0$, assume that $F$ is absolutely continuous and strictly
increasing, and write $g_t(q):=\frac{dG'_t(q)}{dq}$.
We assume that $g_t(q)$ is jointly measurable and uniformly bounded and that
its aggregate density $\int_0^1 \lambda_t g_t(q)dt$
is strictly positive on $(0,1)$. The latter condition only rules out
zero-mass type intervals, which may be deleted from the type space.  Value
atoms can be split by an independent uniform tie-breaking mark.
It is easy to see that any distributions $\bm{G}$ can be fully characterized by a valuation mapping $F$ and the corresponding density mapping $\bm{G'}$, where $\bm{G}'$ can be interpreted as distributions over the index set $[0,1]$. As a result, when the adversary chooses the instance $\bm{G}$ to minimize the ratio of $\frac{\ALG_K(\pi^{\bm{p}^*, \bm{s}^*}, \bm{G})}{\ALG_K(\pi^{\mathrm{off}}, \bm{G})}$, it is equivalent to choosing the valuation mapping $F$ and the density mapping $\bm{G}'$. %We first present the analysis under the regularity assumption that $F$ is absolutely continuous and each $G'_t$ admits a bounded density in $q$. The general case follows by standard smoothing/approximation; atoms can be handled by replacing derivatives with Radon-Nikodym derivatives with respect to the corresponding index measure.
Then, the expected total value collected by the offline optimum can be characterized as follows.
\begin{lemma}\label{lem:Offline}
For any problem instance $\bm{G}$ and the corresponding $F$ and $\bm{G}'$, it holds that
\begin{equation}\label{eqn:Offline}
   \ALG_K(\pi^{\mathrm{off}}, \bm{G}) = \int_{q=0}^1 \mathbb{E}\left[ \min\left\{ \Pois\left(\int_{t=0}^1\left( 1-G'_t(q) \right)\cdot\lambda_tdt\right), K\right\} \right] dF(q), 
\end{equation}
where $\Pois\left(\int_{t=0}^1\left( 1-G'_t(q) \right)\cdot\lambda_tdt\right)$ denotes a Poisson random variable with parameter $\int_{t=0}^1\left( 1-G'_t(q) \right)\cdot\lambda_tdt$, for each $q\in[0,1]$.
\end{lemma}

We now formulate an optimization problem that characterizes the ratio of $\frac{\ALG_K(\pi^{\bm{p}^*, \bm{s}^*}, \bm{G})}{\ALG_K(\pi^{\mathrm{off}}, \bm{G})}$. We introduce a decision variable $U_t^k(q)$ that represents the marginal gain of the decision maker when the value of customer $t$ is realized with index $q$, with $k$ units remaining. We introduce a variable $V^k_t$ to represent the value-to-go for the decision maker at period $t$ when there are $k$ units remaining. The optimization problem, which we denote as $\mathrm{OP}(\bm{\tau}, \bm{s}, \bm{G}')$, can be written as follows with initialization $V^k_1=0$ for each $k=1,\dots,K$ and $V^0_t=0$ for each $t\in[0, 1]$.
\begin{subequations} \label{lp:innerPrimal}
\begin{align}
\min\ &V^K_0& \label{eqn:primalObj}
\\ \mathrm{s.t.\ }&\frac{dV_t^k}{dt}=-\int_{q=\tau_{i(k)}}^1 U^k_{t}(q)dG'_t(q)\cdot\lambda_t ~~~~~~~~~~~~~~~~~~~~\forall t\in[0, 1],k\in[K] \label{eqn:dpValueToGo}
\\ &U^k_{t}(q) =\int_{q'=0}^qdF(q')-(V^k_{t}-V^{k-1}_{t}) ~~~~~~~~~~~~~~\forall t\in[0, 1],q\in[0,1],k\in[K] \label{eqn:dpUtility}
\\ &\int_{q=0}^1 \mathbb{E}\left[ \min\left\{ \Pois\left(\int_{t=0}^1\left( 1-G'_t(q) \right)\cdot\lambda_tdt\right), K\right\} \right] dF(q)=1  \label{eqn:optIs1}
\\ &dF(q)\geq0, ~~~~~~~~~~~~~~~~\ \ \ \ \ \ \ \ \ \ \ \ \ \ \ \ \ \ \ ~~~~~~~~~~~ \forall q\in[0, 1]. \label{eqn:primalNonneg}
\end{align}
\end{subequations}
Note that in the above formulation, the threshold $\tau_{i(k)}$ is an index among $[0,1]$ that represents the price $p_{i(k)}$, where $i(k)\in[m]$ is an index satisfying $s_{i(k)+1} < k \leq s_{i(k)}$, for each $k\in[K]$. Also, we regard the valuation mapping $\{dF(q)\}_{\forall q}$ as a decision variable for the adversary to minimize, with a constraint that the expected total value collected by the offline optimum is normalized to $1$ as in \eqref{eqn:optIs1}. In fact, the optimization problem $\mathrm{OP}(\bm{\tau}, \bm{s}, \bm{G}')$ in \eqref{lp:innerPrimal} concerns the stock-based threshold policies that are \textit{oblivious}, as defined in the following way.
\begin{definition}\label{def:ObviousPolicy}
For a stock-based threshold policy characterized by $\bm{\tau}=(\tau_1,\dots, \tau_m)$, we call $\bm{\tau}$ to be oblivious if $\bm{\tau}$ is decided based on $\bm{G}'$ and is independent of $F$.
\end{definition}
In addition to being oblivious, a class of threshold policies that are of particular interest is the class of monotone threshold policies, defined as follows.
\begin{definition}\label{def:MonotoneThreshold}
For a stock-based threshold policy characterized by $\bm{\tau}=(\tau_1,\dots, \tau_m)$, we call $\bm{\tau}$ to be monotone if it satisfies that $\tau_1\leq\tau_2\leq\dots\leq\tau_m\leq \tau_{m+1}=1$.
\end{definition}
The above definition implies that we will set a low threshold $\tau_1$ at the very beginning. Then, every time the remaining capacity has been decreased to a level $s_i$ for some $i\in[m]$, we switch to a threshold $\tau_i$ that is slightly higher than $\tau_{i-1}$. 
From now on, we restrict to the class of oblivious and monotone threshold policies as defined in \Cref{def:ObviousPolicy} and \Cref{def:MonotoneThreshold}. %In fact, restricting to the class of oblivious and monotone threshold policies will not influence the worst-case guarantee, as we will show later in . 
We denote by $\Gamma$ the set of monotone threshold policies, which include both $\bm{\tau}$ and $\bm{s}$.
The dual problem of the optimization problem $\mathrm{OP}(\bm{\tau}, \bm{s}, \bm{G}')$ in \eqref{lp:innerPrimal} can be given as follows, which we denote as $\mathrm{Dual}(\bm{\tau}, \bm{s}, \bm{G}')$.
\begin{subequations} \label{lp:innerDual}
\begin{align}
\max\ & \theta \  \label{eqn:dualObj}
\\ 
 \mathrm{s.t. }~& \theta\cdot \mathbb{E}\left[ \min\left\{ \Pois\left(\int_{t=0}^1\left( 1-G'_t(q) \right)\cdot\lambda_tdt\right), K\right\} \right] \leq\sum_{k=1}^K\int_{t=0}^1\int_{q'=q}^1 y^k_{t}(q')dq'dt ~~\forall q\in[0,1] \label{eqn:dualCover}
\\ 
& y^k_{t}(q) = \frac{dG'_t(q)}{dq}\cdot\lambda_t\cdot x^k_t ~~~~~~~~~~~~~~~~~~~~~~~~~~~~~~~~~~~~~~~~~\ \ \ \ \forall t\in[0, 1], q\in[\tau_{i(k)},1],k\in[K] \label{eqn:dualUB1}
\\
& y^k_{t}(q) = 0 ~~~~~~~~~~~~~~~~~~~~~~~~~~~~~~~~~~~~~~~~~~~~~~~~~~~~~~\ \ \ \ \ \ \ \   \forall t\in[0, 1], q\in[0, \tau_{i(k)}),k\in[K] \label{eqn:dualUB2}
\\
& x^K_0 = 1, x^k_0 = 0~~~~~~~~~~~~~~~~~~~~~~~~~~~~~~~~~~~~~~~~~~~~~~~~~~~~~~~~ \forall k<K \label{eqn:dualUB3}\\
& \frac{dx^k_t}{dt} =-\int_{q=0}^1 (y^k_{t}(q)-y^{k+1}_{t}(q))dq ~~~~~~~~~~~~~~~~~~~~~~~~~~~~~
~~\forall t\in[0, 1],k\in[K] \label{eqn:dualUpdate}
\\ 
& y^k_{t}(q) \geq0 ~~~~~~~~~~~~~~~~~~~~~~~~~~~~~~~~~~~~~~~~~~~~~~~~\ \ \ \ \ \ \ \ \ \ \ \ \ \forall t\in[0, 1],q\in[0,1],k\in[K]. \label{eqn:dualNonneg}
\end{align}
\end{subequations}
In the above dual formulation, the variable $x^k_t$ can be interpreted as the ex-ante probability of the event that there are $k$ units remaining at period $t$. The variable $y^k_t(q)$ can be interpreted as the ex-ante probability (divided by $dq$) for the decision maker to accept the customer $t$ with value index $q$ when there are $k$ units remaining at period $t$. The dynamics for the variables $x^k_t$ and $y^k_t(q)$ is described in \eqref{eqn:dualUB1}, \eqref{eqn:dualUB2}, \eqref{eqn:dualUB3} and \eqref{eqn:dualUpdate}. A crucial constraint is given in \eqref{eqn:dualCover} which ensures that for each value index $q\in[0,1]$, the number of customers accepted by the decision maker with value index $q'$ satisfying $q'\geq q$ is at least $\theta$ fraction of that of the offline optimum. Then, the ratio $\theta$ gives the competitive ratio guarantee. The dual formulation $\mathrm{Dual}(\bm{\tau}, \bm{s}, \bm{G}')$ can be interpreted as a ``type-cover problem'' and has been derived in \cite{jiang2025tightness} for static-threshold policies.

In the following lemma, we show that the competitive ratio for the oblivious stock-based threshold policies can be obtained from solving the optimization problem $\mathrm{OP}(\bm{\tau}, \bm{s}, \bm{G}')$ in \eqref{lp:innerPrimal} and its dual formulation $\mathrm{Dual}(\bm{\tau}, \bm{s}, \bm{G}')$ in \eqref{lp:innerDual}.
\begin{lemma}\label{lem:TightOblivious}
Fix the class $\Gamma$ of monotone oblivious stock-based threshold policies. 
The worst-case guarantee certified by this more restricted class can be obtained from the following minimax problem
\begin{equation}\label{eqn:minimaxOblivious}
\underline{\gamma}_{K,m}=\min_{\bm{G}'}\max_{(\bm{\tau}, \bm{s})\in\Gamma} \mathrm{OP}(\bm{\tau}, \bm{s}, \bm{G}')=\min_{\bm{G}'}\max_{(\bm{\tau}, \bm{s})\in\Gamma}\mathrm{Dual}(\bm{\tau}, \bm{s}, \bm{G}').
\end{equation}
Moreover, it holds that $\gamma^*_{K,m}\geq \underline{\gamma}_{K,m}$.
\end{lemma}
%Though the current derivation is for oblivious threshold policies, we can follow a similar manner to derive the minimax optimization formulation for threshold policies that are non-oblivious, i.e., the thresholds $\bm{\tau}$ are allowed to depend on both $\bm{G}'$ and $F$. However, we will show later than () the competitive ratios for oblivious stock-based threshold policies and non-oblivious stock-based threshold policies are in fact equivalent to each other in the worst case. Therefore, 
The above \Cref{lem:TightOblivious} implies that $\underline{\gamma}_{K,m}$ gives a valid lower bound for the competitive ratio $\gamma^*_{K, m}$.
From now on, we focus on solving the minimax optimization problem in \eqref{eqn:minimaxOblivious}.

\subsection{Illustration under IID arrivals}\label{sec:IIDcharac}
In order to solve the minimax optimization problem in \eqref{eqn:minimaxOblivious}, we provide further characterization over the structures for $\mathrm{OP}(\bm{\tau}, \bm{s}, \bm{G}')$ and $\mathrm{Dual}(\bm{\tau}, \bm{s}, \bm{G}')$, for fixed $\bm{\tau}, \bm{s}$ and $\bm{G}'$. 
To illustrate our main ideas, we now assume that $\bm{G}'$ is an IID arrivals in that $G'_t = G'_{t'}$ for any $t, t'\in[T]$.

\begin{figure}[!h]
    \centering
    \includegraphics[width=0.6\linewidth]{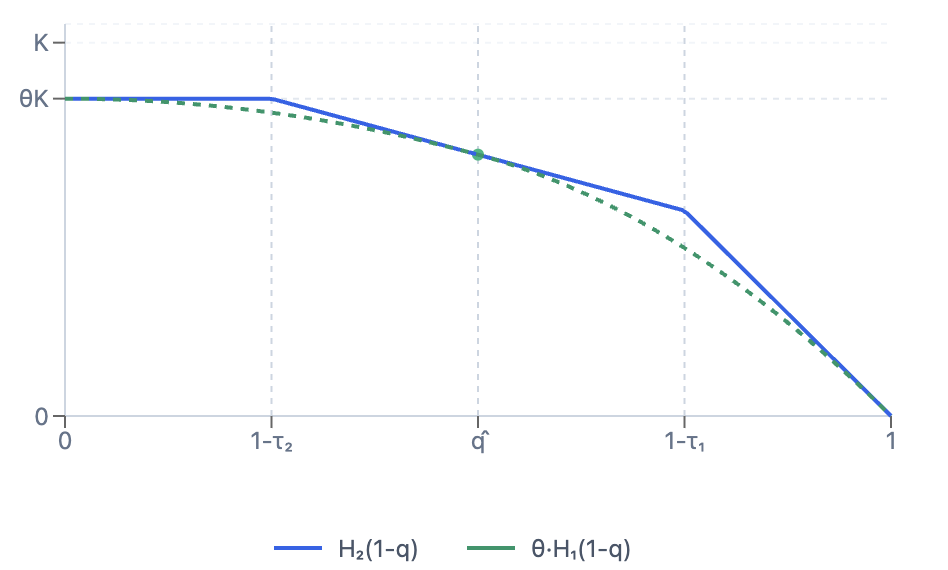}
    \caption{The visualization of $\theta\cdot H_1(q)$ and $H_2(q)$ when $K=2$.}
    \label{fig:K2plot}
\end{figure}

The key point is to analyze the type-cover constraint in \eqref{eqn:dualCover}. We now define a function $H_1$ such that
\[
H_1(q) = \mathbb{E}\left[ \min\left\{ \Pois\left(\int_{t=0}^1\left( 1-G'_t(q) \right)\cdot\lambda_tdt\right), K\right\} \right],~~~\forall q\in[0,1],
\]
which denotes the expected number of customers with value index above $q$ accepted by the offline optimum, for any $q\in[0,1]$. We also denote a function $H_2$ such that
\[
H_2(q) = \sum_{k=1}^K\int_{t=0}^1\int_{q'=q}^1 y^k_{t}(q')dq'dt,~~~\forall q\in[0,1],
\]
for an optimal solution $\{y^k_{t}(q')\}$ to $\mathrm{Dual}(\bm{\tau}, \bm{s}, \bm{G}')$, which denotes the expected number of customers with value index above $q$ accepted by our policy, for any $q\in[0,1]$. Then, the type-cover constraint in \eqref{eqn:dualCover} implies that the following set of constraints are satisfied for a feasible ratio $\theta$,
\begin{equation}\label{eqn:dualCover2}
    \theta\cdot H_1(q)\leq H_2(q).
\end{equation}
It is useful to analyze the two functions $H_1(q)$ and $H_2(q)$. We plot the shapes of the two functions $H_1(q)$ and $H_2(q)$ in \Cref{fig:K2plot}. As we can see, the function $H_1(q)$ is monotone decreasing in $q$, and in fact, it is easy to show that $H_1(q)$ is a concave function over $q$. Moreover, we can show that the function $H_2(q)$ is a piece-wise linear function over $q$, with the turning points given by the specified index thresholds $\bm{\tau}$. 

In order to maximize the ratio $\theta$, it is clear that the constraint in \eqref{eqn:dualCover2} will become binding for some index $q$, which also implies that the function $H_2(q)$ will intersect with the function $\theta\cdot H_1(q)$ at some points, as plotted in . It is easy to see that characterizing these ``intersecting points'' is essential to derive the maximal $\theta$. We now denote by $\mathcal{B}$ the set of all ``intersecting points'' $q$ such that $\theta\cdot H_1(q)=H_2(q)$. Note that $\mathcal{B}$ is a set that depends on the value of $\theta$ and the index thresholds $\bm{\tau}$. From \Cref{fig:K2plot}, we observe that, for any $q\in\mathcal{B}$, it holds that
    \begin{equation}\label{eqn:IIDcondition}
        \theta\cdot H_1(q)=H_2(q)\text{~~and~~}\theta\cdot H_1'(q)=H_2'(q).
    \end{equation}
The condition above is easy to understand as in order for the function $H_2(q)$ to intersect with $\theta\cdot H_1(q)$ at some point $q$, both their values and their derivatives should be equivalent to each other. We will rely on the condition \eqref{eqn:IIDcondition} to provide characterization of $\mathrm{OP}(\bm{\tau}, \bm{s}, \bm{G}')$ and $\mathrm{Dual}(\bm{\tau}, \bm{s}, \bm{G}')$ to solve the minimax optimization problem in \eqref{eqn:minimaxOblivious}.

\subsection{A Set of Conditions for General Distributions}
We now formalize the conditions from the IID setting in the previous \Cref{sec:IIDcharac} under the general non-IID setting. In this section, we allow $\bm{G}'$ to be general non-stationary distributions, and we characterize the set $\mathcal{B}$ using the conditions in \eqref{eqn:IIDcondition}.

For general distributions $\bm{G}'$, index thresholds $\bm{\tau}$ and switching options $\bm{s}$, we denote by $\mathcal{B}(\bm{\tau}, \bm{s}, \bm{G}')\subset[0,1]$ the set of indexes such that the constraint \eqref{eqn:dualCover} is binding for any index $q\in\mathcal{B}(\bm{\tau}, \bm{s}, \bm{G}')$ under an optimal solution for the problem $\mathrm{Dual}(\bm{\tau}, \bm{s}, \bm{G}')$ in \eqref{lp:innerDual}. We then have the following characterization for the set $\mathcal{B}(\bm{\tau}, \bm{s}, \bm{G}')$.
\begin{lemma}\label{lem:CharacBinding}
For any distributions $\bm{G}'$, index thresholds $\bm{\tau}$ and switching options $\bm{s}$, and any index $q\in\mathcal{B}(\bm{\tau}, \bm{s}, \bm{G}')$, suppose that $q\in[\tau_i, \tau_{i+1}]$ for some $i\in[m]$, then it holds that
\begin{equation}\label{eqn:condition1}
\theta^*\cdot\Pr\left[\Pois\left(\int_{t=0}^1\left( 1-G'_t(q) \right)\cdot\lambda_tdt\right) \leq K-1 \right]\geq \sum_{k'=s_{i+1}+1}^K x^{k'*}_1,
\end{equation}
where we denote by $\left\{ \theta^*, x^{k*}_t, y^{k*}_t(q), \forall t\in[0, 1], k\in[K], q\in[0,1] \right\}$ an optimal solution to $\mathrm{Dual}(\bm{\tau}, \bm{s}, \bm{G}')$ in \eqref{lp:innerDual}.
\end{lemma}

The condition \eqref{eqn:condition1} above enables us to derive a simpler problem to consider, as we further show in the following sections.

\section{A Reduced Poisson Optimization Problem and Characterization}\label{sec:WorstCase}

We are now ready to reduce the problem in \eqref{eqn:minimaxOblivious} as a new optimization problem which enjoys a simpler structure and we characterize the worst-case distributions from this optimization problem. We first present the reduction, concerning a finite set of equations. We then characterize the worst-case distributions, showing that customers are arriving in an extreme mode.

\subsection{Reduction to a Poisson Optimization Problem}
We now provide a reduction of the problem in \eqref{eqn:minimaxOblivious}. The condition \eqref{eqn:condition1} we established in \Cref{lem:CharacBinding} will also be useful to characterize the binding indexes, denoted as $\hat{q}\in[\tau_i, \tau_{i+1}]$ for some $i\in[m-1]$, which enables us to reduce the optimization problem $\mathrm{Dual}(\bm{\tau}, \bm{s}, \bm{G}')$ in \eqref{lp:innerDual} with infinite number of constraints into an optimization problem with a finite number of constraints. 
Also, we have the constraint that the index $\hat{q}$ is binding in that the expected number of customers accepted by the offline optimum is equivalent to the expected number of customers accepted by the policy, i.e.,
\begin{equation}\label{eqn:WorstCondition2}
\begin{aligned}
&\gamma^*_K\cdot \mathbb{E}\left[ \min\left\{ \Pois\left(\int_{t=0}^1\left( 1-G'_t(\hat{q}) \right)\cdot\lambda_tdt\right), K\right\} \right]\\ 
=& \sum_{k'=1}^{s_{i+1}}\int_{t=0}^1(1-G'_t(\tau_{i(k')}))\cdot\lambda_t\cdot x^{k'*}_tdt+  \sum_{k'=s_{i+1}+1}^K\int_{t=0}^1(1-G'_t(\hat{q}))\cdot\lambda_t\cdot x^{k'*}_tdt.
\end{aligned}
\end{equation}
From \eqref{eqn:condition1}, and \eqref{eqn:WorstCondition2}, we are actually concerned with the Poisson arrival rate $\left\{\left( 1-G'_t(\hat{q}) \right)\cdot\lambda_t\right\}_{t\in[0,1]}$ for a binding index $\hat{q}\in[\tau_{i}, \tau_{i+1}]$ for some $i\in[m-1]$, as well as the Poisson arrival rate $\left\{\left( 1-G'_t(\tau_{i(k)}) \right)\cdot\lambda_t\right\}_{t\in[0,1]}$ for a threshold $\tau_{i(k)}$ for each $k\in[K]$. We now use a decision variable $\alpha^i_t$ to denote the arrival rate $\left( 1-G'_t(\hat{q}) \right)\cdot\lambda_t$ for the binding index $\hat{q}\in[\tau_{i}, \tau_{i+1}]$, as well as use a decision variable $\beta^i_t$ to denote the arrival rate $\left( 1-G'_t(\tau_i) \right)\cdot\lambda_t$. Then, the reformulation can be given as follows, which we denote by $\PoisOPTRe_K(\bm{s}, M)$, and for notational convenience, by $\PoisOPTRe_K(\bm{s})$ when $M$ is considered fixed.
\begin{subequations} \label{lp:PoissonDualRe}
\begin{align}
\min_{\theta, \bm{\alpha}, \bm{\beta}}\ & \theta \  \label{eqn:PoissonObjRe}
\\ 
 \mathrm{s.t. }~& \theta\cdot \Pr\left[\Pois\left(\int_{t=0}^1\alpha_t^idt\right)\leq K-1\right]\geq\sum_{k'=s_{i+1}+1}^K x_1^{k'} ~~~~ \forall i\in\{1,\dots, m\} \label{eqn:PoissonCover1Re}
\\ 
& \theta\cdot \mathbb{E}\left[\min\left\{\Pois\left(\int_{t=0}^1\alpha_t^idt\right), K \right\}\right]\geq\sum_{k'=1}^{s_{i+1}}\int_{t=0}^1\beta^{i(k')}_tx^{k'}_tdt+\sum_{k'=s_{i+1}+1}^K\int_{t=0}^1\alpha^{i}_tx^{k'}_tdt~ \forall 0\leq i\leq m-1 \label{eqn:PoissonCover2Re}
\\ 
& y^k_{t} = \beta^{i}_t\cdot x^k_t ~~~~\forall s_{i+1}<k\leq s_i, \forall i\in[m], \forall t\in[0, 1] \label{eqn:PoissonUB1Re}
\\
& x^K_0 = 1, ~~~x^k_0 = 0 ~~~~\forall k<K \label{eqn:PoissonUB2Re}\\
& \frac{dx^k_t}{dt} =- (y^k_{t}-y^{k+1}_{t}), ~~~~\forall t\in[0, 1],k\in[K] \label{eqn:PoissonUpdateRe}
\\
&\alpha^m_t=0, ~~~\alpha^{0}_t=\lambda_t, ~~~0\leq\beta^{i+1}_t\leq\alpha^i_t\leq\beta^i_t\leq\lambda_t\leq M, ~~~~\forall t\in[0,1], i\in\{1,\dots, m-1\}.\label{eqn:poissonChain}
\end{align}
\end{subequations}
We have the following result showing that it is sufficient to focus on the problem $\PoisOPTRe_K(\bm{s})$ in \eqref{lp:PoissonDualRe} to derive a lower bound of the competitive ratio guarantee.
\begin{theorem}\label{thm:PoissonOPT}
For any $K$ and $\bm{s}$, it holds that
\begin{equation}\label{eqn:PoissonOPT}
    \PoisOPTRe_K(\bm{s}) \leq \min_{\bm{G}'}\max_{(\bm{\tau}', \bm{s}')\in\Gamma} \mathrm{OP}(\bm{\tau}', \bm{s}', \bm{G}')=\min_{\bm{G}'}\max_{(\bm{\tau}', \bm{s}')\in\Gamma}\mathrm{Dual}(\bm{\tau}', \bm{s}', \bm{G}'),
\end{equation}
where $\PoisOPTRe_K(\bm{s})$ is defined in \eqref{lp:PoissonDualRe}, $\mathrm{OP}(\bm{\tau}', \bm{s}', \bm{G}')$ is defined in \eqref{lp:innerPrimal}, and $\mathrm{Dual}(\bm{\tau}', \bm{s}', \bm{G}')$ is defined in \eqref{lp:innerDual}.
\end{theorem}
\begin{myproof}[Proof sketch of \Cref{thm:PoissonOPT}]
We provide a sketch of our proof and refer the detailed proof to the appendix. The main difficulty for analyzing $\mathrm{OP}(\bm{\tau}', \bm{s}', \bm{G}')$ and $\mathrm{Dual}(\bm{\tau}', \bm{s}', \bm{G}')$ is that they contain continuum of constraints. Although \Cref{lem:CharacBinding} allows us to restrict to one binding type to represent the continuum of constraints within the index interval $[\tau_i, \tau_{i+1}]$, it is not clear whether each index interval for each $i$ admits such a binding type. Therefore, it can happen that for some interval, there is no binding type and we still need to involve a continuum of constraints in the final formulation. 

The nontrivial part of the proof is to show that the above situation cannot happen and every interval admits a binding type. To establish this result, we
parameterize the threshold vector by the lengths of the intervals, using the variables
\[
    z_i=\tau_{i+1}-\tau_i,
    \qquad i=0,\ldots,m,
\]
so that $\bm z$ belongs to a simplex. For each $i$, we define a set consisting of all interval partitions for which the $i$th interval contains a globally worst type-cover constraint. The key observation is that these
sets satisfy the face-covering condition of the KKM lemma \citep{knaster1929beweis}. Intuitively, on a face of the simplex, some intervals have zero length. If a zero-length interval is a bottleneck, then its endpoint is shared with a neighboring nondegenerate interval on the same face, whose worst constraint cannot be
larger than the bottleneck value. Therefore, some nondegenerate interval on
that face must also be a bottleneck. The KKM lemma then guarantees a threshold vector $\bm\tau^*$ for which every interval simultaneously
attains the same global bottleneck value. Consequently, each nondegenerate
interval $[\tau_i^*,\tau_{i+1}^*]$ contains a binding type, while degenerate
intervals are handled through their common endpoints.

For each interval, we select such a binding type and use its tail-arrival
rate as the corresponding $\alpha^i$ control, while the tail-arrival rates
at the thresholds define the $\beta^i$ controls. The binding equality and the condition \eqref{eqn:condition1} in \Cref{lem:CharacBinding} give
the constraints in $\PoisOPTRe_K(\bm s)$. In this way, taking the minimum over $\bm G'$ and applying \Cref{lem:TightOblivious} completes the proof.
\end{myproof}

In the next section, we focus on analyzing the reduced problem $\PoisOPTRe_K(\bm{s})$ in \eqref{lp:PoissonDualRe}, where \Cref{thm:PoissonOPT} implies that $\PoisOPTRe_K(\bm{s})$ forms a lower bound for the competitive ratio guarantee.

\subsection{Establishing the Optimal Solution as Extremal Arrivals}

We now characterize the optimal solution to the relaxed optimization problem $\PoisOPTRe_K$ in \eqref{lp:PoissonDualRe}. We first show in the following lemma that it is optimal to set $\alpha^i_t=\beta^{i+1}_t$ when $t\in[0, w_i]$ and set $\alpha^i_t=\beta^{i}_t$ when $t\in(w_i, 1]$ for some $w_i\in[0,1]$, for each $i\in[m-1]$.
\begin{lemma}\label{lem:WorstAlpha}
Fix any feasible solution of \eqref{lp:PoissonDualRe}.  For every
$i\in[m-1]$, there is another feasible solution with the same $\theta$, the
same $\bm\beta$, and the same total mass
$\int_0^1\alpha_t^idt$, such that
\[
    \hat\alpha_t^i
    =
    \begin{cases}
        \beta_t^{i+1},&t\le w_i,\\
        \beta_t^i,&t>w_i,
    \end{cases}
\]
for some $w_i\in[0,1]$.  The replacement weakly decreases the right-hand
side of the corresponding constraint \eqref{eqn:PoissonCover2Re}.
%For each $i\in[m-1]$, there exists $w_i\in[0,1]$ such that for an optimal solution $\{\theta^*, \alpha^{i*}_t, \beta^{i*}_t, \forall i\in[m], \forall t\in[0,1]\}$ to $\PoisOPTRe_K(\bm{s})$ in \eqref{lp:PoissonDualRe}, it holds that $\alpha^{i*}_t=\beta^{i+1, *}_t$ when $t\in[0, w_i]$ and $\alpha^{i*}_t=\beta^{i, *}_t$ when $t\in(w_i, 1]$.
\end{lemma}

Note that the above \Cref{lem:WorstAlpha} essentially implies that the customers with index among $[\tau_i, \hat{q})$ stop arriving after time $w_i$, since we have $\alpha^{i}_t=\beta^{i}_t$ when $t\in(w_i, 1]$, where $\hat{q}\in[\tau_{i}, \tau_{i+1}]$ is a binding index, for each $i\in[m-1]$. \Cref{lem:WorstAlpha} also implies that the customers with index among $[\hat{q}, \tau_{i+1})$ only starts arriving after time $w_i$, since we have $\alpha^{i}_t=\beta^{i+1}_t$ when $t\in[0, w_i]$. For simplicity, we broadly classify the customers with index among the interval $[\tau_i, \tau_{i+1})$ into two classes, the large class with index among $[\hat{q}, \tau_{i+1})$ and the small class with index among $[\tau_i, \hat{q})$. We can draw the intuition that the customers for the large class start arriving after the customers for the small class stop arriving. We now characterize the worst-case distributions across different intervals of $[\tau_i, \tau_{i+1})$, for different $i\in[m-1]$. We have the following result.

\begin{lemma}\label{lem:WorstBeta} 
For every $\eps>0$, there exists a feasible solution $\{\theta, \alpha^{i}_t, \beta^{i}_t, \forall i\in[m], \forall t\in[0,1]\}$ of
\eqref{lp:PoissonDualRe} whose objective is at most
$\PoisOPTRe_K(\bm{s})+\eps$ and such that, for almost every $t$, there is an
index $r(t)\in\{1,\ldots,m+1\}$ satisfying
\begin{equation}\label{eqn:CutOff}
    \beta_t^i
    =
    \begin{cases}
        \lambda_t,&i<r(t),\\
        0,&i\ge r(t).
    \end{cases}
\end{equation}
Moreover, the $\alpha$ controls may be chosen in the late-filled form of
\Cref{lem:WorstAlpha}. Consequently, restricting to cutoff-valued and
late-filled controls does not change the infimum of
$\PoisOPTRe_K(\bm s)$.
%There exists an optimal solution $\{\theta^*, \alpha^{i*}_t, \beta^{i, *}_t, \forall i\in[m], \forall t\in[0,1]\}$ to $\PoisOPTRe_K(\bm{s})$ in \eqref{lp:PoissonDualRe} such that for each $t\in[0, 1]$, there exists an index $r(t)\in\{1,2,\dots,m+1\}$ and it holds that
%\begin{equation}\label{eqn:CutOff}
%   \beta^{i, *}_t = \left\{\begin{aligned}
%    &0, &i \geq r(t)\\
%   &\lambda_t, & i < r(t)
%    \end{aligned}\right.
%\end{equation}
\end{lemma}
The above \Cref{lem:WorstBeta} implies that we can focus on the form that $\beta^{i}_t\in\{0, \lambda_t\}$, with a single ``cutoff'' across $i\in[m]$, for each time $t\in[0, 1]$. The proof of \Cref{lem:WorstBeta} follows the observation that it is optimal to restrict $\bm{\beta}_t$ to the extreme point of its feasible region, which takes the formulation as given in \eqref{eqn:CutOff} (though the exact extreme point may alter over time). 
%From \Cref{lem:WorstAlpha} and \Cref{lem:WorstBeta}, we know that it is crucial to characterize the shape of the index function $r(t)$. In the following lemma, we show that when $\bm{\nu}$ satisfies that $\nu^k_t=\frac{\nabla_q G'_t(\hat{q})\cdot\lambda_t}{\int_{t=0}^1\nabla_q G'_t(\hat{q})\cdot\lambda_t dt}$ for a binding index $\hat{q}\in(\tau_k, \tau_{k-1}]$ for each $k\in[K]$ in a worst-case distribution $\bm{G}'$, the index function $r(t)$ should be a monotone decreasing function over $t\in[0, 1]$.
%\begin{lemma}\label{lem:WorstMonotone}
%it
%\end{lemma}
Combining \Cref{lem:WorstAlpha} and \Cref{lem:WorstBeta}, we derive the following structural result over an optimal solution to $\PoisOPTRe_K(\bm{s})$ in \eqref{lp:PoissonDualRe}. 
\begin{theorem}\label{thm:WorstStructure}
For any $\eps>0$,
there exists an feasible solution $\{\theta, \alpha^{i}_t, \beta^{i}_t, \forall i\in[m], \forall t\in[0,1]\}$ to $\PoisOPTRe_K(\bm{s})$ in \eqref{lp:PoissonDualRe} such that $\alpha^{i}_t=\beta^{i+1}_t$ when $t\in[0, w_i]$ and $\alpha^{i}_t=\beta^{i}_t$ when $t\in(w_i, 1]$, for some $w_i\in[0, 1]$, for each $i\in[m-1]$, and it holds that $\beta^{i}_t=0$ when $i\geq r(t)$ and $\beta^{i}_t=\lambda_t$ when $i< r(t)$, for an index $r(t)\in\{1,2,\dots, m+1\}$, for each $t\in[0, 1]$.
\end{theorem}

In what follows, we use the structural characterization in \Cref{thm:WorstStructure} to derive our competitive ratio guarantees for stock-based threshold policies under various conditions.

\section{Exact Solution of the Two-Threshold Reduced Poisson Problem}
\label{sec:TightRatio}

In this section, we compute lower bounds of the guarantees of stock-based threshold policies with the help of $\PoisOPTRe_K(\bm{s})$ in \eqref{lp:PoissonDualRe} and the structural results presented in \Cref{thm:WorstStructure}.
Note that a caveat for using \Cref{thm:WorstStructure} is that in general the shape of the ``cutoff'' index function $r(t)$ is hard to characterize. However, we show that when $m = 2$, the index function $r(t)$ possesses a clear structure characterization which enables us to analytically derive the optimal closed formulation. When $m=2$, the formulation of $\PoisOPTRe_K(\bm{s})$ simplifies to 
\begin{subequations} \label{lp:m2PoissonDualRe}
\begin{align}
\min_{\theta, \bm{\alpha}, \bm{\beta}}\ & \theta
\label{eqn:PoissonObjRe}
\\
 \mathrm{s.t. }~& \theta\geq \sum_{k=1}^K x_{t=1}^{k}
\label{eqn:m2PoissonCover1Re}
\\
& \theta\cdot \Pr\left[\Pois\left(\int_{t=0}^1\alpha_tdt\right)\leq K-1\right]
\geq \sum_{k=s_2+1}^K x_{t=1}^{k}
\label{eqn:m2PoissonCover2Re}
\\
& \theta\cdot \mathbb{E}\left[\min\left\{\Pois\left(\int_{t=0}^1\alpha_tdt\right), K \right\}\right]
\geq
\sum_{k=1}^{s_2}\int_{t=0}^1\beta^{2}_t x^{k}_tdt
+\sum_{k=s_2+1}^K\int_{t=0}^1\alpha_t x^{k}_tdt
\label{eqn:m2PoissonCover3Re}
\\
& \theta\cdot \mathbb{E}\left[\min\left\{\Pois\left(\int_{t=0}^1\lambda_tdt\right), K \right\}\right]
\geq
\sum_{k=1}^{s_2}\int_{t=0}^1\beta^{2}_t x^{k}_tdt
+\sum_{k=s_2+1}^K\int_{t=0}^1\beta^1_t x^{k}_tdt
\label{eqn:m2PoissonCover4Re}
\\
& y^k_{t} = \beta^{1}_t\cdot x^k_t
\quad\forall k>s_2,
\qquad
y^k_{t} = \beta^{2}_t\cdot x^k_t
\quad\forall 1\leq k\leq s_2,\quad \forall t\in[0, 1]
\label{eqn:m2PoissonUB1Re}
\\
& x^K_0 = 1,\qquad x^k_0 = 0 \quad\forall k<K
\label{eqn:m2PoissonUB2Re}
\\
& \frac{dx^k_t}{dt} =- (y^k_{t}-y^{k+1}_{t}),
\quad\forall t\in[0, 1],\ k\in[K]
\label{eqn:m2PoissonUpdateRe}
\\
&0\leq\beta^{2}_t\leq\alpha_t\leq\beta^1_t\leq \lambda_t,
\quad\forall t\in[0,1]\label{eqn:m2chainingconstriant}.
\end{align}
\end{subequations}

Note that the characterization in \Cref{thm:WorstStructure} implies that given the parameters $\{\beta^i_t, \forall i\in[m], \forall t\in[0, 1]\}$, the value of $\alpha_t$ is in fact determined by the switching time $w$, which is again determined by the quantity
\[
\Lambda_1 = \int_{t=0}^1 \alpha_t dt.
\]
Moreover, we can show that  it is optimal to place the masses of $\beta^2_t$ in one interval. 
The intermediate control $\alpha$ must contain this $\beta^2$ interval, but
its remaining mass is filled as late as possible, as implied by \Cref{thm:WorstStructure}. Consequently, $\alpha$ control
may be supported either on one terminal interval or on the union of the
$\beta^2$ interval and a later terminal interval.  
Therefore, we can finally show that the entire optimization problem $\PoisOPTRe_K(\bm{s})$ in \eqref{lp:m2PoissonDualRe} can be reduced equivalently to a finite-parameter optimization problem. Then, a simple grid search procedure can be applied to efficiently solve the reduced finite-parameter optimization problem. The final competitive ratio can be obtained by solving the following reduced problem.

We derive our results via the following steps. We treat $M$ as a fixed constant and we uniformize $\lambda_t$ to be $M$, as supported by the following lemma.
\begin{lemma}\label{lem:uniformization}
For any fixed $M>0$, the worst-case competitive ratio over instances with $0\leq \lambda_t\leq M$ is unchanged if we restrict attention to instances
with $\lambda_t\equiv M$.
\end{lemma}
We therefore fix $M$ and assume throughout this section that
$\lambda_t=M$.  It is convenient to use the active time of $\beta^1$.
Then, from \Cref{thm:WorstStructure}, without changing the infimum of $\PoisOPTRe_K(\bm{s})$, we can restrict $\{\alpha_t, \beta^1_t, \beta^2_t\}$ to only take value $0$ or $\lambda_t$, for each $t\in[0, 1]$. The chaining constraints $0\leq\beta^{2}_t\leq\alpha_t\leq\beta^1_t\leq \lambda_t$ imply that $\beta^2_t$ or $\alpha_t$ can take nonzero value only when $\beta^1_t=\lambda_t=M$. We introduce a variable $S$ that denotes the total time that the $\beta^1$ control is active, i.e., $S:=\int_0^1\beta_t^1dt\in[0,M]$. As a result, it is optimal to restrict $\beta^1_t=\lambda_t=M$ on a fixed interval given by $[0, S/M]$. 
\Cref{thm:WorstStructure} further implies that $\alpha_t$ enjoys a late-filled structure. The remaining caveat is to characterize the structure of the $\beta^2$ control. 
The key result of this section is to show that the $\beta^2$ control will only be active on a single interval, as formalize by the following lemma.
\begin{lemma}\label{lem:BetaLemma}
When $m=2$ and assuming $\lambda_t=M$ for every $t$, there exists an optimal solution $\{\beta^{1,*}_t, \alpha^*_t, \beta^{2,*}_t, \forall t\in[0, 1]\}$ to $\PoisOPTRe_K(\bm{s})$ in \eqref{lp:m2PoissonDualRe} such that $\beta^1_t=\lambda_t=M$ for $t\in[0, S/M]$ and $0$ otherwise, for some $S$, $\alpha_t$ enjoys a late-filled structure as specified in \Cref{thm:WorstStructure}, and $\beta^2$ control is active (takes value equivalent to $M$) on a single interval within $[0, S]$ and $0$ otherwise.
\end{lemma}
Since \Cref{lem:BetaLemma} establishes that $\beta^2$ control is active in one interval, this active range
can be parameterized by two variables $a$ and $b$. As a results, it is sufficient to focus on a reduced problem with 4 variables, given by $\Lambda_1$, $S$, $a$, and $b$. We now introduce the feasible set for these variables as
\begin{equation}\label{eq:Z_set}
\mathcal Z_{2,M}
:=
\left\{
(S,a,b,\Lambda_1):
0\leq a\leq b\leq S\leq M,\quad
b-a\leq\Lambda_1\leq S
\right\}.
\end{equation}
We further define the following quantities for any $c>0$,
\[
    P_K(c):=\Pr[\Pois(c)\leq K-1],
    \qquad
    H_K(c):=\mathbb E[\min\{\Pois(c),K\}].
\]
For any $\bm z=(S,a,b,\Lambda_1)\in\mathcal Z_{2,M}$, we define the corresponding $\beta^2$ control as 
\begin{equation}
\label{eq:m2_D2_four_variable}
    \beta^2_{\bm z}(t):=M\cdot\bI_{[a/M,b/M]}(t),
    \qquad t\in[0,S/M].
\end{equation}
and the corresponding normalized $\alpha$ control can be represented as 

\begin{equation}
\label{eq:m2_c_four_variable}
\alpha_{\bm z}(t)
:=
\begin{cases}
M\cdot\bI_{\left[\frac{S-\Lambda_1}{M}, \frac{S}{M}\right]}(t),&S-\Lambda_1\leq a,\\[1mm]
M\cdot\bI_{\left[\frac{a}{M}, \frac{b}{M}\right]\cup\left[\frac{S-\Lambda_1+b-a}{M}, \frac{S}{M}\right]}(t),&S-\Lambda_1>a.
\end{cases}
\end{equation}
Note that it holds
\[
    \beta^2_{\bm z}(t)\leq \alpha_{\bm z}(t)\leq M,
    \text{~~and~~}
    \int_0^{S/M} \alpha_{\bm z}(t)dt=\Lambda_1.
\]
The two cases in \eqref{eq:m2_c_four_variable} have a simple interpretation.
If $S-\Lambda_1\leq a$, the optional $\alpha$ mass begins before the $\beta^2$ interval,
and $\alpha$ is supported on the terminal interval $\left[\frac{S-\Lambda_1}{M}, \frac{S}{M}\right]$.  If $S-\Lambda_1>a$, the
$\beta^2$ interval forces an early $\alpha$ block, while the remaining optional
$\alpha$ mass is placed on the terminal interval $\left[\frac{S-\Lambda_1+b-a}{M}, \frac{S}{M}\right]$.

Now for $\bm z\in\mathcal Z_{2,M}$, let the state dynamics
$\bx^{\bm z}=\{x_t^{k,\bm z}\}$ satisfy
\begin{equation}
\label{eq:m2_active_state_four_variable}
\frac{d x_t^{k,\bm{z}}}{dt}
=
-r_k^{\bm z}(t)x_t^{k,\bm z}
+r_{k+1}^{\bm z}(t)x_t^{k+1,\bm z},
\qquad t\in[0,S/M],
\end{equation}
where we let
\[
    x_0^{K,\bm z}=1,~~~
    x_0^{k,\bm z}=0,~\forall k<K,~~~\text{and~~}
r_k^{\bm z}(t)
=
\begin{cases}
\beta^2_{\bm z}(t),&1\leq k\leq s_2,\\
M,&s_2<k\leq K,
\end{cases}
\]
with the convention that $r_{K+1}^{\bm z}(t)=x_t^{K+1,\bm z}=0$ for any $t$. We further define the quantities
\begin{equation}\label{eq:m2_four_variable1}
\mathcal R_1(\bm z)
:=
\sum_{k=1}^Kx_{\frac{S}{M}}^{k,\bm z},
~\mathcal R_2(\bm z)
:=
\sum_{k=s_2+1}^Kx_{\frac{S}{M}}^{k,\bm z},
~\mathcal R_3(\bm z)
:=
\int_0^{\frac{S}{M}}
\left[
\beta^2_{\bm z}(t)\cdot\sum_{k=1}^{s_2}x_t^{k,\bm z}
+
\alpha_{\bm z}(t)\cdot\sum_{k=s_2+1}^Kx_t^{k,\bm z}
\right]dt
\end{equation}
as well as 
\begin{equation}\label{eq:m2_four_variable2}
\mathcal R_4(\bm z)
:=
\int_0^{\frac{S}{M}}
\left[
\beta^2_{\bm z}(t)\cdot\sum_{k=1}^{s_2}x_t^{k,\bm z}
+
M\cdot\sum_{k=s_2+1}^Kx_t^{k,\bm z}
\right]dt.
\end{equation}
Finally, we can show that the competitive ratio guarantee can be obtained from the following expression
\begin{equation}
\label{eq:m2_Theta_four_variable}
\Theta_{K,s_2,M}(\bm z)
:=
\max
\left\{
\mathcal R_1(\bm z),
\frac{\mathcal R_2(\bm z)}{P_K(\Lambda_1)},
\frac{\mathcal R_3(\bm z)}{H_K(\Lambda_1)},
\frac{\mathcal R_4(\bm z)}{H_K(M)}
\right\}.
\end{equation}
When $\Lambda_1=0$, feasibility forces
$\beta^2_{\bm z}(t)=\alpha_{\bm z}(t)=0$ for any $t$, so we have
$\mathcal R_3(\bm z)=0$. Note that in this boundary case, we interpret the quantity
$\mathcal R_3(\bm z)/H_K(\Lambda_1)$ as zero.
\begin{theorem}
\label{thm:m2-reduced}
When $m=2$, for every fixed $M>0$ and every
$s_2\in\{1,\ldots,K-1\}$, it holds that
\begin{equation}
\label{eq:m2_four_variable_limit}
\PoisOPTRe_K(\bm s)
=
\min_{\bm z\in\mathcal Z_{2,M}}
\Theta_{K,s_2,M}(\bm z),
\end{equation}
where the set $\mathcal Z_{2,M}$ is given in \eqref{eq:Z_set} and the expression $\Theta_{K,s_2,M}(\bm z)$ is defined in \eqref{eq:m2_Theta_four_variable}.
\end{theorem}

\begin{myproof}[Proof of
\Cref{thm:m2-reduced}]
\label{app:m2_self_contained}
Fix $M>0$. By \Cref{lem:uniformization}, throughout the
proof we take $\lambda_t=M$. Further note that from \Cref{lem:BetaLemma}, it is sufficient to consider the control $\beta^2_t$ such that
$\beta^2_t = \beta^2_{\bm{z}}(t)$
for any $t$, for some $\bm{z}\in\mathcal{Z}_{2,M}$. As a result, the corresponding state dynamics $\bx$ can be expressed exactly as $\bx^{\bm{z}}$, for some $\bm{z}\in\mathcal{Z}_{2,M}$, by noting that we always have $\beta^1_t=M$ for $t\in[0, S/M]$ according to \Cref{lem:BetaLemma}. 

We now specify the formulation for the $\alpha$ control. Note that from \Cref{lem:BetaLemma} and \Cref{thm:WorstStructure}, it is optimal to consider the $\alpha$ control that enjoys a late-filled structure. To be specific, we know that there exists a parameter $w$ such that
\[
    \alpha_t
    =
    \begin{cases}
        \beta_t^{2},&t\le w,\\
        \beta_t^1,&t>w.
    \end{cases}
\]
We further have the constraints that 
$\Lambda_1 = \int_{t=0}^1 \alpha_t dt$.
Therefore, when $S-\Lambda_1\leq a$, we know that $\alpha$-control becomes active before the $\beta^2$-control becomes active, i.e., the parameter $w$ is set such that $w=\frac{S-\Lambda_1}{M}$. On the other hand, when $S-\Lambda_1>a$, we know that the $\beta^2$ control becomes active earlier than the time $\frac{S-\Lambda_1}{M}$, which forces the $\alpha$-control to be active in the time interval $\left[ \frac{a}{M}, \frac{b}{M} \right]$ through the chaining constraint \eqref{eqn:m2chainingconstriant}. The constraint that $\Lambda_1 = \int_{t=0}^1 \alpha_t dt$ finally forces the $\alpha$ control to be active in the time interval $\left[ \frac{S-\Lambda_1+b-a}{M}, \frac{S}{M} \right]$. In this way, we show that it is sufficient to consider the control $\alpha$ such that
$\alpha_t = \alpha_{\bm{z}}(t)$
for any $t$.

We finally express both sides of the constraints \eqref{eqn:m2PoissonCover1Re}, \eqref{eqn:m2PoissonCover2Re}, \eqref{eqn:m2PoissonCover3Re}, and \eqref{eqn:m2PoissonCover4Re}. It is clear to see that the left hand sides can be expressed through $P_K(\Lambda_1)$, $H_k(\Lambda_1)$, and $H_K(M)$, while the right hand sides can be expressed through
\[
\sum_{k=1}^{s_2}\int_{t=0}^1\beta^{2}_t x^{k}_tdt
+\sum_{k=s_2+1}^K\int_{t=0}^1\alpha_t x^{k}_tdt = \mathcal{R}_3(\bm{z})
\]
and
\[
\int_0^{\frac{S}{M}}
\left[
\beta^2_{\bm z}(t)\cdot\sum_{k=1}^{s_2}x_t^{k,\bm z}
+
M\cdot\sum_{k=s_2+1}^Kx_t^{k,\bm z}
\right]dt = \mathcal{R}_4(\bm{z}).
\]
In this way, we prove that the maximum $\theta$ to satisfy the constraints \eqref{eqn:m2PoissonCover1Re}, \eqref{eqn:m2PoissonCover2Re}, \eqref{eqn:m2PoissonCover3Re}, and \eqref{eqn:m2PoissonCover4Re} can be selected as the ratio 
\[
\Theta_{K,s_2,M}(\bm z)
:=
\max
\left\{
\mathcal R_1(\bm z),
\frac{\mathcal R_2(\bm z)}{P_K(\Lambda_1)},
\frac{\mathcal R_3(\bm z)}{H_K(\Lambda_1)},
\frac{\mathcal R_4(\bm z)}{H_K(M)}
\right\}.
\]
Our proof is thus completed.
\end{myproof}

The above \Cref{thm:m2-reduced} enables us to reduce computing $\PoisOPTRe_K(\bm{s})$ into solving a simpler optimization problem. Though the optimization problem is still in general non-convex, it only contains 3 free variables, and thus can be solved efficiently by grid search
Note that the problem can be solved only for a finite $M$, while the competitive ratio guarantee is obtained when $M\rightarrow\infty$ (It is direct to see that the value of $\PoisOPTRe_K(\bm{s}, M)$ decreases with $M$). Nevertheless, we show in the following proposition that $\min_{\bm z\in\mathcal Z_{2,M}}\Theta_{K,s_2,M}(\bm z)$ converges to its limit quite quickly. 
\begin{proposition}\label{prop:m2_finite_M}
Let the set $\mathcal Z_{2,M}$ is given in \eqref{eq:Z_set} and the expression $\Theta_{K,s_2,M}(\bm z)$ is defined in \eqref{eq:m2_Theta_four_variable}. Then, for any $s_2\in[K-1]$, for $M$ large enough, it holds that
\[
\gamma^*_{K, 2} \geq \min_{\bm z\in\mathcal Z_{2,M}}
\Theta_{K,s_2,M}(\bm z) - \frac{K-H_K(M)}{H_K(M)},
\]
where we have
\[
\frac{K-H_K(M)}{H_K(M)}
=
\left(
e^{-M}\displaystyle\sum_{j=0}^{K-1}
(K-j)\frac{M^j}{j!}\right)\Bigg/
\left(
K-
e^{-M}\displaystyle\sum_{j=0}^{K-1}
(K-j)\frac{M^j}{j!}\right)
=
O\left(e^{-M}M^{K-1}\right).
\]
\end{proposition}
For $2\le K\le8$, the bound in
\Cref{prop:m2_finite_M} is at most
$1.41\times10^{-4}$ when $M=20$, and at most
$3.85\times10^{-6}$ when $M=25$. We report the computed values in \Cref{tab:m2-Ratios} for $M=10^6$.

\begin{table}[H]
    \centering
    \begin{tabular}{|c|c|c|c|c|c|c|c|}
    \hline
     Value of $K$     & 2 & 3 & 4 & 5 & 6 & 7 & 8 \\
     \hline
     Our Bounds ($m=2$)    & 0.6269 & 0.6732  & 0.7014  & 0.7218  & 0.7376  & 0.7503  & 0.7607 \\
     \hline
    \end{tabular}
    \caption{Our competitive ratios when $m=2$ for $K=2$ up to $K=8$, for $M=10^6$.}
    \label{tab:m2-Ratios}
\end{table}

\section{Structural Characterization for General Thresholds}
\label{sec:GeneralThresholdStructure}

In this section, we develop a finite-dimensional structural characterization of the optimal beta controls in the relaxed Poisson optimization problem $\PoisOPTRe_K(\bm{s})$ in \eqref{lp:PoissonDualRe} for a general number of thresholds $m$ and a general switching vector $\bm{s}$. The formulation in \eqref{lp:PoissonDualRe} is infinite-dimensional because the controls $\{\alpha^i_t,\beta^i_t: i\in[m],\,t\in[0,1]\}$ are functions over time. The purpose of this section is to show that, after fixing the total alpha masses
\[
A_i:=\int_0^1\alpha^i_t\,dt,
\qquad i=1,\dots,m-1,
\]
one can restrict attention to $\beta$-controls described by finitely many switching times.

The structural intuition is the following. 
By \Cref{lem:WorstBeta}, for every time $t$, the vector $(\beta^1_t,\dots,\beta^m_t)$ has a cutoff structure: the first few beta levels are active, and the remaining beta levels are inactive. Therefore, instead of optimizing over $m$ arbitrary functions, it is enough to understand how this cutoff index changes over time. The key difficulty is to characterize the shape of the cutoff index function $r(t)$ from \Cref{lem:WorstBeta}. We prove that the first nontrivial $\beta$ level $\bm{\beta}^2$ can be chosen to be active on a single interval, and that each deeper $\beta$ level has a bounded number of active ``islands'' inside the active time of the previous level.
The final conclusion is that, for every fixed vector $\bm{A}=(A_1,\dots,A_{m-1})$, the fixed-$\bm{A}$ problem reduces to a finite-dimensional optimization over at most
$4m^2-8m-5$ number of
$\beta$-timing variables, independent of $K$. This result gives a clean characterization of the optimal structure of $\PoisOPTRe_K(\bm{s})$ for general $m$ and enables us to reduce solving $\PoisOPTRe_K(\bm{s})$ into solving an optimization problem with a limited number of variables.

\subsection{Active-time Representation}
We use the extreme-point structure from \Cref{lem:WorstBeta} to develop a new representation of $\PoisOPTRe_K(\bm{s})$ using the active time of the $\beta$-level. After the active-time reparametrization, each normalized beta indicator
$D_i$ is either active, taking value $1$, or inactive, taking value $0$. 
%After the active-time reparametrization, each $\beta$-level is either fully active (equaling $M$) or inactive (equaling $0$).
The original formulation uses calendar time $t\in[0,1]$. In this section, it is convenient to re-parametrize time by the cumulative active arrival mass of the first beta level. By \Cref{lem:WorstBeta}, on every non-idle period, we have $\beta^1_t=\lambda_t$, while on idle periods it holds that $\beta^1_t=0$. Define the active-time clock
\[
S:=\int_0^1\beta^1_sds.
\]
Since the state dynamics in \eqref{eqn:PoissonUpdateRe} depend on the controls only through their integrated intensities, this change of variables preserves the objective and the feasibility constraints. After this reparametrization, we relabel the active-time variable again by $t$. Thus, throughout this section, we let $t\in[0,S]$ denotes the active time corresponding to an optimal $\beta$-control, not the original calendar time.

In this active-time scale, $\beta^1\equiv 1$, and from \Cref{lem:WorstBeta}, we can write the normalized control (later proved in the appendix) that
\[
D_i(t):=\frac{\beta^i_t}{\beta^1_t}\in\{0,1\},
\qquad i=1,\dots,m.
\]
Thus, we have that $D_1(t)=1$,
and the chain constraints in \eqref{eqn:poissonChain} imply that
\[
1=D_1(t)\geq D_2(t)\geq\cdots\geq D_m(t)\geq0.
\]
From the previous conditions, it is convenient to consider the mode structure. We call mode $r\in[m]$ is active at time $t$ by meaning that
\[
D_1(t)=\cdots=D_r(t)=1,
\qquad
D_{r+1}(t)=\cdots=D_m(t)=0.
\]
Therefore, in order to characterize the optimal $\beta$-control, it is equivalent to characterize the optimal $D$-control, where the $\alpha$-control follows a late-filled pattern as described in \Cref{lem:WorstAlpha} and \Cref{thm:WorstStructure}.

\subsection{Structural Characterization}
\label{subsec:collapsed_time_islands}

We now describe the structural form of an optimal beta schedule. 
For every level $i$, define its accumulated active time as
\[
L_i(t):=\int_0^tD_i(s)ds,
\qquad
L_i:=L_i(S).
\]
Thus, we have that
\[
S=L_1\ge L_2\ge\cdots\ge L_m\ge0\qquad \text{and}\qquad
L_{i+1}\le A_i\le L_i,
\qquad i=1,\dots,m-1.
\]
The key observation is that $\beta^i_t$ can only be active when the lower level $\beta^{i-1}_t$ is active. Therefore,
for level $i\geq2$, the natural time scale is not the original active time $t$, but the time during which level $i-1$ is active. To describe this, imagine deleting all time intervals on which $D_{i-1}=0$, and concatenating the remaining pieces. The resulting collapsed timeline has length $L_{i-1}$.
For a beta schedule $\bm{D}=(D_1,\dots,D_m)$, define $J_i(\bm{D})$ to be the smallest integer $J$ such that, after this collapse, the active set of $D_i$ can be written as a union of $J$ intervals. Equivalently, $J_i(\bD)$ is the smallest $J$ for which there exist endpoints
\[
0\le a_{i,1}\le b_{i,1}\le a_{i,2}\le b_{i,2}
\le\cdots\le a_{i,J}\le b_{i,J}\le L_{i-1}
\text{~~such that~~}
D_i(t)
=
D_{i-1}(t)\cdot
\bI_{\left\{
L_{i-1}(t)\in
\bigcup_{\ell=1}^{J}[a_{i,\ell},b_{i,\ell}]
\right\}}.
\]
In words, $J_i(\bD)$ is the number of active islands of beta level $i$, measured inside the active time of beta level $i-1$.
The next lemma gives a general safe bound on the number of beta islands.

\begin{lemma}
\label{lem:gen_safe_component}
Fix $m\ge3$ and a vector
$\bA=(A_1,\dots,A_{m-1})$.
There exists an optimal solution to the fixed-$\bA$ problem with
beta indicators $\bD^*$ such that
$J_m(\bD^*)\le2$,
and, for every intermediate level $i=3,\dots,m-1$, it holds
\[
J_i(\bD^*)\le 4(m-i+1)+2.
\]
\end{lemma}
%\begin{lemma}\label{lem:gen_safe_component}
%For every fixed  vector $\bA=(A_1,\dots,A_{m-1})$, there exists an optimal solution to the fixed-$\bA$ problem with $\beta$-indicators $\bD^*$ such that
%\[
%J_i(\bD^*)\le 2(m-i+1),
%\qquad i=2,\dots,m.
%\]
%\end{lemma}
%The lemma controls how often each beta level can turn on and off. For a fixed level $i$, while level $i-1$ is active, there is a one-dimensional marginal gain curve on the Hamiltonian for also activating level $i$. The active islands of $D_i$ are exactly the regions where this curve favors turning level $i$ on.At the bottom level, this marginal gain curve has a Poisson-tail shape, which is monotone or U-shaped. Hence it can generate at most two active islands. Moving upward, the marginal gain curve is obtained by smoothing lower-level information. Such smoothing is variation-diminishing, so it cannot create many new oscillations. This yields the recursive bound $J_i\le 2(m-i+1)$. The formal proof is conducted by a careful analysis of the geometric shape of the corresponding costate in the Hamiltonian.
The lemma controls how often each deeper beta level can turn on and off. For a fixed level $i$, while level $i-1$ is active, the Hamiltonian defines a one-dimensional full-block switching gap that compares the shallower mode with all deeper modes. The derivative of this full-block gap reduces, after a telescoping calculation, to a boundary marginal term.
At the bottom level, this marginal has an explicit Poisson-tail form whose derivative has at most one zero. Hence the bottom beta level can have at most two active islands. For an intermediate level, the corrected
marginal is obtained from lower-level information through a convolution. A fixed-horizontal-level variation-diminishing argument bounds the number of its nonnegative superlevel components. A vanishing strict regularization rules out singular switching intervals and converts
these superlevel bounds into
$J_i(\bD^*)\le4(m-i+1)+2$.
The perturbations are then sent to zero; compactness of the endpoint representation preserves the same island bounds for an exact optimum of the original problem.

%The safe bound in \Cref{lem:gen_safe_component} gives $J_2(\bD^*)\leq 2(m-1)$, which can be further sharpened to a single interval.
\Cref{lem:gen_safe_component} controls all levels from $3$ through $m$.
The first nontrivial $\beta$ level, $D_2$, admits a sharper and separate
characterization in that it can be chosen to be active on a single interval.
\begin{lemma}
\label{lem:gen_D2_connected}
For every fixed vector $\bA=(A_1,\dots,A_{m-1})$, there exists an optimal solution to the fixed-$\bA$ problem with $\beta$ indicators $\bD^*$ such that
$D_2^*$ is active on a single interval, in addition to the conditions in \Cref{lem:gen_safe_component}.
\end{lemma}
The indicator $D_2$ is the first gate into all lower $\beta$ levels. When $D_2=0$, every deeper level $D_3,\dots,D_m$ is also inactive. Thus, if $D_2$ turned on, then off, and then on again, the middle off-period would mean that closing this gate is optimal only temporarily.
The key point is that, on such a $D_2=0$ gap, the marginal comparison between ``reopen $D_2$'' and ``stay in mode $1$'' cannot form the hump needed for a strict internal gap. It is valley-shaped: it can move from decreasing to increasing, but not the reverse. Hence a strict internal $D_2=0$ gap is impossible and the only possible gap is a flat tie. A vanishing tie-breaking perturbation removes these flat ties, and sending the perturbation to zero gives an optimal solution of the original problem with $D_2$ connected. Together, \Cref{lem:gen_safe_component} and \Cref{lem:gen_D2_connected} give the following structural theorem.
\begin{theorem}
\label{thm:gen_threshold_structure}
Fix $m\ge3$ and
$\bA=(A_1,\dots,A_{m-1})$.
There exists an optimal solution to the fixed-$\bA$ version of
$\PoisOPTRe_K(\bm{s})$, with beta indicators $\bD^*$, such that
$D_2^*$ is active on a single interval,
and, for every $i=3,\dots,m-1$, it holds
$J_i(\bD^*)\le4(m-i+1)+2$.
Consequently, the optimal beta control can be represented using timing variables with a number at most
$4m^2-8m-5$.
\end{theorem}

\begin{myproof}[Proof of \Cref{thm:gen_threshold_structure}]
By \Cref{lem:gen_D2_connected}, there exists an optimal fixed-$\bA$
solution in which $D_2^*$ is active on a single interval. This interval
requires two endpoints.
By \Cref{lem:gen_safe_component}, the bottom level satisfies
$J_m(\bD^*)\le2$,
and hence requires at most four endpoints.
For every intermediate level $i=3,\dots,m-1$, we have
$J_i(\bD^*)\le4(m-i+1)+2$.
Since each active interval requires two endpoints in the active time of
level $i-1$, level $i$ requires at most
\[
2J_i(\bD^*)
\le
8(m-i+1)+4
\]
endpoint variables.
Adding the active horizon $S$, the two endpoints of $D_2^*$, the four
endpoints of the bottom level, and the endpoints of the intermediate
levels gives
\[
1+2+4+
\sum_{i=3}^{m-1}\bigl(8(m-i+1)+4\bigr)=
4m^2-8m-5.
\]
Our proof is thus completed.
\end{myproof}

The above \Cref{thm:gen_threshold_structure} gives a nested-island structure. The first nontrivial $\beta$-level $D_2^*$ is one outer interval, shown by \Cref{lem:gen_D2_connected}. Inside the active time of $D_2^*$, the next level $D_3^*$ may have active islands. Inside the active time of $D_3^*$, the next level $D_4^*$ may again have islands, and so on. The number of possible islands is bounded through \Cref{lem:gen_safe_component}.

\subsection{Finite-parameter Reduced Problem}
\label{subsec:finite_parameter_reduced_problem}

We now write the finite-dimensional problem implied by \Cref{thm:gen_threshold_structure}. 
Let $\mathcal{Z}(\bA)$ be the set of endpoint vectors described as follows.
First choose the active horizon
$S\in[0,M]$.
Then choose one interval
$[a_2,b_2]\subseteq[0,S]$,
and set
$D_2(t)=\mathbf 1_{[a_2,b_2]}(t)$,
which implies that
$L_1=S$ and
$L_2=b_2-a_2$.
For the deeper levels, define
\[
N_i:=
\begin{cases}
4(m-i+1)+2, & i=3,\dots,m-1,\\
2, & i=m.
\end{cases}
\]
For each level $i=3,\dots,m$, choose $N_i$ possibly zero-length
intervals in the active time of level $i-1$:
%For each level $i=3,\dots,m$, we define $N_i:=2(m-i+1)$ and choose $N_i$ possibly zero-length intervals in the active time of level $i-1$:
\[
0\le a_{i,1}\le b_{i,1}\le a_{i,2}\le b_{i,2}
\le\cdots\le a_{i,N_i}\le b_{i,N_i}\le L_{i-1}.
\]
Then we define $D_i$ recursively by
\[
D_i(t)
=
D_{i-1}(t)
\mathbf 1_{\left\{
L_{i-1}(t)\in
\bigcup_{\ell=1}^{N_i}[a_{i,\ell},b_{i,\ell}]
\right\}}.
\]
This construction automatically enforces the constraint
$D_1\ge D_2\ge\cdots\ge D_m$.
The total activity of level $i$ is
$L_i=\int_0^S D_i(t)dt$.
Under the endpoint parametrization above, we have that
\[
L_i=\sum_{\ell=1}^{N_i}(b_{i,\ell}-a_{i,\ell}),
\qquad i=3,\dots,m,
\]
which is required to satisfy the feasibility constraints
\[
L_{i+1}\le A_i\le L_i,
\qquad i=1,\dots,m-1.
\]
Given $\bm{z}\in\mathcal{Z}(\bA)$, let $\bD^{\bm{z}}$ be the beta indicators generated by the endpoint vector $\bm{z}$. Let $\bm{\alpha}^{\bm{z}}$ be the alpha controls obtained from $\bD^{\bm{z}}$ by the late-fill rule in \Cref{lem:WorstAlpha} and \Cref{thm:WorstStructure}. Let $\bx^{\bm{z}}$ be the solution to the state equations in \eqref{eqn:PoissonUpdateRe} under these controls.

Define $\Theta_{\bA}(\bm{z}, \bm{s}, M)$ to be the objective value of the fixed-$\bA$ version of $\PoisOPTRe_K(\bm{s}, M)$ evaluated under $(\bD^{\bm{z}}, \bm{\alpha}^{\bm{z}}, \bx^{\bm{z}})$. Equivalently, $\Theta_{\bA}(\bm{z}, \bm{s}, M)$ is the maximum over the normalized left-hand sides of the covering constraints in \eqref{lp:PoissonDualRe}, evaluated under the controls generated by $\bm{z}$.
Then the fixed-$\bA$ problem is equivalent to the finite-dimensional optimization
\[
\min_{\bm{z}\in\mathcal{Z}(\bA)}\Theta_{\bA}(\bm{z}, \bm{s}, M).
\]
By \Cref{thm:gen_threshold_structure}, this finite-dimensional problem has the same value as the original fixed-$\bA$ infinite-dimensional problem.
The number of beta timing variables is
\[
1+2+\sum_{i=3}^m2N_i
=
4m^2-8m-5.
\]
Here the first variable is the active horizon $S$, the next two variables
are the endpoints of $D_2$, the next four variables describe the at most
two bottom-level islands, and the remaining variables are the endpoints
of the intermediate beta islands.
%The number of beta timing variables is
%\[
%1+2+\sum_{i=3}^m 2N_i=
%1+2+\sum_{i=3}^m4(m-i+1)=2(m-1)(m-2)+3.\]
%Here the first variable is \(S\), the next two variables are the endpoints of \(D_2\), and the remaining variables are the endpoints of the deeper beta islands represented by $N_i$.
If one also optimizes over the alpha masses $\bA$, then the outer problem becomes
\[
\PoisOPTRe_K(\bm{s}, M)=\min_{\bA}\ \min_{\bm{z}\in\mathcal{Z}(\bA)}\Theta_{\bA}(\bm{z}, \bm{s}, M),
\]
with $\bA=(A_1,\dots,A_{m-1})$. This adds additional $m-1$ scalar variables to the fixed-$\bA$ reduced problem and gives the precise value of $\PoisOPTRe_K(\bm{s}, M)$.

\subsubsection{Illustration when $m=3$}
\label{subsec:m3_illustration}
When $m=3$, the general structural bound in
\Cref{thm:gen_threshold_structure} becomes
$4m^2-8m-5=
7$.
%When $m=3$, the general structural bound in \Cref{thm:gen_threshold_structure} becomes $2(m-1)(m-2)+3=7$.
The structure has a simple interpretation. The first nontrivial $\beta$ level $D_2$ is active on one interval. Inside the active time of $D_2$, the third level $D_3$ is allowed to have at most two active islands. Thus the $\beta$ mode can be chosen to follow the mode pattern
$1(23232)1\,0$, where it belongs to mode $0$ when even $\beta^1$ becomes inactive.
Equivalently, the active part of the mode sequence has the form
\[
1^{h_1}
2^{h_2}
3^{h_3}
2^{h_4}
3^{h_5}
2^{h_6}
1^{h_7},
\]
followed by idle time $0^{h_0}$. The two $3$-segments represent the two possible active islands of $D_3$ inside the single $D_2$-interval.

Following \Cref{subsec:finite_parameter_reduced_problem}, one convenient endpoint parametrization is to use the 7 timing variables 
$S, a_2,b_2, a_{3,1},b_{3,1}, a_{3,2},b_{3,2}$ for $\beta$ control. Therefore, for $m=3$ and fixed $(A_1,A_2)$, the infinite-dimensional fixed-$\bA$ problem reduces to the finite-dimensional problem
\[
\min_{\bm{z}\in\mathcal Z(A_1,A_2)}\Theta_{A_1,A_2}(\bm{z}, \bm{s}, M),
\]
where
$\bm{z}=(S,a_2,b_2,a_{3,1},b_{3,1},a_{3,2},b_{3,2})$,
$\mathcal Z(A_1,A_2)$ is the set of endpoint vectors satisfying the ordering and feasibility constraints, and $\Theta_{A_1,A_2}(\bm{z}, \bm{s}, M)$ is the objective value of $\PoisOPTRe_K(\bm{s}, M)$ generated by the corresponding $\beta$-schedule and late-filled $\alpha$-controls. To compute the exact value of $\PoisOPTRe_K(\bm{s}, M)$, we also need to minimize over $\bA$, which in total yields $9$ variables to optimize over. 
The problem can be solved by grid search. Finally, to establish the worst case lower bound, we take $M\rightarrow\infty$ (It is direct to see that the value of $\PoisOPTRe_K(\bm{s}, M)$ decreases with $M$). Similar to \Cref{prop:m2_finite_M}, we obtain the following bound regarding the convergence rate.
\begin{proposition}\label{prop:m3_finite_M}
For any $\bm{s}$, it holds that
\[
\gamma^*_{K, 3} \geq \min_{\bA}\ \min_{\bm{z}\in\mathcal{Z}(\bA)}\Theta_{\bA}(\bm{z}, \bm{s}, M) - \frac{K-H_K(M)}{H_K(M)}.
\]
\end{proposition}
The bound of $\frac{K-H_K(M)}{H_K(M)}$ has been discussed in \Cref{prop:m2_finite_M}.
We report the values in \Cref{tab:m3-Ratios} for $M=10^6$.

\begin{table}[H]
    \centering
    \begin{tabular}{|c|c|c|c|c|c|c|}
    \hline
     Value of $K$     & 3 & 4 & 5 & 6 & 7 & 8 \\
     \hline
     Our Ratios ($m=3$)  & 0.6816 & 0.7119  & 0.7324 & 0.7479 & 0.7604 & 0.7707 \\
     \hline
    \end{tabular}
    \caption{Our competitive ratios when $m=3$ for $K=3$ up to $K=8$, for $M=10^6$.}
    \label{tab:m3-Ratios}
\end{table}

%\input{sec_m3}

%In fact, the entire optimization problem $\PoisOPTRe_K(\bm{s})$ in \eqref{lp:m2PoissonDualRe} can be reformulated as a single-parameter optimization problem over $\Lambda$ under the structural results in \Cref{thm:WorstStructure}. It is easy to see that all the other quantities in \eqref{lp:m2PoissonDualRe} can be expressed using $\Lambda$. Then, the final competitive ratio can be computed via solving

\section{Concluding Remarks}

In this paper, we provided a comprehensive analysis of stock-based pricing policies for the multi-unit prophet inequality. By introducing a novel reduction to a finite-dimensional Poisson optimization problem, we overcame the challenge of the infinite-dimensional adversarial strategy space. Our analysis reveals that the worst-case scenarios for these static policies follow a special structure, characterized by extreme arrival intensities. This structural characterization allowed us to lower bound the guarantee for $m=2$ and $m=3$, for any initial inventory level $K$, establishing explicit worst-case guarantees of static thresholds in non-stationary environments. Broadly, our findings offer a rigorous theoretical justification for the widespread industrial reliance on inventory contingent pricing. the numerical results suggest that the loss from restricting to stock-based thresholds may be modest. Future work may extend this framework to settings with reusable resources or more complex combinatorial constraints, further exploring the boundaries of oblivious decision-making.

\bibliographystyle{abbrvnat}
\bibliography{bibliography}

\clearpage

\OneAndAHalfSpacedXI

% Appendix here
% Options are (1) APPENDIX (with or without general title) or
%             (2) APPENDICES (if it has more than one unrelated sections)
% Outcomment the appropriate case if necessary
%
% \begin{APPENDIX}{<Title of the Appendix>}
% \end{APPENDIX}
%
%   or
%

\begin{APPENDICES}
\crefalias{section}{appendix}

\section{Missing Proofs of \Cref{sec:reformulation}}

\subsection{Proof of \Cref{lem:Offline}}

Let $N$ be the realized number of arrivals and let
$Q_1,\ldots,Q_N\in[0,1]$ be their realized types.  Write
$Q_{(1)}\ge Q_{(2)}\ge\cdots$
for the decreasing order statistics, with the convention $Q_{(j)}=0$ for
$j>N$.  For $q\in[0,1]$, define
$N_q:=\sum_{j=1}^N \bI\{Q_j\ge q\}$ as the number of arrivals with an index greater than or equal to $q$.
Pathwise, the following identity holds:
\begin{equation}\label{eqn:pathwise-layer-cake}
    \sum_{j=1}^K F(Q_{(j)})
    =
    \int_{[0,1]} \min\{N_q,K\}\,dF(q).
\end{equation}
Indeed, an increment $dF(q)$ is counted once for each of the top $K$ types
that is at least $q$, and the number of such types is exactly
$\min\{N_q,K\}$.
It is direct to see that $N_q$ is Poisson with mean
$A(q):=\int_0^1 \lambda_t\bigl(1-G'_t(q)\bigr)dt$.
Taking expectations in \eqref{eqn:pathwise-layer-cake} gives
\[
\begin{aligned}
    \ALG_K(\pi^{\mathrm{off}},\bm G)
    =
    \int_{[0,1]}
    \bE\bigl[\min\{N_q,K\}\bigr] dF(q)=
    \int_{[0,1]}
    \bE\left[
        \min\left\{
            \Pois\left(
                \int_0^1\lambda_t(1-G'_t(q))dt
            \right),K
        \right\}
    \right]dF(q),
\end{aligned}
\]
which is exactly \eqref{eqn:Offline}. Our proof is completed.

\subsection{Proof of \Cref{lem:TightOblivious}}

We prove a pointwise strong-duality statement for every fixed 
$(\bm{\tau},\bm{s},\bm{G}')$, and then take the maximization over
$(\bm{\tau},\bm{s})\in\Gamma$ and the minimization over $\bm{G}'$.
Throughout the proof, the arrival-rate process $\bm{\lambda}$ is fixed and
suppressed in the notation. 

Fix $(\bm{\tau},\bm{s},\bm{G}')$, we write
\[
    g_t(q):=\frac{dG'_t(q)}{dq}.
\]
For a valuation mapping $F$, let $\mu=dF$ be the corresponding nonnegative measure on $[0,1]$, so that
\[
    F(q)=\int_{[0,q]} \mu(dr).
\]
For each $q\in[0,1]$, define the prophet's type-cover function
\[
    H_1(q)
    :=
    \mathbb{E}\left[
    \min\left\{
    \Pois\left(\int_{t=0}^1(1-G'_t(q))\lambda_tdt\right),K
    \right\}
    \right].
\]
This is the expected number of type at least $q$ customers selected by the
prophet. Next let $(x,y)$ be the unique solution of the forward equations
\eqref{eqn:dualUB1}--\eqref{eqn:dualUpdate}, with the convention
$y_t^{K+1}(q)\equiv0$. Define the online type-cover function
\[
    H_2(q)
    :=
    \sum_{k=1}^K
    \int_{t=0}^1
    \int_{q'=q}^1 y^k_t(q')dq'dt .
\]
This is the expected number of type at least $q$ customers accepted by the
stock-based threshold policy $(\bm{\tau},\bm{s})$.

We first show that, after eliminating the value-to-go variables in
\eqref{lp:innerPrimal}, the objective of $\mathrm{OP}(\bm{\tau},\bm{s},\bm{G}')$
is exactly
\[
    \int_{[0,1]} H_2(q)\,\mu(dq).
\]
For $k\in[K]$, define
\[
    a_t^k
    :=
    \int_{q=\tau_{i(k)}}^1 g_t(q)\lambda_t\,dq,
    \qquad
    r_t^k
    :=
    \int_{q=\tau_{i(k)}}^1 F(q)g_t(q)\lambda_t\,dq .
\]
Then the backward equations \eqref{eqn:dpValueToGo}--\eqref{eqn:dpUtility}
can be rewritten as
\[
    \frac{dV_t^k}{dt}
    =
    -r_t^k+a_t^k(V_t^k-V_t^{k-1}),
    \qquad k\in[K],
\]
where $V_t^0=0$. Meanwhile, the forward equations for the stock process are
\[
    \frac{dx_t^k}{dt}
    =
    -a_t^kx_t^k+a_t^{k+1}x_t^{k+1},
    \qquad k\in[K],
\]
with the convention $a_t^{K+1}x_t^{K+1}=0$. Therefore, we have that
\[
\begin{aligned}
    \frac{d}{dt}\sum_{k=1}^K x_t^kV_t^k
    =
    \sum_{k=1}^K
    \left(-a_t^kx_t^k+a_t^{k+1}x_t^{k+1}\right)V_t^k
    +
    \sum_{k=1}^K
    x_t^k\left(-r_t^k+a_t^k(V_t^k-V_t^{k-1})\right)  =
    -\sum_{k=1}^K x_t^k r_t^k .
\end{aligned}
\]
Using $x_0^K=1$, $x_0^k=0$ for $k<K$, and $V_1^k=0$ for all $k$, integrating
over $t\in[0,1]$ gives that
\[
    V_0^K
    =
    \sum_{k=1}^K
    \int_{t=0}^1
    x_t^k
    \int_{q'=\tau_{i(k)}}^1
    F(q')g_t(q')\lambda_t\,dq'dt .
\]
By the definition of $y_t^k(q')$, this becomes
\[
    V_0^K
    =
    \sum_{k=1}^K
    \int_{t=0}^1
    \int_{q'=0}^1
    F(q')y_t^k(q')\,dq'dt .
\]
Since $F(q')=\int_{[0,q']}\mu(dq)$, Fubini's theorem yields that
\[
\begin{aligned}
    V_0^K
    =
    \sum_{k=1}^K
    \int_{t=0}^1
    \int_{q'=0}^1
    \left(\int_{[0,q']}\mu(dq)\right)y_t^k(q')\,dq'dt  =
    \int_{[0,1]}
    \left(
    \sum_{k=1}^K
    \int_{t=0}^1
    \int_{q'=q}^1 y_t^k(q')dq'dt
    \right)\mu(dq)  =
    \int_{[0,1]} H_2(q)\,\mu(dq).
\end{aligned}
\]
By \Cref{lem:Offline}, the prophet value is
\[
    \ALG_K(\pi^{\mathrm{off}},\bm{G})
    =
    \int_{[0,1]} H_1(q)\,\mu(dq).
\]
Therefore, after normalizing the prophet value to one, the inner problem
\eqref{lp:innerPrimal} is equivalent to the one-constraint measure linear
program
\begin{equation}\label{eqn:measureLP-L2}
    P(\bm{\tau},\bm{s},\bm{G}')
    :=
    \inf_{\mu\in\mathcal{M}_+([0,1])}
    \left\{
    \int_{[0,1]} H_2(q)\,\mu(dq):
    \int_{[0,1]} H_1(q)\,\mu(dq)=1
    \right\},
\end{equation}
where $\mathcal{M}_+([0,1])$ denotes the cone of finite nonnegative Borel
measures on $[0,1]$.
The dual of \eqref{eqn:measureLP-L2} is
\begin{equation}\label{eqn:measureDual-L2}
    D(\bm{\tau},\bm{s},\bm{G}')
    :=
    \sup_{\theta\in\mathbb{R}}
    \left\{
    \theta:
    \theta H_1(q)\le H_2(q),\quad \forall q\in[0,1]
    \right\}.
\end{equation}
Indeed, after substituting the unique feasible solution $(x,y)$ of
\eqref{eqn:dualUB1}--\eqref{eqn:dualUpdate}, the constraint
$\theta H_1(q)\le H_2(q)$ is exactly the type-cover constraint
\eqref{eqn:dualCover}. Thus \eqref{eqn:measureDual-L2} is exactly
$\mathrm{Dual}(\bm{\tau},\bm{s},\bm{G}')$.
It remains to prove that \eqref{eqn:measureLP-L2} and
\eqref{eqn:measureDual-L2} have the same value. Let
\[
    \mathcal{Q}_+
    :=
    \{q\in[0,1]:H_1(q)>0\}.
\]
We exclude the trivial case $\mathcal{Q}_+=\emptyset$, in which the prophet
value is zero for every valuation mapping and the normalized ratio is
irrelevant. Define
\[
    \rho
    :=
    \inf_{q\in\mathcal{Q}_+}
    \frac{H_2(q)}{H_1(q)} .
\]
Since $H_1(q)\ge0$ and $H_2(q)\ge0$, we have $\rho\ge0$. Moreover, $\rho$ is
finite whenever the normalized problem is nontrivial.
First consider the dual. If $\theta$ is feasible for
\eqref{eqn:measureDual-L2}, then for every $q\in\mathcal{Q}_+$, we have that
$\theta
    \le
    \frac{H_2(q)}{H_1(q)}$.
Hence $\theta\le\rho$. Conversely, by the definition of $\rho$, we have that
\[
    \rho H_1(q)\le H_2(q),\qquad \forall q\in[0,1],
\]
where the inequality is automatic when $H_1(q)=0$. Thus $\theta=\rho$ is dual
feasible, and the dual value is $\rho$.
Now consider the primal. For any feasible $\mu$ in \eqref{eqn:measureLP-L2}, it holds that
\[
\begin{aligned}
    \int_{[0,1]} H_2(q)\mu(dq)
    =
    \int_{\mathcal{Q}_+}
    \frac{H_2(q)}{H_1(q)}H_1(q)\,\mu(dq)
    +
    \int_{[0,1]\setminus\mathcal{Q}_+}
    H_2(q)\,\mu(dq)  
    \ge
    \rho
    \int_{\mathcal{Q}_+} H_1(q)\,\mu(dq) =
    \rho .
\end{aligned}
\]
Therefore the primal value is at least $\rho$. Conversely, for every
$\epsilon>0$, choose $q_\epsilon\in\mathcal{Q}_+$ such that
\[
    \frac{H_2(q_\epsilon)}{H_1(q_\epsilon)}
    \le
    \rho+\epsilon .
\]
The measure
$\mu_\epsilon
    :=
    \frac{1}{H_1(q_\epsilon)}\delta_{q_\epsilon}$
is feasible for \eqref{eqn:measureLP-L2}, and its objective value is
\[
    \int_{[0,1]} H_2(q)\,\mu_\epsilon(dq)
    =
    \frac{H_2(q_\epsilon)}{H_1(q_\epsilon)}
    \le
    \rho+\epsilon .
\]
Letting $\epsilon\downarrow0$ gives that the primal value is at most $\rho$.
Thus, we have that
\[
    P(\bm{\tau},\bm{s},\bm{G}')
    =
    D(\bm{\tau},\bm{s},\bm{G}').
\]
Equivalently, it holds that
\[
    \mathrm{OP}(\bm{\tau},\bm{s},\bm{G}')
    =
    \mathrm{Dual}(\bm{\tau},\bm{s},\bm{G}')
\]
for every fixed $(\bm{\tau},\bm{s},\bm{G}')$.
%If one insists that $F$ be absolutely continuous, replace the atomic measure $\delta_{q_\epsilon}$ above by a smooth density supported on a sufficiently small one-sided neighborhood of $q_\epsilon$. Under the regularity assumption that each $G'_t$ has a bounded density, both $H_1$ and $H_2$ are continuous in $q$. Hence the objective value of the smoothed measure converges to $H_2(q_\epsilon)/H_1(q_\epsilon)$. Therefore allowing Stieltjes measures does not change the value of the problem.
Since the equality
\[
    \mathrm{OP}(\bm{\tau},\bm{s},\bm{G}')
    =
    \mathrm{Dual}(\bm{\tau},\bm{s},\bm{G}')
\]
holds pointwise for every $(\bm{\tau},\bm{s},\bm{G}')$, we may maximize over
$(\bm{\tau},\bm{s})\in\Gamma$ and then minimize over $\bm{G}'$ to obtain
\[
    \inf_{\bm{G}'}
    \sup_{(\bm{\tau},\bm{s})\in\Gamma}
    \mathrm{OP}(\bm{\tau},\bm{s},\bm{G}')
    =
    \inf_{\bm{G}'}
    \sup_{(\bm{\tau},\bm{s})\in\Gamma}
    \mathrm{Dual}(\bm{\tau},\bm{s},\bm{G}').
\]
This proves the equality part of \eqref{eqn:minimaxOblivious}. Since $\Gamma$ restricts to the class of oblivious and monotone stock-based threshold policies, it is clear to see that $\gamma^*_{K,m}\ge \underline{\gamma}_{K,m}$, completing the proof.

\subsection{Proof of \Cref{lem:CharacBinding}}

For fixed distributions $\bm{G}'$ and index thresholds $\bm{\tau}$, we now define a function  
\[
H_1(q) =\mathbb{E}\left[ \min\left\{ \Pois\left(\int_{t=0}^1\left( 1-G'_t(q) \right)\cdot\lambda_tdt\right), K\right\} \right],~~~\forall q\in[0,1],
\]
which denotes the expected number of customers with value index above $q$ accepted by the offline optimum, for any $q\in[0,1]$. We also denote a function $H_2$ such that
\[
H_2(q) = \sum_{k=1}^K\int_{t=0}^1\int_{q'=q}^1 y^{k*}_{t}(q')dq'dt,~~~\forall q\in[0,1],
\]
for an optimal solution $\left\{ \theta^*, x^{k*}_t, y^{k*}_t(q), \forall t\in[0, 1], k\in[K], q\in[0,1] \right\}$ to $\mathrm{Dual}(\bm{\tau}, \bm{G}')$, which denotes the expected number of customers with value index above $q$ accepted by our policy, for any $q\in[0,1]$. The constraint \eqref{eqn:dualCover} implies that
\begin{equation}\label{eqn:112801}
    \theta^*\cdot H_1(q)\leq H_2(q)
\end{equation}
for any $q\in[0,1]$. 

We now focus on an index $\hat{q}\in\mathcal{B}(\bm{\tau}, \bm{s}, \bm{G}')$ such that $\hat{q}\in[\tau_i, \tau_{i+1}]$ for some $i\in[m-1]$.
Set $\tau_0:=0$ and $\tau_{m+1}:=1$, and define
\[
r
:=
\max\left\{
j\in\{0,\ldots,m\}:\tau_j\le\hat q
\right\}.
\]
Since $\hat q<1$, we have $\hat q<\tau_{r+1}$.
Moreover, because
$\hat q\in[\tau_i,\tau_{i+1}]$, we have $r\ge i$.
Note that $\hat{q}\in\mathcal{B}(\bm{\tau}, \bm{s}, \bm{G}')$ implies that
\begin{equation}\label{eqn:112802}
    \theta^*\cdot H_1(\hat{q}) = H_2(\hat{q}).
\end{equation}
Choose $h>0$ sufficiently small that
$q+h<\tau_{r+1}$. The condition in \eqref{eqn:112801} gives that
\[
\theta^*\cdot H_1(\hat{q}+h)\le H_2(\hat{q}+h).
\]
Subtracting the two relations give that
\begin{equation}\label{eqn:112803}
    H_2(\hat{q})-H_2(\hat{q}+h)
\leq
\theta^*\cdot
\bigl(H_1(\hat{q})-H_1(\hat{q}+h)\bigr)..
\end{equation}
It is easily seen that $H_1$ is continuously differentiable.
We now provide further characterization on $H'_1(\hat{q})$ and $H_2(\hat{q})-H_2(\hat{q}+h)$. In what follows, we define 
\[
\nabla_q G'_t(\hat{q})=\frac{G'_t(\hat{q}+h)-G'_t(\hat{q})}{h}.
\]

\noindent\textbf{Characterization of $H'_1(\hat{q})$}: By plugging in the probability formulation for Poisson random variables, we know that
\[\begin{aligned}
H_1(\hat{q}) =& \underbrace{\sum_{k=1}^K k\cdot \frac{\left( \int_{t=0}^1\left( 1-G'_t(\hat{q}) \right)\cdot\lambda_tdt \right)^k}{k!}\cdot\exp\left( -\int_{t=0}^1\left( 1-G'_t(\hat{q}) \right)\cdot\lambda_tdt\right)}_{I(\hat{q})} \\
&+\underbrace{\sum_{k=K+1}^{\infty} K\cdot \frac{\left( \int_{t=0}^1\left( 1-G'_t(\hat{q}) \right)\cdot\lambda_tdt \right)^k}{k!}\cdot\exp\left( -\int_{t=0}^1\left( 1-G'_t(\hat{q}) \right)\cdot\lambda_tdt\right)}_{II(\hat{q})}.
\end{aligned}\]
It is clear to see that the term $\frac{\left( \int_{t=0}^1\left( 1-G'_t(\hat{q}) \right)\cdot\lambda_tdt \right)^k}{k!}\cdot\exp\left( -\int_{t=0}^1\left( 1-G'_t(\hat{q}) \right)\cdot\lambda_tdt\right)$ exactly express the probability that
\[
 \Pois\left(\int_{t=0}^1\left( 1-G'_t(q) \right)\cdot\lambda_tdt\right) = k ,
\]
which gives the above expression for $H_1(\hat{q})$. We now take the derivative over $\hat{q}$ and analyze each part separately. 

We first consider the part $I$ that concerns $k\leq K$. We have that
\[\begin{aligned}
\nabla_q I=\frac{I(\hat{q}+h)-I(\hat{q})}{h} = &o(h)+ \left( -\int_{t=0}^1 \nabla_q G'_t(\hat{q})\cdot\lambda_tdt \right)\cdot\exp\left( -\int_{t=0}^1\left( 1-G'_t(\hat{q}) \right)\cdot\lambda_tdt\right)\\
&\cdot\left( \sum_{k=0}^{K-1}  \frac{\left( \int_{t=0}^1\left( 1-G'_t(\hat{q}) \right)\cdot\lambda_tdt \right)^{k}}{k!}- K\cdot \frac{\left( \int_{t=0}^1\left( 1-G'_t(\hat{q}) \right)\cdot\lambda_tdt \right)^{K}}{K!}  \right).
\end{aligned}\]
Therefore, we know that
\[\begin{aligned}
\nabla_q I = &o(h) + \left( -\int_{t=0}^1 \nabla_q G'_t(\hat{q})\cdot\lambda_tdt \right)\\
&\cdot\left( \Pr\left[\Pois\left(\int_{t=0}^1\left( 1-G'_t(q) \right)\cdot\lambda_tdt\right) \leq K-1 \right] - K\cdot \Pr\left[\Pois\left(\int_{t=0}^1\left( 1-G'_t(q) \right)\cdot\lambda_tdt\right) = K \right] \right).
\end{aligned}\]
We then consider the $k$ such that $k\geq K+1$. We have that
\[\begin{aligned}
\nabla_q II
=&\frac{II(\hat{q}+h)-II(\hat{q})}{h}= o(h) +  \left( -\int_{t=0}^1 \nabla_q G'_t(\hat{q})\cdot\lambda_tdt \right)\cdot\exp\left( -\int_{t=0}^1\left( 1-G'_t(\hat{q}) \right)\cdot\lambda_tdt\right)\\
&\cdot\left( \sum_{k=K+1}^{\infty}K\cdot  \frac{\left( \int_{t=0}^1\left( 1-G'_t(\hat{q}) \right)\cdot\lambda_tdt \right)^{k-1}}{(k-1)!}- \sum_{k=K+1}^{\infty} K\cdot \frac{\left( \int_{t=0}^1\left( 1-G'_t(\hat{q}) \right)\cdot\lambda_tdt \right)^{k}}{k!}  \right)\\
=&o(h) + \left( -\int_{t=0}^1 \nabla_q G'_t(\hat{q})\cdot\lambda_tdt \right)\cdot\exp\left( -\int_{t=0}^1\left( 1-G'_t(\hat{q}) \right)\cdot\lambda_tdt\right)\cdot K\cdot \frac{\left( \int_{t=0}^1\left( 1-G'_t(\hat{q}) \right)\cdot\lambda_tdt \right)^{K}}{K!}.
\end{aligned}\]
Therefore, we get that 
\[\begin{aligned}
&\frac{H_1(\hat{q}+h)-H_1(\hat{q})}{h}=\nabla_q I +\nabla_q II \\=& o(h) + \left( -\int_{t=0}^1 \nabla_q G'_t(\hat{q})\cdot\lambda_td_t \right)\cdot\exp\left( -\int_{t=0}^1\left( 1-G'_t(\hat{q}) \right)\cdot\lambda_tdt\right)\cdot \sum_{k=0}^{K-1}  \frac{\left( \int_{t=0}^1\left( 1-G'_t(\hat{q}) \right)\cdot\lambda_tdt \right)^{k}}{k!}\\
=&o(h) + \left( -\int_{t=0}^1 \nabla_q G'_t(\hat{q})\cdot\lambda_tdt \right)\cdot\Pr\left[\Pois\left(\int_{t=0}^1\left( 1-G'_t(q) \right)\cdot\lambda_tdt\right) \leq K-1 \right].
\end{aligned}\]
\noindent\textbf{Characterization of $H_2(\hat{q})-H_2(\hat{q}+h)$}: Following the threshold policy rules, there is no policy threshold in $(\hat q,\hat q+h]$.
Hence, every type in this interval is accepted precisely in states
$k>s_{r+1}$. We have that
\[
H_2(\hat{q})-H_2(\hat{q}+h)=
\int_0^1
\lambda_t
\bigl(G'_t(\hat{q}+h)-G'_t(\hat{q})\bigr)
\cdot \sum_{k=s_{r+1}+1}^K x_t^{k*}dt
\nonumber\ge
\sum_{k=s_{r+1}+1}^K x_1^{k*}\cdot
\int_0^1
\lambda_t
\bigl(G'_t(\hat{q}+h)-G'_t(\hat{q})\bigr)dt,
\]
where the second inequality follows by noting that from the update dynamics in \eqref{eqn:dualUB3} and \eqref{eqn:dualUpdate}, we have that
\[
\sum_{k'=s_{r+1}+1}^K x^{k'*}_t = 1- \int_{t'=0}^{t}\int_{q=0}^1 y^{(s_{r+1}+1)*}_{t'}(q)dqdt' \geq 1- \int_{t'=0}^{1}\int_{q=0}^1 y^{(s_{r+1}+1)*}_{t'}(q)dqdt' = \sum_{k'=s_{r+1}+1}^K x^{k'*}_1 
\]
for any $t\in[0, 1]$.
In this way, we show that
\[
H_2(\hat{q})-H_2(\hat{q}+h)\geq h\cdot\sum_{k'=s_{r+1}+1}^Kx^{k'*}_t\cdot\int_{t=0}^1\nabla_q G'_t(\hat{q})\cdot \lambda_t dt\geq h\cdot\sum_{k'=s_{i+1}+1}^Kx^{k'*}_t\cdot\int_{t=0}^1\nabla_q G'_t(\hat{q})\cdot \lambda_t dt,
\]
where the second inequality follows from $r\geq i$.
Therefore, the condition in \eqref{eqn:112803} implies that, by taking $h\rightarrow0$, 
\[
\theta^*\cdot \left(\int_{t=0}^1 \nabla_q G'_t(\hat{q})\cdot\lambda_tdt\right) \cdot\Pr\left[\Pois\left(\int_{t=0}^1\left( 1-G'_t(\hat{q}) \right)\cdot\lambda_tdt\right) \leq K-1 \right] \geq \sum_{k'=s_{i+1}+1}^Kx^{k'*}_t\cdot\int_{t=0}^1\nabla_q G'_t(\hat{q})\cdot \lambda_t dt,
\]
Therefore, we have that
\begin{equation}\label{eqn:112805}
\theta^*\cdot\Pr\left[\Pois\left(\int_{t=0}^1\left( 1-G'_t(\hat{q}) \right)\cdot\lambda_tdt\right) \leq K-1 \right]\geq \sum_{k'=s_{i+1}+1}^K x^{k'*}_1.
\end{equation}
%A further note is that the above proof holds even when $i=m$ such that $\hat{q}=1$. To see this, we can extend the range of the index to $[0, 1+\eps]$ for some small $\eps>0$. Then, every derivations above follows under standard smooth regularity condition of $G'_t$. In this way, we establish \eqref{eqn:112805} even for $i=m$ and $\hat{q}=1$. 
Our proof is thus completed.

\section{Missing Proofs of \Cref{sec:WorstCase}}

\subsection{Proof of \Cref{thm:PoissonOPT}}
Note that from strong duality in \Cref{lem:TightOblivious}, we always have that $\mathrm{OP}(\bm{\tau}', \bm{s}', \bm{G}')=\mathrm{Dual}(\bm{\tau}', \bm{s}', \bm{G}')$. Thus, it is sufficient to show that
\[
\PoisOPTRe_K(\bm{s})\leq \max_{\bm{s}'} \min_{\bm{G}'}\max_{ \bm{\tau}'}\mathrm{Dual}(\bm{\tau}', \bm{s}', \bm{G}') \leq \min_{\bm{G}'}\max_{(\bm{\tau}', \bm{s}')\in\Gamma}\mathrm{Dual}(\bm{\tau}', \bm{s}', \bm{G}').
\]
The second inequality holds from the definition. We only need to prove the first inequality.
In order to prove the equality, we first fix $\bm{s}$ and show that it is optimal to have at least one binding type in each interval $[\tau_i, \tau_{i+1}]$ for each $i$. To be specific, we denote $
\mathcal{T}_m
:=
\left\{
\bm\tau\in[0,1]^m:
0\le\tau_1\le\cdots\le\tau_m\le1
\right\}$ to be the feasible set for the monotone threshold $\bm{\tau}$
and define the $m+1$ threshold cells 
\begin{equation}
    I_i(\bm\tau):=[\tau_i,\tau_{i+1}],
    \qquad i=0,\ldots,m.
    \label{eqn:thresholdCells}
\end{equation}
For each index $q\in[0, 1]$, we further define
\[
H_1(q)=\mathbb{E}\left[ \min\left\{ \Pois\left(\int_{t=0}^1\left( 1-G'_t(q) \right)\cdot\lambda_tdt\right), K\right\} \right]
\text{~~and~~}
H_2(q, \bm{\tau}) = \sum_{k=1}^K\int_{t=0}^1\int_{q'=q}^1 y^k_{t}(q', \bm{\tau})dq'dt
\]
where we denote by $(\bx(\bm{\tau}), \bm{y}(\bm{\tau}))$ the trajectory and the acceptance probabilities under the threshold policies $(\bm{\tau}, \bm{s})$. Under the standing regularity convention, it is easy to show that $H_2(q;\bm\tau)$ is jointly continuous in $(q,\bm\tau)$.
Then, it is clear to see that
\begin{equation}
    \Dual(\bm\tau,\bm s,\bm G')
    =
    \inf_{\{q:H_1(q)>0\}}R(q;\bm\tau)\text{~~~with~~~}R(q;\bm\tau):=
    \frac{H_2(q;\bm\tau)}{H_1(q)}.
    \label{eqn:DualAsMinimumRatio}
\end{equation}
We now introduce a vanishing regularization. For $\eps>0$, define the quantity
\[
    \mathcal R_\eps(q;\bm\tau)
    :=
    \frac{H_2(q;\bm\tau)+\eps}{H_1(q)+\eps},
\]
and, for the cell $I_i(\bm\tau)$, define the bottleneck of cell $i$ by
\begin{equation}\label{eqn:cellBottleneck}
    c_i^\eps(\bm\tau)
    :=
    \min_{q\in I_i(\bm\tau)}
    \mathcal R_\eps(q;\bm\tau),
    \qquad i=0,\ldots,m.
\end{equation}
In other words, $c_i^\eps(\bm\tau)$ refers to the minimum regularized type-cover ratio within the interval $I_i(\bm{\tau})$ under the thresholds $\bm{\tau}$. We now obtain the following result.
\begin{lemma}\label{lem:OneinInterval}
For every fixed environment $\bm{G}'$ and switching vector
$\bm{s}$, for every $\eps>0$, there exists a monotone threshold vector
$\bm{\tau}^\eps\in\mathcal{T}_m$ and a scalar $\theta^\eps$ such that
\begin{equation}
    c^\eps_0(\bm\tau^\eps)
    =
    c^\eps_1(\bm\tau^\eps)
    =
    \cdots
    =
    c^\eps_m(\bm\tau^\eps)
    =
    \theta^\eps.
    \label{eqn:allCellBottlenecksEqual}
\end{equation}
%Moreover, it holds that
%\begin{equation}
%\theta=\inf_{\{q:H_1(q)>0\}}
%R(q;\bm\tau^*)=\Dual(\bm\tau^*,\bm{s},\bm{G}').
%\label{eqn:thetaIsDualValue}
%\end{equation}
\end{lemma}
Note that the above \Cref{lem:OneinInterval} essentially implies that for each interval $[\tau^\eps_i, \tau^\eps_{i+1}]$ there exists one regularized minimizer $q$ (when the interval is degenerate such that $\tau^\eps_i=\tau^\eps_{i+1}$, the regularized minimizer $q$ is the threshold $\tau_i^\eps$ itself), as long as the thresholds $\bm{\tau}^\eps$ are chosen to satisfy the conditions in \Cref{lem:OneinInterval} for each fixed $\bm{G}'$ and $\bm{s}$.

Now for each fixed $\bm{G}'$ and $\bm{s}$, we prove that 
\[
\PoisOPTRe_K(\bm{s})\leq \max_{ \bm{\tau}'}\mathrm{Dual}(\bm{\tau}', \bm{s}, \bm{G}').
\]
For any $\eps>0$, 
fix $\bm{\tau}^\eps$ as the monotone thresholds satisfying the conditions in \Cref{lem:OneinInterval}. Importantly, we know that such $\bm{\tau}^\eps$ exists for every $\bm{G}'$ and $\bm{s}$.
Further from \Cref{lem:OneinInterval}, for each interval $[\tau^\eps_i, \tau^\eps_{i+1}]$ for $i\in[m-1]$, we know that there exists one regularized minimizer $q_i^\eps$. 
Therefore, we know that 
\[
\theta^\eps\cdot\eps+\theta^\eps\cdot \mathbb{E}\left[ \min\left\{ \Pois\left(\int_{t=0}^1\left( 1-G'_t(q^\eps_i) \right)\cdot\lambda_tdt\right), K\right\} \right] =\eps+\sum_{k=1}^K\int_{t=0}^1\int_{q'=q^\eps_i}^1 y^{k}_{t}(q', \bm{\tau}^*)dq'dt.
\]
Moreover, we have the following result as a direct corollary of the condition \eqref{eqn:condition1} in \Cref{lem:CharacBinding},
\[
\theta^\eps \cdot\Pr\left[\Pois\left(\int_{t=0}^1\left( 1-G'_t(q^\eps_i) \right)\cdot\lambda_tdt\right) \leq K-1 \right] \geq \sum_{k'=s_{i+1}+1}^K x^{k'}_1(\bm{\tau}^\eps) dt,~~\forall i=1,\dots,m-1.
\]
When $i=m$, for every $q\in[\tau^\eps_m, 1]$, every customer with an index at least $q$ is accepted as long as there is inventory left, which implies that
\[
H_2(q, \bm{\tau}^\eps)=\int_0^1
\lambda_t
\cdot \left(1-G'_t(q)\right)
\cdot \sum_{k=1}^K x_t^{k}(\bm{\tau}^\eps)dt
\nonumber\ge
\sum_{k=1}^K x_1^{k}(\bm{\tau}^\eps)\cdot
\int_0^1
\lambda_t\cdot \left(1-G'_t(q)\right)dt.
\]
Also, we have that
\[
H_1(q)=\mathbb{E}\left[ \min\left\{ \Pois\left(\int_{t=0}^1\left( 1-G'_t(q) \right)\cdot\lambda_tdt\right), K\right\} \right]\leq \int_{t=0}^1\left( 1-G'_t(q) \right)\cdot\lambda_tdt.
\]
As a result, we have that
\[
H_2(q, \bm{\tau}^\eps) +\eps \geq \sum_{k=1}^K x_1^{k}(\bm{\tau}^\eps)\cdot
\int_0^1
\lambda_t\cdot \left(1-G'_t(q)\right)dt+\eps \geq \sum_{k=1}^K x_1^{k}(\bm{\tau}^\eps)\cdot (H_1(q)+\eps),
\]
which implies that
\[
\theta^\eps=\min_{q\in[\tau^\eps_m, 1]} \mathcal{R}_\eps(q;\bm{\tau}^\eps) = \min_{q\in[\tau^\eps_m, 1]} \frac{H_2(q, \bm{\tau}^\eps) +\eps}{H_1(q)+\eps} \geq \sum_{k=1}^K x_1^{k}(\bm{\tau}^\eps).
\]

We now let $\eps\rightarrow0$ and there exists a sequence $\eps^n\rightarrow0$ such that $\bm{\tau}^{\eps^n}\rightarrow\bm{\tau}^*$, $\theta^{\tau^n}\rightarrow\theta^*$, and $\bm{q}^{\eps^n}\rightarrow\bm{q}^*$, for some $\bm{\tau}^*$, $\theta^*$, and $\bm{q}^*$.
Now, we do the change of variables, by using a decision variable $\alpha^i_t$ to denote the arrival rate $\left( 1-G'_t(q^*_i) \right)\cdot\lambda_t$ for the index $q^*\in[\tau_{i}^*, \tau_{i+1}^*]$, as well as using a decision variable $\beta^i_t$ to denote the arrival rate $\left( 1-G'_t(\tau_i^*) \right)\cdot\lambda_t$. Then, we recover exactly the constraints \eqref{eqn:PoissonCover1Re} and \eqref{eqn:PoissonCover2Re} that
\[
\theta^*\cdot \Pr\left[\Pois\left(\int_{t=0}^1\alpha_t^idt\right)\leq K-1\right]\geq\sum_{k'=s_{i+1}+1}^K x_1^{k'}(\bm{\tau}^*) dt ~~~~ \forall i\in\{1,\dots, m\},
\]
and
\[
\theta^*\cdot \mathbb{E}\left[\min\left\{\Pois\left(\int_{t=0}^1\alpha_t^idt\right), K \right\}\right]\geq\sum_{k'=1}^{s_{i+1}}\int_{t=0}^1\beta^{i(k')}_tx^{k'}_t(\bm{\tau}^*)dt+\sum_{k'=s_{i+1}+1}^K\int_{t=0}^1\alpha^{i}_tx^{k'}_t(\bm{\tau}^*)dt~~ \forall 0\leq i\leq m-1.
\]
Moreover, the constraints \eqref{eqn:PoissonUB1Re}, \eqref{eqn:PoissonUB2Re}, and \eqref{eqn:PoissonUpdateRe} can be recovered from the dynamics described in \eqref{eqn:dualUB1}, \eqref{eqn:dualUB2}, \eqref{eqn:dualUB3}, and \eqref{eqn:dualUpdate}. Therefore, we know that we construct a feasible solution to $\PoisOPTRe_K(\bm{s})$ from a feasible solution to the problem $\max_{\bm{\tau}'}\mathrm{Dual}(\bm{\tau}', \bm{s}, \bm{G}')$, for fixed $\bm{G}'$ and $\bm{s}$. Since $\PoisOPTRe_K(\bm{s})$ takes the min objective, we know that 
\[
\PoisOPTRe_K(\bm{s})\leq \max_{ \bm{\tau}'}\mathrm{Dual}(\bm{\tau}', \bm{s}, \bm{G}')
\]
for every fixed $\bm{G}'$ and $\bm{s}$. 
Therefore, for every $\bm{s}$, we show that
\[
\PoisOPTRe_K(\bm{s})\leq \min_{\bm{G}'}\max_{ \bm{\tau}'}\mathrm{Dual}(\bm{\tau}', \bm{s}, \bm{G}')\leq  \max_{\bm{s}'} \min_{\bm{G}'}\max_{ \bm{\tau}'}\mathrm{Dual}(\bm{\tau}', \bm{s}', \bm{G}') \leq \min_{\bm{G}'}\max_{(\bm{\tau}', \bm{s}')\in\Gamma}\mathrm{Dual}(\bm{\tau}', \bm{s}', \bm{G}').
\]
Our proof is thus completed.

\subsection{Proof of \Cref{lem:OneinInterval}}
We use the well-known KKM lemma \citep{knaster1929beweis} to complete our proof. We first parameterize the threshold vector by its $m+1$ cell lengths:
\[
    z_i:=\tau_{i+1}-\tau_i,
    \qquad i=0,\ldots,m.
\]
Then, we have that
$z_i\ge0$ and
$\sum_{i=0}^m z_i=1$.
Thus, we know that the vector $\bm{z}=(z_0,\ldots,z_m)$ belongs to the simplex
\[
\Delta_m
:=
\left\{
\bm z\in\mathbb R_+^{m+1}:
\sum_{i=0}^m z_i=1
\right\}.
\]
Conversely, we can obtain that
$\tau_i(\bm z)
=
\sum_{j=0}^{i-1}z_j$.
For each $i=0,\ldots,m$, define the bottleneck set
\begin{equation}
    \mathcal{E}_i^\eps
    :=
    \left\{
    \bm z\in\Delta_m:
    c_i^\eps(\bm\tau(\bm z))
    =
    \min_{0\le j\le m}
    c_j^\eps(\bm\tau(\bm z))
    \right\}.
    \label{eqn:KKMSets}
\end{equation}
It is easy to see that
their union is $\Delta_m$. Moreover, since the function $\mathcal{R}_\eps$ is continuous,
it is direct to see that the function $c_i^\eps(\bm{\tau})$ is a continuous function for each $i$, which implies that the set $\mathcal{E}_i^\eps$ is a closed set for each $i$. 

We now verify the face-covering condition of the KKM lemma. For a set $S$ such that
$\varnothing\ne S\subseteq\{0,\ldots,m\}$, consider the face
\[
\Delta_S
:=
\left\{
\bm z\in\Delta_m:
z_j=0\text{ for all }j\notin S
\right\}.
\]
Take any $\bm z\in\Delta_S$, and let $j(\bm{z}, S)$ be a cell attaining the
smallest bottleneck such that
$c_{j(\bm{z}, S)}^\eps(\bm\tau(\bm z))
=
\min_{0\leq h\leq m}c_h^\eps(\bm\tau(\bm z))$. It is clear to see that
if $j(\bm{z}, S)\in S$, then $\bm{z}\in\mathcal{E}_{j(\bm{z}, S)}^\eps$.
Suppose otherwise $j(\bm{z}, S)\notin S$. Then from definition, we know that $z_{j(\bm{z}, S)}=0$, so cell $I_{j(\bm{z}, S)}$ is degenerate in that $\tau_{j(\bm{z}, S)+1}=\tau_{j(\bm{z}, S)}$.
The common endpoint must belong to
the closure of a nondegenerate cell $I_i$ with $i\in S$ ($i\notin S$ would imply $I_i$ is degenerate).  Since the
cellwise cover expressions coincide at the common finite endpoint, we must have that
\[
c_i^\eps(\bm\tau(\bm z))
\le
c_{j(\bm{z}, S)}^\eps(\bm\tau(\bm z)).
\]
Since from definition of $j(\bm{z}, S)$, we know that $c_{j(\bm{z}, S)}^\eps$ is the global minimum, equality holds in the above formula, and as a result, we have that
$\bm{z}\in\mathcal E_i^\eps$.

In summary, we have shown that
\[
\Delta_S
\subseteq
\bigcup_{i\in S}\mathcal E_i^\eps
\qquad
\forall\,\varnothing\ne S\subseteq\{0,\ldots,m\}.
\]
The KKM lemma therefore gives that
$\bigcap_{i=0}^m\mathcal E_i^\eps\ne\varnothing$.
Choose $\bm{z}^\eps$ in this intersection and set
$\bm\tau^\eps=\bm\tau(\bm z^\eps)$.  Then every cell bottleneck
equals the global bottleneck, proving
\eqref{eqn:allCellBottlenecksEqual}. Our proof is thus completed.

%At least one cell has positive length.  On every nondegenerate cell, the cellwise score is the actual normalized type-cover ratio. Degenerate endpoint scores are limits of the corresponding nondegenerate scores and cannot lie below the adjacent actual cell score. Hence the common value is exactly the global type-cover minimum, proving \eqref{eqn:thetaIsDualValue}. Our proof is thus completed.

\subsection{Proof of \Cref{lem:WorstAlpha}}
We denote by $\{\theta, \alpha^{i}_t, \beta^{i}_t, \forall i\in[m], \forall t\in[0,1]\}$ a feasible solution to $\PoisOPTRe_K(\bm{s})$ in \eqref{lp:PoissonDualRe}. We now make adjustment to the values of $\{\alpha^{i}_t, \forall i\in[m-1], \forall t\in[0,1]\}$, while keeping the values of the other variables in the solution fixed. 

Note that according to the update in \eqref{eqn:PoissonUB1Re}, \eqref{eqn:PoissonUB2Re}, and \eqref{eqn:PoissonUpdateRe}, the values of the variables $\{x^k_t, y^k_t,  \forall k\in[K], \forall t\in[0, 1]\}$ depend only on the variables $\{\beta^{i}_t, \forall i\in[m], \forall t\in[0,1]\}$, which would remain fixed as long as the variables $\{\beta^{i}_t, \forall i\in[m], \forall t\in[0,1]\}$ are fixed. Therefore, we can define a new set of solutions to $\PoisOPTRe_K(\bm{s})$ in \eqref{lp:PoissonDualRe}, which we denote by $\{\hat{\theta}, \hat{\alpha}^{i}_t, \hat{\beta}^{i}_t, \forall i\in[m], \forall t\in[0,1]\}$, and it holds that
\[
\hat{\theta}=\theta,~\hat{\beta}^{i}_t = \beta^{i}_t, \forall i\in[m], \forall t\in[0, 1].
\]
For each $i\in[m-1]$, we define the variable $\hat{\alpha}^{i}_t$ such that ($\hat{\alpha}^{m}_t=0$ for any $t$)
\[
\hat{\alpha}^{i}_t=\beta^{i+1}_t, \forall t\in[0, w_i]\text{~~and~~}\hat{\alpha}^{i}_t=\beta^{i}_t, \forall t\in(w_i, 1]
\]
and the value of $w_i$ is selected such that 
\[
\int_{t=0}^1 \alpha^{i}_t dt = \int_{t=0}^1 \hat{\alpha}^{i}_t dt.
\]
Since from the constraint that $\beta^{i+1}_t\leq\alpha^{i}_t\leq\beta^{i}_t$, we know that
\[
\int_{t=0}^1 \beta^{i+1}_t dt \leq \int_{t=0}^1 \alpha^{i}_t dt \leq \int_{t=0}^1 \beta^{i}_t dt, 
\]
which implies that we must be able to select a value $w_i\in[0, 1]$ such that
\[
\int_{t=0}^1 \hat{\alpha}^{i}_t dt=\int_{t=0}^{w_i} \beta^{i+1}_t dt + \int_{t=w_i}^1 \beta^{k}_t dt = \int_{t=0}^1 \alpha^{i}_t dt.
\]
Therefore, we know that the newly defined solution $\{\hat{\theta}, \hat{\alpha}^{i}_t, \hat{\beta}^{i}_t, \forall i\in[m], \forall t\in[0,1]\}$ will keep the values of the variables $\{x^k_t, y^k_t,  \forall k\in[K], \forall t\in[0, 1]\}$ fixed and the left-hand-side of the constraints \eqref{eqn:PoissonCover1Re} and \eqref{eqn:PoissonCover2Re} fixed. The right-hand-side of the constraint \eqref{eqn:PoissonCover1Re} also remains fixed.

It only remains to show that the right hand of the constraint \eqref{eqn:PoissonCover2Re} becomes smaller under the newly defined solution $\{\hat{\theta}, \hat{\alpha}^{i}_t, \hat{\beta}^{i}_t, \forall i\in[m], \forall t\in[0,1]\}$. We note that for any $0\leq t_1< t_2\leq 1$, for any $k\in[K]$, we have that
\[
\sum_{k'=k}^K x^{k'}_{t_1} \geq \sum_{k'=k}^K x^{k'}_{t_2}.
\] 
Then, for each $i\in[m-1]$, we have that
\[\begin{aligned}
&\sum_{k'=s_{i+1}+1}^K\int_{t=0}^1\alpha^{i}_tx^{k'}_tdt - \sum_{k'=s_{i+1}+1}^K\int_{t=0}^1\hat{\alpha}^{i}_tx^{k'}_tdt=\int_{t=0}^{w_i}(\alpha^{i}_t-\hat{\alpha}^{i}_t)\cdot\sum_{k'=s_{i+1}+1}^Kx^{k'}_tdt+\int_{t=w_i}^{1}(\alpha^{i}_t-\hat{\alpha}^{i}_t)\cdot\sum_{k'=s_{t+1}+1}^Kx^{k'}_tdt\\
&\geq (\sum_{k'=s_{i+1}+1}^Kx^{k'}_{w_i})\cdot\int_{t=0}^{w_i}(\alpha^{i}_t-\hat{\alpha}^{i}_t)dt+(\sum_{k'=s_{i+1}+1}^Kx^{k'}_{w_i})\cdot\int_{t=w_i}^{1}(\alpha^{i}_t-\hat{\alpha}^{i}_t)dt=(\sum_{k'=s_{i+1}+1}^Kx^{k'}_{w_i})\cdot\int_{t=0}^{1}(\alpha^{i}_t-\hat{\alpha}^{i}_t)dt=0
\end{aligned}\]
where the inequality follows from the fact that $\alpha^{i}_t-\hat{\alpha}^{i}_t\geq 0$ and $\sum_{k'=s_{i+1}+1}^Kx^{k'}_t\geq \sum_{k'=s_{i+1}+1}^Kx^{k'}_{s_k}$, when $t\leq w_i$, and $\alpha^{i}_t-\hat{\alpha}^{i}_t\leq 0$ and $\sum_{k'=s_{i+1}+1}^Kx^{k'}_t\leq \sum_{k'=s_{i+1}+1}^Kx^{k'}_{w_i}$, when $t> w_i$. In this way, we show that under the newly defined solution $\{\hat{\theta}, \hat{\alpha}^{i}_t, \hat{\beta}^{i}_t, \forall i\in[m], \forall t\in[0,1]\}$, the right-hand-side of the constraint \eqref{eqn:PoissonCover2Re} can only become smaller. Therefore, we know that $\{\hat{\theta}, \hat{\alpha}^{i}_t, \hat{\beta}^{i}_t, \forall i\in[m], \forall t\in[0,1]\}$ is still a feasible solution to $\PoisOPTRe_K(\bm{s})$ in \eqref{lp:PoissonDualRe} with the objective value $\hat{\theta}$ equals the optimal objective value of $\PoisOPTRe_K(\bm{s})$ in \eqref{lp:PoissonDualRe}. As a result, we know that $\{\hat{\theta}, \hat{\alpha}^{i}_t, \hat{\beta}^{i}_t, \forall i\in[m], \forall t\in[0,1]\}$ is a feasible solution to $\PoisOPTRe_K(\bm{s})$ in \eqref{lp:PoissonDualRe} with the same objective, which completes our proof.

\subsection{Proof of \Cref{lem:WorstBeta}}

For each $t\in[0, 1]$, we define a normalized variable 
\[
u^i_t = \frac{\beta^i_t}{\lambda_t}
\]
for each $i\in[m]$. Then, the feasible region for $\bm{u}_t=(u^1_t,\dots, u^m_t)$ can be given as
\begin{equation}\label{eqn:122201}
U = \left\{ (u^1,\dots,u^m): 0\leq u^m\leq \dots\leq u^1\leq 1 \right\}.
\end{equation}
An important fact we can obtain, which leads to our final conclusion, is that the extreme points of the set $U$ can be given as 
\begin{equation}\label{eqn:122202}
\mathrm{Ext}(U) = \left\{ (\underbrace{1, \dots,1}_{m'\text{~entries}},\underbrace{0,\dots,0)}_{m-m'\text{~entries}}: m' = 0, 1,\dots, m \right\}.
\end{equation}
The set of extreme points in \eqref{eqn:122202} follows from observation that any point in $U$ with some coordinate strictly between its neighbors admits a $\pm \epsilon$ perturbation preserving feasibility, hence is not extreme; the only points without such directions are the $0/m$ step vectors.

From now on, we use $u^i_t$ to represent the variable $\beta^i_t$ following the relationship $u^i_t = \frac{\beta^i_t}{\lambda_t}$.
Then, given any feasible solution $\{\theta, \alpha^{i}_t, u^{i}_t, \forall i\in[m], \forall t\in[0,1]\}$ to $\PoisOPTRe_K(\bm{s})$ in \eqref{lp:PoissonDualRe} satisfying the condition in \Cref{lem:WorstAlpha}, we construct another set of feasible solution $\{\hat{\theta}, \hat{\alpha}^{i}_t, \hat{u}^{i}_t, \forall i\in[m], \forall t\in[0,1]\}$ such that $\hat{\theta}\leq\theta+o(1)$ and $\hat{\bm{u}}_t\in\mathrm{Ext}(U)$ for each $t\in[0, 1]$. In order to construct such a $\{\hat{\theta}, \hat{\alpha}^{i}_t, \hat{u}^{i}_t, \forall i\in[m], \forall t\in[0,1]\}$, we partition the time interval $[0, 1]$ into $N+1$ evenly distributed discrete points with $N$ intervals given by $I_n=\left[\frac{n-1}{N}, \frac{n}{N}\right]$ for each $n\in[N]$. For each $i\in[m]$, we define the average quantity
\begin{equation}\label{eqn:122203}
    \bar{u}^i_n = \frac{\int_{t\in I_n} \beta^i_tdt}{\int_{t\in I_n}\lambda_t dt}.
\end{equation}
It is clear to see that for each $n\in[N]$, we have that $\bar{\bm{u}}_n\in U$, which implies that $\bar{\bm{u}}_n$ can be written as the convex combination of the extreme points of $U$, i.e., 
\begin{equation}\label{eqn:122204}
    \bar{\bm{u}}_n = \sum_{p=1}^P \pi^p_n\cdot \bm{u}_{\mathrm{Ext}}^p
\end{equation}
where $P=|\mathrm{Ext}(U)|$, $\mathrm{Ext}(U)=\left\{ \bm{u}_{\mathrm{Ext}}^p, p=1,\dots, P \right\}$, and it holds that $\sum_{p=1}^P\pi^p_n=1$ for each $n\in[N]$.

On each interval $I_n$, we define $\hat{\bm{u}}$ to cycle through the extreme points $\{ \bm{u}_{\mathrm{Ext}}^p, p=1,\dots, P \}$ in $\mathrm{Ext}(U)$, using duration $\frac{\pi^p_n}{N}$. Concretely, subdivide the interval $I_n$ into consecutive sub-intervals $I_{n, p}$ with length $|I_{n,p}| = \frac{\pi^p_n}{N}$ and set
\begin{equation}\label{eqn:122205}
    \hat{\bm{u}}_t = \bm{u}_{\mathrm{Ext}}^p, \text{~~for~}t\in I_{n,p}.
\end{equation}
We define $\hat{\bm{\alpha}}$ according to the conditions described in \Cref{lem:WorstAlpha} with respect to $\hat{\bm{u}}$, which guarantees that the switching time $\hat{w}_i$ for $\hat{\bm{\alpha}}$ is equivalent to the switching time $w_i$ for $\bm{\alpha}$, for each $i\in[m]$. 

By construction, we get exact preservation of the $\alpha$-integral:
\begin{equation}\label{eqn:122206}
    \int_{t=0}^1 \hat{\alpha}^{i}_t dt = \int_{t=0}^1 \alpha^{i}_t dt
\end{equation}
for each $i\in[m]$. Without loss of generality, we can assume that $w_i=\frac{k_i}{N}$ for some integer $k_i$, for each $i\in[m]$. Then, a even stronger condition we can get is that 
\begin{equation}\label{eqn:122207}
    \int_{t\in I_n} \hat{\alpha}^{i}_t dt= \int_{t\in I_n} \alpha^{i}_t dt
\end{equation}
for each $i\in[m]$ and $n\in[N]$.
We let $(\hat{\bx}, \hat{\by})$ be the dynamics driven by $(\hat{\bm{\alpha}}, \hat{\bm{\beta}})$ and let $(\bx, \by)$ be the dynamics driven by $(\bm{\alpha}, \bm{\beta})$. We write the dynamics as
\[
\frac{d\bm{x}}{dt} = F(\bm{x}, \bm{u})\text{~~and~~}\frac{d\hat{\bm{x}}}{dt} = F(\hat{\bm{x}}, \hat{\bm{u}}).
\]
As we can directly tell from the formulation, the function $F(\bx, \bm{u})$ is linear in $\bm{u}$ and Lipschitz continuous in $\bx$. On each interval $I_n$, we compare the true evolution under the original control $\bm{u}$ versus the chattering control $\hat{\bm{u}}$. Because the vector field for $F$ is bounded, a standard estimate gives a local error of order $O(1/N^2)$ per interval, and thus a global $O(1/N)$ error over $[0,1]$. Explicitly, it is easy to show that
\[
\|\bm{x} - \hat{\bm{x}}\|_{\infty}\leq C_1\cdot\frac{1}{N}
\]
for some constant $C_1$. Then, the right-hand-side of constraints \eqref{eqn:PoissonCover1Re} and \eqref{eqn:PoissonCover2Re} under the solution $(\hat{\bm{x}}, \hat{\bm{u}})$ deviate at most $O(1/N)$ from those under the solution $(\bm{x}, \bm{u})$. Therefore, from the feasibility of $\theta$ to the solutions $(\bm{x}, \bm{\alpha}, \bm{u})$ and from the conditions in \eqref{eqn:122206} and \eqref{eqn:122207}, we know that by setting 
\[
\hat{\theta} = \theta+O(1/N),
\]
we have that $(\hat{\theta}, \hat{\bm{x}}, \hat{\bm{\alpha}}, \hat{\bm{u}})$ forms a feasible solution to $\PoisOPTRe_K(\bm{s})$ in \eqref{lp:PoissonDualRe} while ensuring that $\hat{\bm{u}}_t\in\mathrm{Ext}(U)$ for each $t\in[0, 1]$ and $\hat{\theta} = \theta+O(1/N)$. Our proof is completed by letting $N\rightarrow\infty$.

\section{Missing Proofs of \Cref{sec:TightRatio}}

\subsection{Proof of \Cref{lem:uniformization}}
Fix any instance $\mathcal I=(\lambda,G',F)$ with $0\leq \lambda_t\leq M$.
We construct an equivalent constant-rate instance. Add a bottom type $\bot$
with value $F(\bot)=0$, lower than all genuine types, and define the new
type distribution by
\[
    \widetilde G'_t
    =
    \left(1-\frac{\lambda_t}{M}\right)\delta_{\bot}
    +
    \frac{\lambda_t}{M}G'_t .
\]
The new instance has calendar arrival rate $M$ and type distribution
$\widetilde G'_t$.

By the thinning property of Poisson processes, the non-bottom arrivals in the
new instance form exactly the same marked Poisson process as the arrivals in
the original instance. The remaining arrivals are bottom-type dummy customers
with value zero. Therefore the prophet's expected reward is unchanged, since
dummy arrivals add no value.

We now compare the best stock-based threshold policy in the two instances. Any
policy in the original instance can be used in the new instance by rejecting
the bottom type and applying the same thresholds to all genuine types, so the
optimal policy reward in the new instance is at least that in the original
instance. Conversely, in the new instance, accepting a bottom-type customer
earns zero reward and can only reduce future inventory. Hence there is an
optimal stock-based threshold policy that rejects all bottom-type customers.
After removing these rejected dummy arrivals, this policy induces a feasible
stock-based threshold policy in the original instance with the same reward.
Thus the optimal policy reward is also unchanged.

Hence both the prophet reward and the optimal stock-based policy reward are
the same in $\mathcal I$ and in its constant-rate representation. Therefore
the competitive ratio of every variable-rate instance can be matched by a
constant-rate instance with $\lambda_t\equiv M$. Since constant-rate
instances are a subclass of the original instance class, the two worst-case
competitive ratios are equal.

\subsection{Proof of \Cref{lem:BetaLemma}}

For notational simplicity, we write $s:=s_2$. We first introduce a normalized version of the $\alpha$ and $\beta$ control for convenience. Note that we already restrict the active time of $\beta^1$ control to be among $[0, S/M]$. Therefore, in the following derivations, we normalize the $\beta^1$ control to always equals $1$ among the interval $[0, S]$. Under such a normalization, we can define
\begin{equation}
\label{eq:m2sc_normalized_controls}
    D(t):=\frac{\beta_t^2}{\beta_t^1},
    \qquad
    c(t):=\frac{\alpha_t}{\beta_t^1}
\end{equation}
to be the normalized representation of the $\beta^1$ and $\alpha$ control. In this way, all the controls are defined on $[0, S]$ and satisfy the chaining constraints
\begin{equation}
\label{eq:m2sc_control_bounds}
    0\leq D(t)\leq c(t)\leq1
\end{equation}
for any $t$.
Moreover, under the constraint
$\Lambda_1:=\int_0^1\alpha_tdt$, 
we have the requirement that
\begin{equation}
\label{eq:m2sc_alpha_mass_active}
    \int_0^S c(t)dt=\Lambda_1.
\end{equation}
The state equations can be written in the following way under the normalized representation
\begin{equation}
\label{eq:m2sc_state}
    \dot x_t^k
    =-r_k(t)x_t^k+r_{k+1}(t)x_t^{k+1},
    \qquad k=1,\ldots,K,
\end{equation}
where we have 
\begin{equation}
\label{eq:m2sc_rates}
    r_k(t)=
    \begin{cases}
        D(t),&1\leq k\leq s,\\
        1,&s<k\leq K,
    \end{cases}
\end{equation}
with the convention $r_{K+1}x^{K+1}=0$, and with the same initial condition
$x_0^K=1$ and $x_0^k=0$ for $k<K$.  

Define the probability mass above the switching state and its boundary mass
by
\begin{equation}
\label{eq:m2sc_H_h}
    H(t):=\sum_{k=s+1}^Kx_t^k,
    \qquad
    h(t):=x_t^{s+1}.
\end{equation}
Summing \eqref{eq:m2sc_state} over $k=s+1,\ldots,K$ gives that
\begin{equation}
\label{eq:m2sc_H_derivative}
    H'(t)=-x_t^{s+1}.
\end{equation}
The upper states evolve independently of $(D,c)$.  In fact, since the $\beta^1$-control is active (equals $1$) in the entire interval $[0, S]$, it is direct to show that
\begin{equation}
\label{eq:m2sc_h_formula}
    x_t^{s+1}=e^{-t}\frac{t^{K-s-1}}{(K-s-1)!},
\end{equation}
so $h(t)>0$ for every $t>0$.
Now for fixed $(S,\Lambda_1)$, the four right-hand sides in
\eqref{eqn:m2PoissonCover1Re}, \eqref{eqn:m2PoissonCover2Re}, \eqref{eqn:m2PoissonCover3Re}, and \eqref{eqn:m2PoissonCover4Re} can be written as
\begin{subequations}
\label{eq:m2sc_Fs}
\begin{align}
    &F_1(D,c)
    :={\sum_{k=1}^Kx_S^k},
    &&F_2(D,c)
    :=H(S),
    \nonumber\\
    &F_3(D,c)
    :=\int_0^S
      \left[
        D(t)\sum_{k=1}^s x_t^k+c(t)H(t)
      \right]dt,
    &&F_4(D,c)
    :=\int_0^S
      \left[
        D(t)\sum_{k=1}^s x_t^k+H(t)
      \right]dt.
    \nonumber
\end{align}
\end{subequations}
Thus the fixed-$(S,\Lambda_1)$ problem is to minimize $\theta$ subject to the following constraints
\begin{equation}
\label{eq:m2sc_fixed_constraints}
\begin{aligned}
    &\theta\geq F_1(D,c),
    &&\theta\cdot\Pr[\Pois(\Lambda_1)\leq K-1]
        \geq F_2(D,c),\\
    &\theta\cdot\bE[\min\{\Pois(\Lambda_1),K\}]
        \geq F_3(D,c),
    &&\theta\cdot\bE[\min\{\Pois(M),K\}]
        \geq F_4(D,c),
\end{aligned}
\end{equation}
together with the constraints given in \eqref{eq:m2sc_state}, \eqref{eq:m2sc_control_bounds}, and
\eqref{eq:m2sc_alpha_mass_active}. It is without loss of generality to assume that
$0<S\leq M$ and $0<\Lambda_1<S$ since the boundary cases are handled easily.

We now prove that the fixed-$(S,\Lambda_1)$ problem has an exact optimum for
which $D$ is the indicator of one interval, which directly proves our argument that $\beta^2$ control is active on a single interval. We introduce the
regularization below is to remove flat switching ties and all perturbation terms
will be sent to zero at the end.
We
define the integral of the state $x_t^{s+1}$ to be
\begin{equation}
\label{eq:m2sc_U}
    U(t):=\int_0^t x_v^{s+1}dv,
    \qquad
    U_{\max}:=U(S).
\end{equation}
By \eqref{eq:m2sc_h_formula}, we know that the function $U$ is strictly increasing on $[0,S]$. To derive the regularization, 
we choose a countable dense set $\{u_n:n\geq1\}\subset(0,U_{\max})$ and define
\begin{equation}
\label{eq:m2sc_selector_H}
    H_\star(u):=\sum_{n=1}^\infty2^{-n}(u-u_n)^+,
    \qquad
    r_\star(u):=\sum_{n=1}^\infty2^{-n}\bI_{\{u>u_n\}}.
\end{equation}
It is clear to see that the series defining $H_\star$ converges uniformly. Each summand is convex,
and every interval contains at least one $u_n$ and therefore a
strict kink.  Hence, we know that the function $H_\star$ is strictly convex. At every point distinct
from the countable set $\{u_n\}$, termwise differentiation gives that
\begin{equation}
\label{eq:m2sc_selector_derivative}
    H_\star'(u)=r_\star(u).
\end{equation}
The function $r_\star$ is bounded and strictly increasing, but it changes
only through jumps.  Consequently, we know that $r_\star'(u)=0$ wherever its derivative
exists, and in particular almost everywhere.
Now we 
choose the parameters
\[
    \delta>0,
    \qquad
    \varepsilon>0,
    \qquad
    \zeta>0,
    \qquad
    0<\tau<\zeta/2,
\]
and define a term
\begin{equation}
\label{eq:m2sc_selector_psi}
    \psi(t):=-\zeta U(t)+\tau H_\star(U(t)).
\end{equation}
We keep the constraints \eqref{eq:m2sc_fixed_constraints} unchanged and
minimize the regularized objective given by
\begin{equation}
\label{eq:m2sc_regularizer}
    \theta
    +\delta\bigl[F_1+F_2+F_3+F_4\bigr]
    +\varepsilon\int_0^S D(t)\,dt
    +\int_0^S\psi(t)D(t)\,dt.
\end{equation}
The perturbation vanishes uniformly with
$(\delta,\varepsilon,\zeta,\tau)$.  Indeed, we have the bounded relationships that
$0\leq F_1,F_2\leq1$, $0\leq F_3,F_4\leq S\leq M$,
$0\leq U\leq1$, and $0\leq H_\star(U)\leq U\leq1$.  Therefore the
absolute value of the added terms is bounded by
\[
    \delta(2+2M)+\varepsilon M+M(\zeta+\tau),
\]
which tends to zero with the regularization parameters.

We can see that the pointwise control set
$\{(D,c):0\leq D\leq c\leq1\}$ is compact and convex. The state dynamics in \eqref{eq:m2sc_state}
are affine in $(D,c)$ and Lipschitz in $\bx$, while all terminal and running
expressions are continuous. The standard existence theorem (e.g., Theorem 4.1 of Chapter III in \cite{fleming2012deterministic}) for optimal
control problems with compact convex control sets therefore gives an optimal
ordinary control for the
fixed-$(S,\Lambda_1)$ problem (regularized or unregularized version).

We introduce the auxiliary states
\[
    \dot a(t)=c(t),
    \qquad
    \dot g_3(t)=D(t)\cdot\sum_{k=1}^s x_t^k+c(t)\cdot H(t),
    \qquad
    \dot g_4(t)=D(t)\cdot\sum_{k=1}^s x_t^k+H(t),
\]
with zero initial values.  Then the alpha-mass restriction is the endpoint
equality $a(S)=\Lambda_1$, while $F_3(D, c)=g_3(S)$ and $F_4(D, c)=g_4(S)$.
Thus all constraints are endpoint constraints, and the normal Pontryagin
minimum principle (e.g. \cite{vinter2000optimal}) applies.  The relevant constraint qualification is
quite direct here. The equality $a(S)=\Lambda_1$ has a nonzero
gradient in the $a$-coordinate, while none of the four covering inequalities
depends on that coordinate.  Moreover, the endpoint direction
$\Delta\theta=1$, with all other endpoint coordinates fixed, is tangent to
the equality and strictly relaxes every active covering inequality.  Hence
the endpoint constraint qualification holds, abnormal multipliers are
excluded, and the objective multiplier can be normalized to one.

Let the multipliers of the two terminal covering constraints be absorbed
into positive coefficients $\rho_1,\rho_2>0$, and let the multipliers of the
two running covering constraints be absorbed into positive coefficients
$\eta_0,\eta_1>0$. Their strict positivity follows from the
$\delta(F_1+F_2+F_3+F_4)$ term, which have been included in the definition of $\rho_1, \rho_2, \eta_0, \eta_1$.  Let $\kappa$ be the multiplier of
\eqref{eq:m2sc_alpha_mass_active}, and let $p_k$ be the costate of $x^k$,
with $p_0\equiv0$.  Up to constants independent of the controls, the regularized objective can he written as
\begin{equation}
\label{eq:m2sc_scalarized}
\begin{aligned}
    \mathcal J=&
    \rho_2F_1(D,c)+\rho_1F_2(D,c)
    +\eta_1F_3(D,c)+\eta_0F_4(D,c)+\kappa\left(\int_0^S c(t)dt-\Lambda_1\right)\\
    &+\varepsilon\int_0^S D(t)dt
    +\int_0^S\psi(t)D(t)dt.
\end{aligned}
\end{equation}
The terminal costate conditions are therefore can be written as
\begin{equation}
\label{eq:m2sc_terminal_costate}
    p_k(S)=
    \begin{cases}
        \rho_2,&1\leq k\leq s,\\
        \rho_1+\rho_2,&s<k\leq K.
    \end{cases}
\end{equation}
We further define
\begin{equation}
\label{eq:m2sc_mu}
    C:=\eta_0+\eta_1,
    \qquad
    \mu_k(t):=C+p_{k-1}(t)-p_k(t),
    \qquad 1\leq k\leq s.
\end{equation}
The Hamiltonian can be written as
\begin{equation}
\label{eq:m2sc_Hamiltonian}
\begin{aligned}
    \mathcal H(t)=&
    \sum_{k=1}^Kp_k(t)
       \bigl[-r_k(t)x_t^k+r_{k+1}(t)x_t^{k+1}\bigr]+C D(t)\sum_{k=1}^s x_t^k
      +\eta_1c(t)H(t)\\
      &+\eta_0H(t)
      +\kappa c(t)+(\varepsilon+\psi(t))D(t).
\end{aligned}
\end{equation}
For $1\leq k\leq s$, the coefficient of $x^k$ in the Hamiltonian is
$D(C+p_{k-1}-p_k)=D\mu_k$. Hence, we have that
\begin{equation}
\label{eq:m2sc_lower_costate}
    p_k'(t)=-D(t)\cdot\mu_k(t),
    \qquad 1\leq k\leq s.
\end{equation}
Using the definition in \eqref{eq:m2sc_mu}, we know that the condition \eqref{eq:m2sc_lower_costate} is equivalent to
\begin{equation}
\label{eq:m2sc_mu_dynamics}
    \mu_1'(t)=D(t)\cdot\mu_1(t),
    \qquad
    \mu_k'(t)=D(t)\cdot(\mu_k(t)-\mu_{k-1}(t)),
    \quad k=2,\ldots,s.
\end{equation}
Further note that the coefficient of $c$ in \eqref{eq:m2sc_Hamiltonian} is
\begin{equation}
\label{eq:m2sc_chi}
    \chi(t):=\kappa+\eta_1H(t).
\end{equation}
By \eqref{eq:m2sc_H_derivative}, we have that
\begin{equation}
\label{eq:m2sc_chi_derivative}
    \chi'(t)=-\eta_1\cdot x_t^{s+1}<0,
    \qquad t>0.
\end{equation}
For a fixed value of $D$, pointwise minimization over $D\leq c\leq1$
gives the following principle to determine the value of the $c$ control,
\begin{equation}
\label{eq:m2sc_c_rule}
    c(t)=D(t)+(1-D(t))\bI_{\{\chi(t)<0\}},
\end{equation}
except possibly at the unique zero of $\chi$. A direct observation is that the optional part of the
$\alpha$ control is placed as late as possible, which corresponds to the findings in \Cref{thm:WorstStructure}.

We then define the aggregated term
\begin{equation}
\label{eq:m2sc_Sigma}
    \Sigma(t):=\sum_{k=1}^s x_t^k\mu_k(t).
\end{equation}
The terms of the Hamiltonian that depend on $(D,c)$ can be represented as
\[
    D\Sigma+\chi c+(\varepsilon+\psi)D.
\]
For fixed $D$, we have that
\[
    \min_{D(t)\leq c(t)\leq1}\chi(t)\cdot c(t)
    =\chi^-(t)+D(t)\cdot\chi^+(t),
    ~~~\text{with the notation~~}
    z^+:=\max\{z,0\},\quad z^-:=\min\{z,0\}.
\]
Therefore, after minimizing over $c$, the difference between choosing
$D=1$ and choosing $D=0$ can be represented by
\begin{equation}
\label{eq:m2sc_gap}
    \Phi(t):=\Sigma(t)+\chi(t)^++\varepsilon+\psi(t).
\end{equation}
The Pontryagin minimum condition \citep{vinter2000optimal} yields that
\begin{equation}
\label{eq:m2sc_D_sign_rule_weak}
    \Phi(t)<0\Longrightarrow D(t)=1,
    \qquad
    \Phi(t)>0\Longrightarrow D(t)=0.
\end{equation}
We next compute the derivative of $\Sigma$.  The lower-state dynamics are
\[
    \dot x_t^k=D(t)\cdot(x_t^{k+1}-x_t^k),
    \quad 1\leq k<s,
    \qquad
    \dot x_t^s=h(t)-D(t)\cdot x_t^s.
\]
Combining these equations with \eqref{eq:m2sc_mu_dynamics} gives that
\begin{align}
    \Sigma'(t)
    =&(x_t^{s+1}-D(t)x_t^s)\cdot\mu_s(t)
      +\sum_{k=1}^{s-1}D(t)\cdot(x^{k+1}_t-x_t^k)\cdot\mu_k(t)
      +D(t)\cdot x_t^1\cdot\mu_1(t)\nonumber\\
      &+\sum_{k=2}^sD(t)\cdot x_t^k\cdot(\mu_k(t)-\mu_{k-1}(t))=x_t^{s+1}\cdot\mu_s(t).
    \label{eq:m2sc_Sigma_derivative}
\end{align}
Note that all internal terms cancel in the second equality. Since we have
$U'(t)=x_t^{s+1}$ and $H_\star'=r_\star$ almost everywhere, we also have
\[
    \psi'(t)=x_t^{s+1}\cdot\bigl[-\zeta+\tau r_\star(U(t))\bigr]
    \quad\text{almost entirely}.
\]
Together with \eqref{eq:m2sc_chi_derivative}, this yields that
\begin{equation}
\label{eq:m2sc_gap_derivative}
    \Phi'(t)
    =x_t^{s+1}\cdot\left[
        \mu_s(t)-\eta_1\bI_{\{\chi(t)>0\}}
        -\zeta+\tau r_\star(U(t))
    \right]
    \quad\text{almost entirely}.
\end{equation}
We further prove the following result, which removes the existence of fractional $D$ controls and active tie arcs.
\begin{claim}\label{claim:ExFrac}
It holds that 
\begin{equation}
\label{eq:m2sc_no_active_tie}
    \bigl|\{t\in[0,S]:D(t)>0,\ \Phi(t)=0\}\bigr|=0.
\end{equation}
\end{claim}
Note that the minimized Hamiltonian is affine in $D$.  Therefore, the conditions
\eqref{eq:m2sc_D_sign_rule_weak} and \eqref{eq:m2sc_no_active_tie} imply that
\begin{equation}
\label{eq:m2sc_D_strict_rule}
    D(t)=\bI_{\{\Phi(t)<0\}}.
\end{equation}
We now show that the strict negative set of $\Phi$ is connected. Consider a
nondegenerate open interval $I$ on which $D=0$. From
\eqref{eq:m2sc_lower_costate}, we know that every lower-block costate is constant on $I$.
In particular, we know that $\mu_s$ is constant on $I$.  Since the function $U$ is strictly
increasing, we may use $u=U(t)$ as the time coordinate on $I$.  Equations
\eqref{eq:m2sc_H_derivative} and \eqref{eq:m2sc_U} imply that
\[
    H(t)+U(t)=H(0),
\]
and as a result, we have that
\[
    \chi(t(u))=\kappa+\eta_1\bigl(H(0)-u\bigr),
\]
where $t(u)$ transfer the time coordinate on the interval $I$ to the original time indexed by $t\in[0, S]$.
Moreover, the condition \eqref{eq:m2sc_Sigma_derivative} gives that
$d\Sigma/du=\mu_s$ on $I$. Substituting these identities into
\eqref{eq:m2sc_gap}, we obtain, up to an additive constant independent of
$u$, that
\begin{equation}
\label{eq:m2sc_gap_convex_form}
    \Phi(t(u))
    =
    (\mu_s-\zeta)u
    +\bigl[\kappa+\eta_1(H(0)-u)\bigr]^+
    +\tau H_\star(u)
    +\text{constant}.
\end{equation}
The first term is affine, the second term is convex, and the third term is
strictly convex.  Hence, we know that the function $u\mapsto\Phi(t(u))$ is strictly convex on every
interval on which $D=0$.

Suppose that the set $\{\Phi<0\}$ could not be represented
by one interval. Then two consecutive components of $\{\Phi<0\}$ would
be separated by a nondegenerate open gap $I$.  (If two components were
separated only by a single point, their union would already agree almost
everywhere with one interval.) By \eqref{eq:m2sc_D_strict_rule}, we have that $D=0$ on
$I$, and by definition of the gap, $\Phi\geq0$ on $I$.  Since it is clear to see that $\Phi$ is
continuous, it equals zero at both endpoints of $I$.  Strict convexity from
\eqref{eq:m2sc_gap_convex_form} then implies that $\Phi<0$ at every
interior point of $I$, contradicting $\Phi\geq0$. In this way, we prove that the set $\{\Phi<0\}$ must be able to be represented
by one interval, which implies that the $D$ control is active on an interval.
Thus there exist
$0\leq a\leq b\leq S$ such that
\begin{equation}
\label{eq:m2sc_D_interval}
    D(t)=\bI_{[a,b]}(t)
    \quad\text{for almost every }t\in[0,S].
\end{equation}
Note that here, the interval is allowed to have zero length.

We have proved the desired result with regularization. We now remove the regularization and show that the active set of the $D$ control can still be represented by one interval.
Choose a sequence of regularization parameters converging to zero, with
$0<\tau_n<\zeta_n/2$, and let
$(\theta_n,D_n,c_n,\bx_n)$ be an optimum of the $n$-th regularized
fixed-$(S,\Lambda_1)$ problem that satisfies the desired property that $D_n$ is active on a single interval, i.e., there exists $0\leq a_n\leq b_n\leq S$ such that
\[
    D_n=\bI_{[a_n,b_n]}.
\]
The perturbation terms are uniformly
bounded in absolute value by numbers $e_n\downarrow0$.  Let
$(\theta^\star,D^\star,c^\star,\bx^\star)$ be an optimum of the
unregularized fixed-$(S,\Lambda_1)$ problem. Optimality of the regularized
solution gives that
\begin{align*}
    \theta_n
    &\leq
    \left[\text{regularized objective at }(\theta_n,D_n,c_n,\bx_n)\right]
    +e_n\leq
    \left[\text{regularized objective at }
      (\theta^\star,D^\star,c^\star,\bx^\star)\right]+e_n\\
    &\leq\theta^\star+2e_n.
\end{align*}
Passing to a subsequence, assume that
$a_n\to a$, $b_n\to b$, and $\theta_n\to\bar\theta$. Then, we have that
\[
    D_n\longrightarrow D:=\bI_{[a,b]}
    \quad\text{in }L^1[0,S].
\]
Since the coefficients in the linear state equations converge, the integral equations for $\bx_n-\bx$ gives that
\begin{equation}
\label{eq:m2sc_state_convergence}
    \sup_{0\leq t\leq S}\|\bx_n(t)-\bx(t)\|
    \longrightarrow0.
\end{equation}
Consequently, the terminal expressions $F_1,F_2$ and the running
expressions $F_3,F_4$ all converge.  Passing to the limit in
\eqref{eq:m2sc_fixed_constraints} shows that
$(\bar\theta,D,c,\bx)$ is feasible for the unregularized fixed problem.
The preceding objective bound gives $\bar\theta\leq\theta^\star$, while
optimality of $\theta^\star$ gives the reverse inequality.  Hence
$\bar\theta=\theta^\star$, and the limit control is an exact optimum with
$D=\bI_{[a,b]}$.

The boundary cases are simpler.  If $\Lambda_1=0$, then
$0\leq D\leq c$ and $\int_0^S c=0$ imply $D=c=0$.  If
$\Lambda_1=S$, then $c=1$ almost everywhere.  In this case one removes $c$
from the control variables; the same switching-function argument applies
with the term $\chi^+$ deleted.  The active-tie level becomes
$\zeta-\tau r_\star(U)>0$, and the switching gap on every $D=0$ interval is
still the sum of an affine function and the strictly convex term
$\tau H_\star$.  Thus $D$ is again the indicator of one interval.  The case
$S=0$ is immediate. In this way, we complete our proof of
\Cref{lem:BetaLemma}.

\subsubsection{Proof of \Cref{claim:ExFrac}}
Suppose, to the contrary, that this set has positive length. Since we have
$\{D>0\}=\bigcup_{n\geq1}\{D\geq1/n\}$, there are a number $d_0>0$ and a
positive-length subset $E$ on which $D\geq d_0$. By splitting $E$ if
necessary, we may also assume that the sign of $\chi$ is fixed on $E$.

An absolutely continuous function that is constant on the set $E$ has derivative
zero at almost every point of the set $E$. Hence we know that $\Phi'(t)=0$ at almost every
point of $E$. Since we have $h(t)>0$ for $t>0$, the condition
\eqref{eq:m2sc_gap_derivative} gives that
\begin{equation}
\label{eq:m2sc_flat_top_marginal}
    \mu_s(t)
    =\eta_1\bI_{\{\chi(t)>0\}}
     +\zeta-\tau r_\star(U(t))
    =:g(t)
    \quad\text{for a.e. }t\in E.
\end{equation}
Because from definition we have $0\leq r_\star\leq1$ and $\tau<\zeta/2$, we know that
\begin{equation}
\label{eq:m2sc_g_positive}
    g(t)\geq\zeta-\tau>\zeta/2>0.
\end{equation}
Moreover, we have that $g'(t)=0$ at almost every point of $E$.  Indeed, the indicator in
\eqref{eq:m2sc_flat_top_marginal} is constant on $E$.  Also, because $U$
is strictly increasing, for each $n$ there is at most one time $t_n$ such
that $U(t_n)=u_n$, and
\[
    r_\star(U(t))
    =\sum_{n=1}^\infty2^{-n}\bI_{\{t>t_n\}},
\]
with the terms for which $u_n>U(S)$ omitted. Thus, we know that $r_\star(U(t))$ is a
pure-jump function of $t$ and has derivative zero almost everywhere, which implies $g'(t)=0$ at almost every point of $E$.

We now use one elementary fact, repeatedly: if an absolutely continuous function
$f$ agrees with a function $g$ on a set $E$, and both derivatives exist at a
density point of $E$, then $f'=g'$ at that point. This follows by taking a
sequence of points of $E$ converging to the density point and comparing the
two difference quotients.

Applying this fact to \eqref{eq:m2sc_flat_top_marginal} gives that
$\mu_s'=g'=0$ almost everywhere on $E$. If $s=1$, then the first equation in
\eqref{eq:m2sc_mu_dynamics} gives that
\[
    0=\mu_1'=D\mu_1=Dg>0,
\]
which is a contradiction. Now suppose $s\geq2$. From the last equation in
\eqref{eq:m2sc_mu_dynamics} and $D\geq d_0$, we obtain
\[
    \mu_{s-1}=\mu_s=g
    \quad\text{a.e. on }E.
\]
Applying the same argument to $\mu_{s-1}=g$ gives that
$\mu_{s-1}'=0$, and hence $\mu_{s-2}=\mu_{s-1}=g$.  Repeating this step
through the lower block yields the following relationship
\[
    \mu_1=\mu_2=\cdots=\mu_s=g
    \quad\text{a.e. on }E.
\]
The first equation in \eqref{eq:m2sc_mu_dynamics} now gives again
\[
    0=\mu_1'=D\mu_1=Dg>0,
\]
which is impossible. This completes the proof of \eqref{eq:m2sc_no_active_tie}.

\subsection{Proof of \Cref{prop:m2_finite_M}}
Fix $\bm z=(S,a,b,\Lambda_1)\in\mathcal Z_{2,M}$. Under the change of
variables $u=Mt$, the controls and state equations are expressed on
the active-time interval $[0,S]$ with normalized rates in $\{0,1\}$.
Consequently, for fixed $\bm z$, the quantities
$\mathcal R_1(\bm z),\ldots,\mathcal R_4(\bm z)$ do not otherwise depend on $M$.

The limiting objective is therefore the same expression as
$\Theta_{K,s_2,M}(\bm z)$, except that its last term
$\mathcal R_4(\bm z)/H_K(M)$ is replaced by $\mathcal R_4(\bm z)/K$. Since
$H_K(M)\le K$ and $\mathcal R_4(\bm z)$ is the expected number of allocated
units, we have $0\le \mathcal R_4(\bm z)\le K$ and hence
\[
\begin{aligned}
0
\le
\Theta_{K,s_2,M}(\bm z)
-
\max\left\{
R_1(\bm z),
\frac{R_2(\bm z)}{P_K(\Lambda_1)},
\frac{R_3(\bm z)}{H_K(\Lambda_1)},
\frac{R_4(\bm z)}{K}
\right\}
\le
\mathcal R_4(\bm z)
\left(\frac{1}{H_K(M)}-\frac{1}{K}\right)
\le
\frac{K-H_K(M)}{H_K(M)}.
\end{aligned}
\]
Finally, it is direct to see that the value of $\PoisOPTRe_K(\bm{s}, M)$ decreases with $M$. Combining \Cref{lem:TightOblivious}, \Cref{thm:PoissonOPT}, \Cref{thm:m2-reduced}, and minimizing over $\mathcal Z_{2,M}$ proves the result.

\section{Missing Proofs of \Cref{sec:GeneralThresholdStructure}}
\label{app:finite_beta_structure_proofs}

\subsection{Proof of \Cref{lem:gen_safe_component}}
\label{app:proof_safe_component}

We first introduce the notation used throughout the proof. Recall that we work in active time, so we have
\[
1=D_1(t)\ge D_2(t)\ge\cdots\ge D_m(t)\ge0,
\qquad t\in[0,S],
\]
Here, we do not require at this stage that the controls are binary. We simply regard the $D$-control as an equivalent representation of the $\beta$ control. In the
original calendar-time formulation, we define
$\zeta(t):=\int_0^t\beta_s^1ds$.
On the set $\{\beta^1>0\}$, divide every $\beta$ and $\alpha$ intensity by
$\beta^1$ and use $d\zeta=\beta^1dt$. By the nesting constraints, all $\beta$ and $\alpha$
intensities vanish whenever $\beta^1=0$.  Therefore, the state equations,
all running expressions, all $\alpha$ masses, and all terminal expressions
are preserved by this change of variables; periods with $\beta^1=0$
contribute only idle time.  Conversely, any active-time control of length
$S\le\int_0^1\lambda_tdt$ can be lifted to calendar time by placing its
active clock on a set of total $\lambda$-mass $S$ and setting
$\beta^1=\lambda$ there and all $\beta$ controls to zero elsewhere.  Thus the
active-time formulation is exactly equivalent to the original fixed-$\bA$
problem and we use this convenient representation throughout the proof.

%The compatible bang--bang argument in \Cref{claim:compatible_bang_bang} below proves that an optimal solution of every strictly regularized active-time problem is cutoff-valued almost everywhere. After that claim, the controls may be interpreted in terms of the discrete modes used throughout the remainder of the proof.  

With the conventions $s_{m+1}:=0$ and $D_{m+1}(t):=0$, we define the normalized $\alpha$ controls
\begin{equation}
\label{eqn:c_control}
c_i(t):=\frac{\alpha_t^i}{\beta_t^1},
\qquad i=1,\dots,m-1.
\end{equation}
Then, we have that
\[
D_{i+1}(t)\le c_i(t)\le D_i(t),
\qquad
\int_0^S c_i(t)\,dt=A_i.
\]
For notational convenience, we set $c_0(t)\equiv1$.
For each $i=1,\dots,m$, define the upper-tail mass
\[
R_i(t):=\sum_{k=s_{i+1}+1}^{K}x_t^k.
\]
For $i=1,\dots,m-1$, mass exits this upper tail only through the boundary state $s_{i+1}+1$, and only while level $i$ is active. Hence, we have
\begin{equation}
\label{eq:tail_derivative_safe}
R_i'(t)=-D_i(t)x_t^{s_{i+1}+1}\le0.
\end{equation}
Thus each $R_i$ is nonincreasing. Once the $\beta$ controls and the masses $\bA$ are fixed, \Cref{lem:WorstAlpha} and \Cref{thm:WorstStructure} imply that the optional $\alpha$ mass is late-filled. In the Hamiltonian derivation below, we first treat each $c_i$ as an independent control subject to its pointwise bounds and mass constraint; the late-fill rule is recovered from the corresponding switching coefficient.
Each state $k$ belongs to the unique block
\[
B_i:=\{s_{i+1}+1,\dots,s_i\},
\qquad i=1,\dots,m.
\]
If $k\in B_i$, its active-time transition rate is
$r_k(t)=D_i(t)$.
The state equations are
\begin{equation}
\label{eq:active_state_dynamics_safe}
\dot x_t^k
=
-r_k(t)x_t^k+r_{k+1}(t)x_t^{k+1},
\qquad k=1,\dots,K,
\end{equation}
with
\[
r_{K+1}(t)x_t^{K+1}:=0,
\qquad
x_0^K=1,
\qquad
x_0^k=0\quad(k<K).
\]
For $i=0,\dots,m-1$, define the running expression
\begin{equation}
\label{eq:def_G_safe}
G_i(\bD,\bm c,\bx)
:=
\int_0^S
\left[
\sum_{k=1}^{s_{i+1}}r_k(t)x_t^k
+
\sum_{k=s_{i+1}+1}^{K}c_i(t)x_t^k
\right]dt,
\end{equation}
where the second sum is empty for $i=0$ because $s_1=K$. The terminal expressions are $R_i(S)$, $i=1,\dots,m$.

The original fixed-$\bA$ problem is generally nonconvex. We therefore use Pontryagin/KKT conditions only as first-order necessary conditions, not as a sufficient argument. To remove singular switching arcs, we first study a family of strictly regularized problems and later let all perturbations vanish.
If $A_i=0$ for some $i$, then from
$0=\int_0^S c_i(t)dt$ and
$D_{i+1}(t)\le c_i(t)$,
we have that
$c_i(t)=D_{i+1}(t)=0$ almost entirely. 
The nesting constraints then force every deeper level to be inactive. In that case we apply the argument to the remaining active levels and reinsert the deleted intervals with zero length. We may therefore assume below that all nonredundant Poisson denominators are positive.

%Let $\Theta(\bD,\bm c)$ denote the original fixed-$\bA$ objective.
We choose fixed constants
$a_r>0$ for $r=1,\dots,m$,
and
$b_r>0$ for
$r=0,\dots,m-1$.
Choose positive regularization parameters
\[
\delta>0,
\qquad
\varepsilon_2,\dots,\varepsilon_m>0,
\qquad
\zeta>0,
\qquad
0<\tau<\frac{\zeta}{2}.
\]
We consider constructing the following regularizer.
\begin{claim}
\label{claim:singular_slope}
For every $U_{\max}>0$, there exists a Lipschitz strictly convex function
$H_\star:[0,U_{\max}]\to\mathbb R$
and a bounded strictly increasing function
$r_\star:[0,U_{\max}]\to[0,1]$
such that
\[
H_\star'(u)=r_\star(u)
\quad\text{for a.e. }u,
\text{~~and~~}
r_\star'(u)=0
\quad\text{for a.e. }u.
\]
\end{claim}
Let $\bar{x}(t)$ denote the trajectory of the boundary state
$x_t^{s_2+1}$ generated by the always-active top block $B_1$. This
trajectory is independent of $D_2,\dots,D_m$. Define
\begin{equation}
\label{eq:def_U_selector}
U(t):=\int_0^t\bar x(s)\,ds,
\qquad
U_{\max}:=U(M).
\end{equation}
For every $t>0$, $\bar x(t)>0$, so $U$ is strictly increasing.  Choose
$r_\star$ and $H_\star$ from \Cref{claim:singular_slope} on
$[0,U_{\max}]$, and define
\begin{equation}
\label{eq:def_omega_selector}
\omega_{\zeta,\tau}(t)
:=-\zeta U(t)+\tau H_\star(U(t)).
\end{equation}
Then we know that
\begin{equation}
\label{eq:omega_derivative}
\omega_{\zeta,\tau}'(t)
=
\bar x(t)\bigl[-\zeta+\tau r_\star(U(t))\bigr].
\end{equation}
Because $0\le r_\star\le1$ and $\tau<\zeta/2$, we have that
\begin{equation}
\label{eq:selector_positive_level}
\zeta-\tau r_\star(U(t))\ge\zeta-\tau>\frac{\zeta}{2}>0.
\end{equation}
Moreover, from the properties of $r_\star$, it holds that
\begin{equation}
\label{eq:singular_composition_derivative}
\frac{d}{dt}r_\star(U(t))=0
\quad\text{for a.e. }t\in[0,M].
\end{equation}
Consider the regularized problem obtained by replacing the
objective $\theta$ by
\begin{equation}
\label{eq:strict_regularized_objective}
\begin{aligned}
\theta
+
\delta\sum_{r=1}^{m}a_rR_r(S)
+
\delta\sum_{r=0}^{m-1}b_rG_r(\bD,\bm c,\bx)
+
\sum_{i=2}^{m}\varepsilon_i\int_0^S D_i(t)dt
+
\int_0^M\omega_{\zeta,\tau}(t)D_2(t)dt,
\end{aligned}
\end{equation}
where all controls are extended by zero after their active horizon. Since all perturbations are uniformly
bounded, we have
\[
\sup_{\bD,\bm c}
\left|
\delta\sum_ra_rR_r(S)
+
\delta\sum_rb_rG_r
+
\sum_i\varepsilon_i\int D_i
+
\int\omega_{\zeta,\tau}D_2
\right|
\longrightarrow0
\]
as $\delta,\bm\varepsilon,\zeta,\tau\downarrow0$ with
$0<\tau<\zeta/2$.
The regularization contains three conceptually distinct parts.
First, the terms multiplied by $\delta$ assign a small positive weight
to every terminal and running covering expression. Their only purpose is
to make all effective coefficients (defined later) strictly positive.
Second, the penalties
$\varepsilon_i\int_0^S D_i(t)dt$
break zero-level ties between two adjacent modes. They are used to exclude
zero-valued singular switching arcs.
Third, the selector
$\int_0^M\omega_{\zeta,\tau}(t)D_2(t)dt$
is used only at the outer $\beta$ level. Its linear part
$-\zeta U(t)$ makes every outer active tie occur at a strictly positive
marginal level, which is further needed to prove the bang--bang argument. Its strictly
convex part $\tau H_\star(U(t))$ makes the outer switching gap strictly
convex on every mode-$1$ interval, which forces $D_2$ to be connected.
The slope of $H_\star$ has derivative zero almost everywhere, so this
second part does not disturb the later equations.

Let
$\varpi
:=
(\delta,\bm\varepsilon,\zeta,\tau)$
denote the vector of regularization parameters, and let
$(P_{\varpi})$ denote the corresponding regularized fixed-$\bA$
problem.
The argument has two stages. First, we fix $\varpi$ and prove that an
optimal solution of $(P_{\varpi})$ has the desired bang--bang and
finite-island structure. The resulting component bounds will not depend
on $\varpi$. Second, after these uniform bounds have been established,
we choose a sequence
\[
\varpi_n
=
(\delta_n,\bm\varepsilon_n,\zeta_n,\tau_n)
\longrightarrow0
\]
and pass to a limit of the corresponding structured optimal solutions.
Since the perturbations vanish uniformly, the limiting control is an
exact optimum of the original fixed-$\bA$ problem and retains the same
component bounds.

For every fixed $\varpi$, the relaxed regularized problem
$(P_{\varpi})$ attains its minimum. This follows from the standard
existence theorem for relaxed finite-dimensional optimal-control
problems: after extending the controls by zero to the common interval
$[0,M]$, the pointwise control set is compact and convex, the state
dynamics are affine in the controls and Lipschitz in the state, the
active horizon belongs to the compact interval $[0,M]$, and all
objective and constraint functionals are continuous under the standard
weak control and uniform state convergence.

Fix one optimal relaxed solution of $(P_{\varpi})$, and denote its
active horizon by $S$. 
Fix the active horizon $S$ of such an optimum.  Introduce auxiliary states
\[
\dot a_i(t)=c_i(t),~~~ a_i(0)=0,
~~~ i=1,\dots,m-1,
\text{~~and~~}
\dot g_r(t)
=
\sum_{k=1}^{s_{r+1}}r_k(t)x_t^k
+
\sum_{k=s_{r+1}+1}^{K}c_r(t)x_t^k,
~~~ r=0,\dots,m-1.
\]
Then the alpha-mass restrictions are the endpoint equalities
$a_i(S)=A_i$, and $G_r=g_r(S)$.  We apply the standard Pontryagin
minimum principle with endpoint equalities and inequalities; see, e.g.,
\citep{vinter2000optimal}. Every nonredundant covering inequality has the form
$F_j(\bD,\bm c,\bx)-\theta d_j\le0$ for $d_j>0$.
We can see that the endpoint constraint qualification holds, the Pontryagin
multiplier is normal, and the objective multiplier can be normalized to
one. 
Let
$\lambda_r\ge0$ for $r=1,\dots,m$
be the multipliers of the terminal covering constraints, and let
$\gamma_r\ge0$ for $r=0,\dots,m-1$
be the multipliers of the running covering constraints.  Define the
effective coefficients
$\rho_r:=\lambda_r+\delta a_r>0$ for
$r=1,\dots,m$ and 
$\eta_r:=\gamma_r+\delta b_r>0$ for $r=0,\dots,m-1$.
Let
$\kappa_i$ for $i=1,\dots,m-1$
be the multipliers of the fixed alpha-mass constraints.  After removing
terms independent of $(\bD,\bm c,\bx)$, the scalarized first-order objective
is
\begin{equation}
\label{eq:scalarized_obj_safe}
\begin{aligned}
\mathcal J_{\delta,\bm\varepsilon,\zeta,\tau}
=
\sum_{r=1}^{m}\rho_rR_r(S)
+
\sum_{r=0}^{m-1}\eta_rG_r(\bD,\bm c,\bx)
+
\sum_{i=1}^{m-1}\kappa_i
\left(
\int_0^S c_i(t)dt-A_i
\right)
+
\sum_{i=2}^{m}\varepsilon_i\int_0^S D_i(t)dt +
\int_0^M\omega_{\zeta,\tau}(t)D_2(t)dt.
\end{aligned}
\end{equation}
Let $p_k(t)$ be the costate associated with $x_t^k$, with convention
$p_0(t):=0$.
For block $B_i$, define
\begin{equation}
\label{eq:def_Ci}
C_i:=\eta_0+\eta_1+\cdots+\eta_{i-1}.
\end{equation}
If $k\in B_i$, the transition rate $r_k=D_i$ appears directly in the running expressions with total coefficient $C_i$. Consequently, the marginal value of activating the transition $k\to k-1$ is
\begin{equation}
\label{eq:def_mu}
\mu_k(t)
:=
C_i+p_{k-1}(t)-p_k(t),
\qquad k\in B_i.
\end{equation}
For $k\in B_i$, the coefficient of $x_t^k$ in the Hamiltonian is
$D_i(t)\mu_k(t)
+
\sum_{\ell=i}^{m-1}\eta_\ell c_\ell(t)$.
Therefore the costate equation is
\begin{equation}
\label{eq:corrected_costate}
p_k'(t)
=
-D_i(t)\mu_k(t)
-
\sum_{\ell=i}^{m-1}\eta_\ell c_\ell(t),
\qquad k\in B_i.
\end{equation}
The terminal conditions are
\begin{equation}
\label{eq:terminal_costate_safe}
p_k(S)=\sum_{r=i}^{m}\rho_r,
\qquad k\in B_i.
\end{equation}
For $i=1,\dots,m-1$, define the $\alpha$ switching coefficient
\begin{equation}
\label{eqn:def_chi}
\chi_i(t):=\kappa_i+\eta_iR_i(t).
\end{equation}
By \eqref{eq:tail_derivative_safe}, we have that
\begin{equation}
\label{eq:chi_derivative_safe}
\chi_i'(t)
=
-\eta_iD_i(t)x_t^{s_{i+1}+1}\le0.
\end{equation}
Thus $\chi_i$ is nonincreasing and has at most one zero cutoff. Pointwise minimization of $\chi_i(t)c_i(t)$ over
$D_{i+1}(t)\le c_i(t)\le D_i(t)$
gives the expression, away from the cutoff set $\{\chi_i=0\}$, that
\begin{equation}
\label{eq:pointwise_c_rule}
c_i(t)
=
D_{i+1}(t)
+
\bigl(D_i(t)-D_{i+1}(t)\bigr)
\bI_{\{\chi_i(t)<0\}}.
\end{equation}
Since
$\chi_i$ is nonincreasing, the rule in \eqref{eq:pointwise_c_rule} uses the lower feasible bound before
its cutoff and the upper feasible bound after its cutoff.  If $\chi_i=0$
on an interval, we choose the optional $\alpha$ mass (the $c$-control) as late as possible on
that interval. Following an argument similar to \Cref{lem:WorstAlpha}, we know that every regularized
problem admits an optimal representative whose $\alpha$ controls are
late-filled. We use such a representative in the compactness argument at
the end of the proof.

For each block, define the full block score
\begin{equation}
\label{eq:block_score}
S_i(t)
:=
\sum_{k=s_{i+1}+1}^{s_i}x_t^k\mu_k(t).
\end{equation}
For each cutoff mode $r\in\{1,\dots,m\}$, define its optimized Hamiltonian
score, up to terms common to all modes, by
\begin{equation}
\label{eq:mode_hamiltonian}
\begin{aligned}
\widehat{\mathcal H}_r(t)
=
\sum_{j=2}^{r}S_j(t)
+
\sum_{j=1}^{r-1}\chi_j(t)
+
\bI_{\{r<m\}}\chi_r(t)^-+
\sum_{j=2}^{r}\varepsilon_j
+
\omega_{\zeta,\tau}(t)\bI_{\{r\ge2\}},
\end{aligned}
\end{equation}
where
$u^-:=\min\{u,0\}$ and
$u^+:=\max\{u,0\}$.
Indeed, in cutoff mode $r$, we have $c_j=1$ for $j<r$, the optimized value of
$\chi_rc_r$ is $\chi_r^-$ when $r<m$, and all deeper alpha controls are
zero.
Define the adjacent mode gap
$\widehat A_j(t)
:=
\widehat{\mathcal H}_{j+1}(t)-\widehat{\mathcal H}_j(t)$.
Then we have
\begin{equation}
\label{eq:adjacent_gap}
\widehat A_j(t)
=
S_{j+1}(t)
+
\chi_j(t)^+
+
\chi_{j+1}(t)^-
+
\varepsilon_{j+1}
+
\omega_{\zeta,\tau}(t)\bI_{\{j=1\}},
\end{equation}
where $\chi_m^-:=0$. We further
define
\begin{equation}
\label{eq:Gamma_safe}
\Gamma_i(t)
:=
\mu_{s_{i+1}}(t)
-
\eta_i\bI_{\{\chi_i(t)>0\}},
~~~ i=1,\dots,m-1,
\text{~~and~~}
\vartheta_i(t):=x_t^{s_{i+1}+1}\Gamma_i(t),
~~~
\vartheta_m(t):=0,
\end{equation}
and obtain the following expression.

\begin{claim}
\label{claim:adjacent_gap_derivative}
Along every relaxed regularized control, we have that
\begin{equation}
\label{eq:adjacent_gap_derivative}
\widehat A_j'(t)
=
D_j(t)\vartheta_j(t)
-
D_{j+2}(t)\vartheta_{j+1}(t)
+
\omega_{\zeta,\tau}'(t)\bI_{\{j=1\}}
\end{equation}
for $j=1,\dots,m-1$, with $D_{m+1}:=0$.
\end{claim}
We also have the following result.
\begin{claim}
\label{claim:ordinary_time_plateau}
Fix a block $B_i$ and a set $E$ of positive measure on which it holds that
$D_i(t)\ge d_0>0$.
Let a differentiable function $g$ satisfy that $g'=0$ almost
everywhere on $E$, and suppose
\[
\mu_{s_i}(t)=g(t)
\quad\text{for a.e. }t\in E.
\]
Then, after discarding a null subset of $E$, we have that
$\mu_k(t)=g(t)$ for all $k\in B_i$, 
and, if $i<m$, it holds that
\begin{equation}
\label{eq:ordinary_plateau_boundary}
\frac{D_{i+1}(t)}{D_i(t)}\mu_{s_{i+1}}(t)
-
\eta_i\frac{c_i(t)}{D_i(t)}
=
g(t).
\end{equation}
If $i=m$, the same assumptions force $g(t)=0$ almost everywhere on
$E$.
\end{claim}
Note that \Cref{lem:WorstBeta} implies that the $\beta$-control can take binary values for the original problem. We can further show that for the regularized problem, every normal optimal solution also satisfies that $D_i(t)\in\{0,1\}$ for almost entirely $t$ and all $i=1,\dots,m$. Thus every regularized problem admits an optimal cutoff-valued
representative. We formalize the binary value of $D$ control in the following claim.
\begin{claim}\label{claim:compatible_bang_bang}
Every normal optimal solution of the regularized relaxed problem satisfies
\[
D_i(t)\in\{0,1\}
\qquad\text{for a.e. }t,
\qquad i=1,\dots,m.
\]
\end{claim}
%\end{comment}
For $i=1,\dots,m-1$, define the full switching gap between cutoff mode
$i$ and the deeper envelope by
\begin{equation}
\label{eq:full_switching_gap}
\widehat\Phi_i(t)
:=
\min_{r=i+1,\dots,m}
\left[
\widehat{\mathcal H}_r(t)-\widehat{\mathcal H}_i(t)
\right].
\end{equation}
For $r>i$, we let that
\[
\widehat\Psi_{i,r}
:=
\widehat{\mathcal H}_r-
\widehat{\mathcal H}_i
=
\sum_{j=i}^{r-1}\widehat A_j.
\]
Summing \eqref{eq:adjacent_gap_derivative} gives the exprssion that
\[
\widehat\Psi_{i,r}'
=
D_i\vartheta_i
+
\sum_{j=i+1}^{r-1}(D_j-D_{j+1})\vartheta_j
-
D_{r+1}\vartheta_r+
\omega_{\zeta,\tau}'\bI_{\{i=1\}}.
\]
Note that the actual $D$ control has a cutoff
mode. If that mode is $i$, every deeper $D_j$ is zero and every branch
has derivative $\vartheta_i$. If the actual mode is $q>i$, then the
branch $r=q$ is minimizing, due to the fact that
$D_i=\cdots=D_q=1$ and $D_{q+1}=0$,
all internal terms cancel and the same derivative remains. On a
positive-measure branch-tie set, tied absolutely continuous branches
have equal derivatives almost everywhere. With the definitions in \eqref{eq:Gamma_safe}, this proves that for $i=2,\dots,m-1$, in the active time of level $i$, we have
\begin{equation}
\label{eq:switch_safe}
\frac{d}{dL_i}\widehat\Phi_i(t)
=
x_t^{s_{i+1}+1}\Gamma_i(t)
\quad\text{almost entirely}
\end{equation}
Finally, substituting
\eqref{eq:omega_derivative}, for the outer gap, we have that
\begin{equation}
\label{eq:outer_switch_regularized}
\widehat\Phi_1'(t)
=
x_t^{s_2+1}
\left[
\Gamma_1(t)-\zeta+\tau r_\star(U(t))
\right]
\quad\text{almost entirely}
\end{equation}
%\begin{claim}[Derivative of the full switching gap]
%\label{claim:alpha_cutoff_contribution}
%For $i=2,\dots,m-1$, in the active time of level $i$,
%\begin{equation}
%\label{eq:switch_safe}
%\frac{d}{dL_i}\widehat\Phi_i(t)
%=
%x_t^{s_{i+1}+1}\Gamma_i(t)
%\quad\text{a.e.}
%\end{equation}
%For the outer gap,
%\begin{equation}
%\label{eq:outer_switch_regularized}
%\widehat\Phi_1'(t)
%=
%x_t^{s_2+1}
%\left[
%\Gamma_1(t)-\zeta+\tau r_\star(U(t))
%\right]
%\quad\text{a.e.}
%\end{equation}
%\end{claim}

We can further see that, 
%\begin{claim}[Positivity of the boundary mass]
%\label{claim:positive_boundary_mass}
for $i\ge2$, if the active time of level $i$ up to time $t$ is greater than 0, i.e., if $L_{i-1}(t)>0$, then we must have
\begin{equation}\label{eqn:16071703}
x_t^{s_i+1}>0,
\end{equation}
which follows directly from the state-transition dynamics that once a level becomes active, the probability mass on that level will always be positive.
%This holds for both relaxed and cutoff-valued controls.
%\end{claim}
%Indeed, %put $d:=K-(s_i+1)$. if $K-(s_i+1)=0$, then $x_t^K=\exp(-\int_0^tD_1)>0$. Suppose $K-(s_i+1)\ge1$.  Iterating the formula for the triangular pure-death system gives that $x_t^{s_i+1}$ as an integral over ordered transition times $0<u_{K-(s_i+1)}<\cdots<u_1<t$.  Every survival factor in that integral is strictly positive.  Whenever $D_{i-1}(u)>0$, nesting gives $r_k(u)\ge D_{i-1}(u)$ for every transition above $s_i$. Consequently, restricting the iterated integral to transition times in $\{u<t:D_{i-1}(u)>0\}$ gives the strictly positive lower bound
%\[
%\frac{c(t)}{(K-(s_i+1))!}
%\left(
%\int_0^tD_{i-1}(u)du
%\right)^{K-(s_i+1)}
%>0,
%\]
%where $c(t)>0$ is a lower bound for the finitely many survival factors on $[0,t]$.  Since the integral equals $L_{i-1}(t)>0$, the result follows.
%\end{proof}
We next eliminate active singular switching arcs in every level below $D_2$.
\begin{claim}
%[No active singular arcs below level $2$]
\label{claim:no_active_singular_arc}
For every regularized optimum and every $i=3,\dots,m$,
%\[
%\operatorname{Leb}_{L_{i-1}}
%\left(
%\{D_i=1\}\cap\{\widehat\Phi_{i-1}=0\}
%\right)
%=0.
%\]
%Consequently, 
after changing the control on a null set if necessary, we have that
\begin{equation}
\label{eq:strict_switch_rule_deep}
D_i=1
\quad\Longleftrightarrow\quad
\widehat\Phi_{i-1}<0
\quad\text{a.e. in }L_{i-1}\text{-time}.
\end{equation}
\end{claim}
We now derive the block representation used in the remaining proof. For every level, define the remaining future activity
\[
B_i(t):=\int_t^S D_i(s)ds.
\]
Deleting intervals on which $D_i=0$ and reversing the resulting clock does not change the number of connected components. We therefore write $z=B_i(t)\in[0,L_i]$ and use any right-continuous generalized inverse $t_i(z)$ of this collapsed reverse clock.
\begin{comment}
\begin{claim}[Block propagation]
\label{claim:block_propagation}
Let
\[
n_i:=s_i-s_{i+1}.
\]
In the remaining $L_i$-active clock,
\begin{equation}
\label{eq:correct_block_propagation}
\begin{aligned}
\mu_{s_i}(t_i(z))
={}&
C_i\Pi_{n_i-1}(z)
-
\rho_i\pi_{n_i-1}(z)
\\
&+
\int_0^z
\pi_{n_i-1}(z-y)W_i(y)\,dy,
\end{aligned}
\end{equation}
where
\[
\pi_\ell(u)
:=
e^{-u}\frac{u^\ell}{\ell!}\bI_{\{u\ge0\}},
\]
\[
\Pi_\ell(z)
:=
e^{-z}\sum_{r=0}^{\ell}\frac{z^r}{r!},
\]
and
\begin{equation}
\label{eq:def_Wi}
W_i(y)
:=
D_{i+1}(t_i(y))\mu_{s_{i+1}}(t_i(y))
-
\eta_i c_i(t_i(y)).
\end{equation}
For $i=m$, the forcing term is absent and is interpreted as zero.
\end{claim}
The bottom level is the homogeneous special case of \Cref{claim:block_propagation}.
\end{comment}
\begin{claim}
\label{claim:bottom_marginal_formula}
In the remaining $L_i$-active clock, it holds that
\begin{equation}
\label{eq:correct_block_propagation}
\mu_{s_i}(t_i(z))
=
C_i\Pi_{s_i-s_{i+1}-1}(z)
-
\rho_i\pi_{s_i-s_{i+1}-1}(z)+
\int_0^z
\pi_{s_i-s_{i+1}-1}(z-y)W_i(y)dy,
\end{equation}
where we define $\pi_\ell(u)
:=
e^{-u}\frac{u^\ell}{\ell!}\bI_{\{u\ge0\}}$, $\Pi_\ell(z)
:=
e^{-z}\sum_{r=0}^{\ell}\frac{z^r}{r!}$, and
\begin{equation}
\label{eq:def_Wi}
W_i(y)
:=
D_{i+1}(t_i(y))\mu_{s_{i+1}}(t_i(y))
-
\eta_i c_i(t_i(y)).
\end{equation}
Specifically, when $i=m$, it holds that
\begin{equation}
\label{eq:bottom_q_safe}
\mu_{s_m}(t)
=
C_m\Pi_{s_m-1}(B_m(t))
-
\rho_m\pi_{s_m-1}(B_m(t)),
\end{equation}
Moreover, as a function of $z=B_m(t)$, the right-hand side of \eqref{eq:bottom_q_safe} is monotone or valley-shaped.
\end{claim}
Regarding the expression given in \Cref{claim:bottom_marginal_formula}, we can obtain the following intermediate result.
%We replace the invalid moving-threshold argument by the following fixed-level convolution result.
\begin{claim}
\label{claim:variation_diminishing_step}
Let $W\in L^1(0,L)$ and, for $z\in[0,L]$, let
\begin{equation}
\label{eq:general_erlang_Q}
Q(z)
=
C\Pi_\ell(z)
-
\rho\pi_\ell(z)
+
\int_0^z\pi_\ell(z-y)W(y)\,dy,
\end{equation}
where $C,\rho\ge0$.  Fix $\zeta\ge0$.  Suppose that, for every
$a>\zeta$, the positivity set
$\{y\in(0,L): W(y)>a\}$
has at most $r$ connected components. Then every superlevel set
$\{z\in[0,L]:Q(z)>\zeta\}$
has at most $r+2$ connected components, for every real number $\zeta$.
%\[
%\#\operatorname{comp}\{W>a\}\le r.
%\]
%Then
%\[
%\#\operatorname{comp}\{Q>\zeta\}\le r+2.
%\]
\end{claim}
Finally, we show the following result.
\begin{claim}
\label{claim:no_repeated_birth}
Fix $i\in\{2,\dots,m-1\}$ and $a\ge0$. Then, in $L_i$-time, it holds that
\begin{equation}
\label{eq:Wi_superlevel_identity}
\{W_i>a\}
=
\{D_{i+1}=1\}
\cap
\{\mu_{s_{i+1}}>\eta_i+a\}.
\end{equation}
Moreover, inside every connected component of the set
$\{\mu_{s_{i+1}}>\eta_i+a\}$,
the set $\{W_i>a\}$ has at most one connected component.
\end{claim}
For $i=2,\dots,m$, let $P_i$ be a uniform upper bound on the number of connected components of every superlevel set
$\{\mu_{s_i}>\xi\}$, for any value $\xi$,
in the active time of level $i-1$ ($L_{i-1}$ time).
%let $P_i$ be a uniform upper bound such that, for every $\zeta\ge0$,
%\[
%\#\operatorname{comp}_{L_{i-1}}
%\{\mu_{s_i}>\zeta\}
%\le P_i.
%\]
At the bottom, \Cref{claim:bottom_marginal_formula} implies
\begin{equation}
\label{eq:Pm_bound}
P_m\le2.
\end{equation}
Now suppose $P_{i+1}$ is known. For every $a\ge0$, \Cref{claim:no_repeated_birth} gives that the number of connected components of the superlevel set $\{W_i>a\}$ is upper bounded by $P_{i+1}$.
%\[
%\#\operatorname{comp}\{W_i>a\}\le P_{i+1}.
%\]
Applying \Cref{claim:variation_diminishing_step} to \eqref{eq:correct_block_propagation} gives that
\[
P_i\le P_{i+1}+2.
\]
Backward induction therefore yields that
\begin{equation}
\label{eq:Pi_final}
P_i\le2(m-i+1),
\qquad i=2,\dots,m.
\end{equation}
%When $D_i=0$, all deeper controls and the corresponding alpha controls vanish, and \eqref{eq:corrected_costate} shows that $\mu_{s_i}$ is constant. Hence inserting or deleting the inactive intervals when passing between the $L_{i-1}$- and $L_i$-clocks cannot increase the number of strict superlevel components.

We now convert the superlevel bounds into beta-island bounds for the regularized problem. Fix $i\in\{3,\dots,m-1\}$. By \eqref{eq:switch_safe} and \eqref{eq:Gamma_safe}, we have that
\[
\frac{d}{dL_{i-1}}\widehat\Phi_{i-1}
=
x^{s_i+1}
\left[
\mu_{s_i}
-
\eta_{i-1}\bI_{\{\chi_{i-1}>0\}}
\right].
\]
Before the cutoff of $\chi_{i-1}$, the positive derivative set is
$\{\mu_{s_i}>\eta_{i-1}\}$,
and after the cutoff it is
$\{\mu_{s_i}>0\}$.
Each of these sets has at most $P_i$ components. Therefore, we know that the number of connected components of the set $\left\{
\frac{d}{dL_{i-1}}\widehat\Phi_{i-1}>0
\right\}$ is upper bounded by $2P_i$.
%\[
%\#\operatorname{comp}
%\left\{
%\frac{d}{dL_{i-1}}\widehat\Phi_{i-1}>0
%\right\}
%\le2P_i.
%\]
The complement of a union of at most $2P_i$ open intervals has at most $2P_i+1$ interval components. On each such complementary interval, we have that
\[
\frac{d}{dL_{i-1}}\widehat\Phi_{i-1}\le0
\quad\text{a.e.},
\]
so $\widehat\Phi_{i-1}$ is nonincreasing and can enter the strict negative half-line at most once. There may additionally be one negative component already present at the left endpoint. Hence, we know that
\[
J_i\le2P_i+2.
\]
Using \eqref{eq:Pi_final}, we obtain that
\begin{equation}
\label{eq:intermediate_J_bound}
J_i
\le
4(m-i+1)+2,
\qquad i=3,\dots,m-1.
\end{equation}
By \Cref{claim:no_active_singular_arc}, the active set agrees almost everywhere with the strict negative set of the full switching gap, so this is a bound on the actual $\beta$ islands.

For the bottom level, at the beginning of the $L_{m-1}$-clock the bottom block is empty. Hence
$S_m=0$,
and the regularized adjacent gap satisfies that
\[
\widehat\Phi_{m-1}(0)
=
\chi_{m-1}(0)^+
+
\varepsilon_m
>0.
\]
Before the cutoff of $\chi_{m-1}$, the bottom gap decreases only where
$\mu_{s_m}<\eta_{m-1}$,
and after the cutoff it decreases only where
$\mu_{s_m}<0$.
By \Cref{claim:bottom_marginal_formula}, $\mu_{s_m}$ is monotone or valley-shaped, so each of these two strict sublevel sets has at most one component. The bottom gap can therefore decrease on at most two intervals. Starting from a strictly positive value, each decreasing interval can create at most one strict negative component. By \Cref{claim:no_active_singular_arc}, we have that
\begin{equation}
\label{eq:bottom_J_bound}
J_m\le2.
\end{equation}
The preceding bounds hold uniformly for every regularized pure-mode optimum.
We carry out the vanishing-regularization passage jointly with the
connectedness argument for $D_2$ in the proof of
\Cref{lem:gen_D2_connected}.  The common endpoint-compactness argument at the
end of that proof produces one exact optimum satisfying all the bounds above. Our proof of \Cref{lem:gen_safe_component} is thus completed.

%We now let all the perturbations vanish. As formalized in the proof of \Cref{lem:gen_D2_connected}, such an operation preserves the bound on $J_i$ for $i=3,\dots, m$. Our proof of \Cref{lem:gen_safe_component} is thus completed.

\subsubsection{Prof of \Cref{claim:singular_slope}}
Choose a countable dense set (can be the set of all rational numbers)
$\{u_n:n\ge1\}\subset(0,U_{\max})$.
Then, we define
\begin{equation}
\label{eq:dense_kink_H}
H_\star(u)
:=
\sum_{n=1}^{\infty}
2^{-n}(u-u_n)^+,
\qquad
(x)^+:=\max\{x,0\},
\end{equation}
and define
\begin{equation}
\label{eq:dense_kink_r}
r_\star(u)
:=
\sum_{n=1}^{\infty}
2^{-n}\mathbf 1_{\{u>u_n\}}.
\end{equation}
The series defining $H_\star$ converges uniformly, because we have
$0\le 2^{-n}(u-u_n)^+\le 2^{-n}U_{\max}$.
Thus \(H_\star\) is continuous and Lipschitz. Each summand in
\eqref{eq:dense_kink_H} is convex, so \(H_\star\) is convex.
Moreover, every nondegenerate interval contains some $u_n$. Hence
$H_\star$ has a positive kink inside every nondegenerate interval and
therefore cannot be affine on any such interval. A convex function with
no nondegenerate affine interval is strictly convex.
At every point different from all $u_n$, termwise differentiation gives
$H_\star'(u)=r_\star(u)$.
Thus this identity holds almost everywhere.
Also, we have that
$0\le r_\star(u)\le\sum_{n=1}^{\infty}2^{-n}=1$.
If $u<v$ for arbitrary two points $u$ and $v$, density gives some $u_n\in(u,v)$, and therefore we have
\[
r_\star(v)-r_\star(u)\ge2^{-n}>0.
\]
Hence $r_\star$ is strictly increasing.
Finally, the derivative measure of $r_\star$ is
\[
dr_\star
=
\sum_{n=1}^{\infty}2^{-n}\delta_{u_n},
\]
which is purely atomic and therefore singular with respect to Lebesgue
measure. Consequently, we have
\[
r_\star'(u)=0
\quad\text{for a.e. }u,
\]
which completes our proof. 

\subsubsection{Proof of \Cref{claim:adjacent_gap_derivative}}
We first differentiate a block score.  For $j=2,\dots,m-1$, it holds that
\begin{equation}
\label{eq:block_score_derivative}
\begin{aligned}
S_j'(t)
={}&
D_{j-1}(t)x_t^{s_j+1}\mu_{s_j}(t)
+
x_t^{s_{j+1}+1}
\left[
\eta_jc_j(t)
-
D_{j+1}(t)\mu_{s_{j+1}}(t)
\right].
\end{aligned}
\end{equation}
For $j=m$, the same formula holds with the second line omitted.
For an interior state $k>s_{j+1}+1$ of $B_j$, subtracting the two
costate equations in \eqref{eq:corrected_costate} gives that
\[
\mu_k'
=
D_j(\mu_k-\mu_{k-1}).
\]
At the lower boundary $k=s_{j+1}+1$, subtraction across the two adjacent
blocks gives
\[
\mu_{s_{j+1}+1}'
=
D_j\mu_{s_{j+1}+1}
-
D_{j+1}\mu_{s_{j+1}}
+
\eta_jc_j.
\]
Differentiating $S_j=\sum_{k\in B_j}x^k\mu_k$ and substituting the state
and marginal equations cancels every internal term, leaving precisely
\eqref{eq:block_score_derivative}.
Next, by \eqref{eq:chi_derivative_safe}, we have
\[
(\chi_j^+)'
=
-\eta_jD_jx^{s_{j+1}+1}\bI_{\{\chi_j>0\}},
\text{~~and~~}
(\chi_{j+1}^-)'
=
-\eta_{j+1}D_{j+1}x^{s_{j+2}+1}
\bI_{\{\chi_{j+1}<0\}}.
\]
Using \eqref{eq:pointwise_c_rule}, it holds that
\[
\eta_{j+1}c_{j+1}
-
\eta_{j+1}D_{j+1}\bI_{\{\chi_{j+1}<0\}}
=
\eta_{j+1}D_{j+2}\bI_{\{\chi_{j+1}>0\}}.
\]
For $j=m-1$, all terms involving $\chi_m$, $c_m$, and $\eta_m$ are
absent, and the same calculation gives the formula with
$D_{m+1}\vartheta_m=0$.  Differentiating \eqref{eq:adjacent_gap} and
using these identities yields \eqref{eq:adjacent_gap_derivative}.

\subsubsection{Proof of \Cref{claim:ordinary_time_plateau}}
For $j=1,\dots,n_i$ where $n_i=s_i-s_{i+1}$, set
\[
M_j(t):=\mu_{s_{i+1}+j}(t).
\]
The ordinary-time marginal equations obtained from
\eqref{eq:corrected_costate} are
\begin{equation}
\label{eqn:16071701}
\dot{M}_j
=
D_i(M_j-M_{j-1}),
\qquad j=2,\dots,n_i,
\end{equation}
and, when $i<m$, we have that
\begin{equation}
\label{eqn:16071702}
\dot{M}_1
=
D_iM_1-D_{i+1}\mu_{s_{i+1}}+\eta_i c_i.
\end{equation}
For $i=m$, the first equation is instead $\dot{M}_1=D_mM_1$.

Since $M_{n_i}=g$ on $E$, it is easy to see that
$\dot{M}_{n_i}=g'=0$ almost everywhere on $E$. Because we also have the condition that $D_i\ge d_0>0$, from \eqref{eqn:16071701},
we obtain $M_{n_i-1}=M_{n_i}=g$ there. Restricting to the full
set on which this equality holds and repeating the argument yields that
\[
M_1=\cdots=M_{n_i}=g
\quad\text{a.e. on }E.
\]
Substitution into \eqref{eqn:16071702} gives
\eqref{eq:ordinary_plateau_boundary}.  In the bottom block it gives
$0=D_mg$. Since $D_m\ge d_0$, this forces $g=0$. Our proof is thus completed.

\subsubsection{Proof of \Cref{claim:compatible_bang_bang}}
For a feasible solution, denote
\[
\nu_r(t):=D_r(t)-D_{r+1}(t),
\qquad r=1,\dots,m,
\qquad D_{m+1}:=0.
\]
Then we have $\nu_r(t)\ge0$, $\sum_r\nu_r(t)=D_1(t)=1$, and $D_i(t)=\sum_{r=i}^{m}\nu_r(t)$. %After pointwise minimization over the $\alpha$ controls, linearity of the Hamiltonian gives that
%\begin{equation}
%\label{eq:mode_convex_combination}
%\widehat{\mathcal H}(D,c^*)
%=
%\sum_{r=1}^{m}\lambda_r\widehat{\mathcal H}_r.
%\end{equation}
Indeed, after pointwise minimization over the alpha
controls, \eqref{eq:pointwise_c_rule} gives that
\[
c_i^*(t)
=
D_{i+1}(t)+(D_i(t)-D_{i+1}(t))\cdot\bI_{\{\chi_i(t)<0\}}
=
\sum_{r>i}\nu_r(t)+\nu_i(t)\cdot\bI_{\{\chi_i(t)<0\}},
\]
which is the same convex combination of the $\alpha$ choices in the cutoff modes. 
When $\chi_i(t)=0$, the Hamiltonian is independent of $c_i(t)$, so the identity remains valid. The Pontryagin minimum condition \citep{vinter2000optimal} therefore implies that
\begin{equation}
\label{eq:supported_modes_minimize}
\nu_r(t)>0
\quad\Longrightarrow\quad
\widehat{\mathcal H}_r(t)=\min_q\widehat{\mathcal H}_q(t)
\quad\text{a.e.}
\end{equation}
Suppose that $D$ is nonbinary on a set of positive measure. Since there
are only finitely many mode supports and sign patterns of the functions
$\chi_i(t)$, there is a positive-measure set $E$ on which both are fixed.
Let $a<b$ be the two largest supported modes. By restricting to one of
the sets $\{\nu_a\ge1/n,\nu_b\ge1/n\}$, we may also assume
\[
\nu_a,\nu_b\ge d_0>0
\quad\text{on }E.
\]
Since no supported mode lies strictly between $a$ and $b$, we have that
\begin{equation}
\label{eq:support_block_equal_D}
D_{a+1}=\cdots=D_b=:d\ge d_0,
\qquad
D_{b+1}=0
\quad\text{on }E.
\end{equation}
Modes $a$ and $b$ are tied on $E$. Therefore, differentiating their gap yields value $0$ and
summing \eqref{eq:adjacent_gap_derivative} from $j=a$ to $b-1$ gives
\begin{equation}
\label{eq:bangbang_tied_derivative}
0
=
D_a\vartheta_a
+
\omega_{\zeta,\tau}'\bI_{\{a=1\}}
\quad\text{a.e. on }E.
\end{equation}
If $a=1$, then $D_a=1$, and the condition in
\eqref{eq:omega_derivative} and
\eqref{eq:bangbang_tied_derivative} give that
\begin{equation}
\label{eq:bangbang_positive_outer_level}
\mu_{s_2}
=
\eta_1\bI_{\{\chi_1>0\}}
+
\zeta-\tau r_\star(U)
=:
g.
\end{equation}
By \eqref{eq:selector_positive_level}, we have $g>\zeta/2>0$.  On the selected
sign component of $E$, the condition in \eqref{eq:singular_composition_derivative} establishes that
$g'=0$ almost everywhere.
If $a\ge2$, then $D_a>0$, and we have that
\begin{equation}
\label{eq:bangbang_inner_level}
\mu_{s_{a+1}}
=
\eta_a\bI_{\{\chi_a>0\}}
=:
g,
\end{equation}
where $g$ is either $0$ or the positive constant $\eta_a$.

Apply \Cref{claim:ordinary_time_plateau} first to block $B_{a+1}$.
For every intermediate supported block $i=a+1,\dots,b-1$,
from \eqref{eq:support_block_equal_D}, we have that
\[
D_i=D_{i+1}=c_i=d.
\]
Thus, the condition in \eqref{eq:ordinary_plateau_boundary} establishes that
\begin{equation}
\label{eq:positive_plateau_increment}
\mu_{s_{i+1}}=g_i+\eta_i,
\end{equation}
where $g_i$ is the propagated top marginal in block $B_i$.  The new
comparison function again has derivative zero almost everywhere on $E$.
Consequently, after the first intermediate block the propagated marginal
is strictly positive, and each further block adds the positive number
$\eta_i$.

At the last supported block $B_b$, we have $D_b=d$ and $D_{b+1}=0$.  If its
propagated top marginal is strictly positive, then
\eqref{eq:ordinary_plateau_boundary} gives that
$-\eta_b\frac{c_b}{d}>0$
when $b<m$, while the bottom-block conclusion of
\Cref{claim:ordinary_time_plateau} gives a positive quantity equal to
zero when $b=m$.  Both are impossible.

The only remaining case is
\[
a\ge2,
\qquad b=a+1,
\qquad g=0.
\]
\Cref{claim:ordinary_time_plateau} gives that every marginal in $B_b$ is zero.  If
$b<m$, it also gives $c_b=0$, and pointwise $\alpha$ minimization yields
$\chi_b\ge0$; when $b=m$, the $\chi_b^-$ term is absent.  Moreover,
$g=0$ and $\eta_a>0$ imply $\chi_a\le0$.  Therefore
\[
\widehat{\mathcal H}_b-
\widehat{\mathcal H}_a
=
S_b+\chi_a^++\bI_{\{b<m\}}\chi_b^-+\varepsilon_b
=
\varepsilon_b>0,
\]
contradicting the tie.  Hence no fractional minimizing face has positive
measure. Our proof is thus completed.

\subsubsection{Proof of \Cref{claim:no_active_singular_arc}}
Suppose instead that a positive-measure set in the $L_{i-1}$-clock
satisfies $D_i=1$ and $\widehat{\Phi}_{i-1}=0$.  Pull it back to ordinary
active time and denote the pullback by $E$.  Since $D_i=1$ implies
$D_{i-1}=1$, the ordinary-time and $L_{i-1}$ measures agree on $E$, so
$E$ has positive measure. From $\widehat{\Phi}_{i-1}=0$, we know that $\frac{d}{dL_{i-1}}\widehat{\Phi}_{i-1}=0$ almost entirely on $E$. Then, from \eqref{eq:switch_safe} and \eqref{eqn:16071703}, we know that $\Gamma_i(t)=0$ almost entirely on $E$. Thus, from the expression \eqref{eq:Gamma_safe}, we know that
\[
\mu_{s_i}
=
\eta_{i-1}\bI_{\{\chi_{i-1}>0\}}
\quad\text{a.e. on }E.
\]
Split $E$ at the unique cutoff of $\chi_{i-1}$.  On a positive-measure
subset, we have that
\begin{equation}
\label{eq:singular_level_xi}
\mu_{s_i}=\xi,
\qquad
\xi\in\{0,\eta_{i-1}\}.
\end{equation}
Apply \Cref{claim:ordinary_time_plateau} to block $B_i$, with
$D_i=1$ and comparison function $g\equiv\xi$.

If $\xi>0$, the condition that requires $g=0$ for $i=m$ in \Cref{claim:ordinary_time_plateau} rules out the case
$i=m$. If $i<m$, from \eqref{eq:ordinary_plateau_boundary} and the conditions that $D_i=1$, $g\equiv\xi$, we have that
\begin{equation}\label{eq:boundary_identify}
D_{i+1}(t)\mu_{s_{i+1}}(t) - \eta_i c_i(t) = \xi
\end{equation}
which shows that $D_{i+1}$ cannot be
zero. Hence we have $D_{i+1}=c_i=1$ and
$\mu_{s_{i+1}}=\xi+\eta_i>0$.
Applying the same argument successively to all lower active blocks
reaches the bottom block ($i=m$) with a strictly positive value for the corresponding $\mu$, which is
impossible by the condition that requires $g=0$ for $i=m$ in \Cref{claim:ordinary_time_plateau}.

Now suppose $\xi=0$. If $D_{i+1}=1$ on a positive-measure subset, then we have
$c_i=1$ and \eqref{eq:boundary_identify} gives
$\mu_{s_{i+1}}=\eta_i>0$, where $\eta_i>0$ comes from the definition that $\eta_i=\gamma_i+\delta b_i$ with $\gamma_i\geq0$, $\delta>0$, and $b_i>0$ from regularization. This would reduce to the preceding case that results in the bottom block ($i=m$) with a strictly positive value for the corresponding $\mu$, which is
impossible by the condition that requires $g=0$ for $i=m$ in \Cref{claim:ordinary_time_plateau}. Therefore
$D_{i+1}=0$ almost everywhere on the selected set. Since $\eta_i>0$,
\eqref{eq:boundary_identify} gives that $c_i=0$.  Pointwise $\alpha$ minimization \eqref{eq:pointwise_c_rule} then
implies $\chi_i\ge0$, while \eqref{eq:singular_level_xi} and
$\eta_{i-1}>0$ imply $\chi_{i-1}\le0$.  All marginals in $B_i$ are
zero, so $S_i=0$.  Since the actual mode is exactly $i$, from \eqref{eq:adjacent_gap}, we have that
\[
\widehat A_{i-1}
=
S_i+\chi_{i-1}^++\bI_{\{i<m\}}\chi_i^-+\varepsilon_i
=
\varepsilon_i>0,
\]
contradicting Hamiltonian minimization that forces $\widehat A_{i-1}\leq 0$.  Both cases are impossible.

Finally, in $L_{i-1}$-time, $\widehat\Phi_{i-1}<0$ forces the minimizing
mode to be at least $i$, hence $D_i=1$.  The reverse implication holds
outside the null active-tie set just excluded.  This proves
\eqref{eq:strict_switch_rule_deep}.

\subsubsection{Proof of \Cref{claim:bottom_marginal_formula}}
For $j=1,\dots,s_i-s_{i+1}$, define
\[
M_j(z):=\mu_{s_{i+1}+j}(t_i(z)).
\]
In the reverse block clock, the formula in \eqref{eq:corrected_costate} gives that
\[
M_1'(z)=W_i(z)-M_1(z),
\text{~~and~~}
M_j'(z)=M_{j-1}(z)-M_j(z),
\qquad j=2,\dots,s_i-s_{i+1}.
\]
From \eqref{eq:terminal_costate_safe}, we know that
\[
M_1(0)=C_i-\rho_i,
\qquad
M_j(0)=C_i
\qquad j=2,\dots,s_i-s_{i+1}.
\]
Solving this linear equation gives that
\[
M_{s_i-s_{i+1}}(z)
=
C_i\Pi_{s_i-s_{i+1}-1}(z)
-
\rho_i\pi_{s_i-s_{i+1}-1}(z)
+
\int_0^z\pi_{s_i-s_{i+1}-1}(z-y)W_i(y)dy.
\]
Since $M_{s_i-s_{i+1}}=\mu_{s_i}$, we know that \eqref{eq:correct_block_propagation} holds.

The formula \eqref{eq:bottom_q_safe} follows from \eqref{eq:correct_block_propagation} with $i=m$. Define the function
\[
F_m(z)
:=
C_m\Pi_{s_m-1}(z)-\rho_m\pi_{s_m-1}(z).
\]
If $s_m=1$, then we have
\[
F_m(z)=(C_m-\rho_m)e^{-z},
\]
which is monotone. Suppose $s_m\ge2$. Using the relationship
\[
\frac{d}{dz}\Pi_{s_m-1}(z)
=
-\pi_{s_m-1}(z)
\text{~~and~~}
\frac{d}{dz}\pi_{s_m-1}(z)
=
\pi_{s_m-2}(z)-\pi_{s_m-1}(z),
\]
we obtain that
\begin{equation}
\label{eq:bottom_derivative_shape}
F_m'(z)
=
(\rho_m-C_m)\pi_{s_m-1}(z)
-
\rho_m\pi_{s_m-2}(z).
\end{equation}
For $z>0$, it is direct to verify that
\[
\frac{\pi_{s_m-1}(z)}{\pi_{s_m-2}(z)}
=
\frac{z}{s_m-1}.
\]
Thus $F_m'$ has at most one zero. Since its sign near zero is nonpositive, $F_m$ is either decreasing or decreases first and then increases. Hence it is monotone or valley-shaped. Our proof is thus completed.

\subsubsection{Proof of \Cref{claim:variation_diminishing_step}}
The proof has two steps. First, we bound the number of sign changes of
$Q-a$ for every $a>\zeta$. We then convert this sign-change bound
into a component bound at the original level $\zeta$.

Fix $a>\zeta$. For a real-valued function $v$, let
$\operatorname{sc}(v)$
denote its number of strict sign changes: namely, the largest possible
number $N-1$ for which there exist
\[
y_1<y_2<\cdots<y_N
\]
such that every $v(y_j)$ is nonzero and consecutive selected values
have opposite signs.
The expression $Q-a$ can be viewed as an Erlang convolution with
three pieces of input: a constant input $C-a$ before time $0$; a negative mass $-\rho$ at time $0$; the input $W-a$ after time $0$.
For an expression using ordinary integrable functions, let $T>1$ and
$0<\epsilon<1$, and define
\begin{equation}
\label{eq:erlang_input_approx}
v_{T,\epsilon}(y)
:=
(C-a)\bI_{[-T,-\epsilon)}(y)
+
\left(C-a-\frac{\rho}{\epsilon}\right)
\bI_{[-\epsilon,0)}(y)
+
(W(y)-a)\bI_{(0,L)}(y).
\end{equation}
By assumption, the positive set of $W-a$ has at most $r$
components. Therefore $W-a$ has at most $2r$ strict sign changes
on $(0,L)$. The two constant pieces on the negative half-line can
add at most one sign change between them, and the transition across
\(0\) can add at most one more. Hence we have that
\begin{equation}
\label{eq:input_sign_change_bound}
\operatorname{sc}(v_{T,\epsilon,a})
\le
2r+2.
\end{equation}
The Erlang kernel $\pi_\ell$ is the $(\ell+1)$-fold convolution of
the one-sided exponential kernel. By Schoenberg's variation-diminishing theorem \citep{schoenberg1948variation},
convolution with $\pi_\ell$ cannot increase the number of strict sign
changes and we have that
\begin{equation}
\label{eq:erlang_variation_diminishing}
\operatorname{sc}
\bigl(
\pi_\ell*v_{T,\epsilon,a}
\bigr)
\le
2r+2.
\end{equation}

We now identify the limit of this convolution. For $z\in[0,L]$, we have
\begin{equation}
(\pi_\ell*v_{T,\epsilon,a})(z)
=
(C-a)\int_{-T}^{0}\pi_\ell(z-y)dy
-
\frac{\rho}{\epsilon}
\int_{-\epsilon}^{0}\pi_\ell(z-y)dy
+
\int_0^z
\pi_\ell(z-y)\bigl(W(y)-a\bigr)dy.
\label{eq:approx_convolution_expansion}
\end{equation}
For the first term, the change of variables $u=z-y$ gives
\[
\int_{-\infty}^{0}\pi_\ell(z-y)\,dy
=
\int_z^\infty\pi_\ell(u)\,du
=
\Pi_\ell(z).
\]
Therefore, as $T\to\infty$, we have that
\[
\int_{-T}^{0}\pi_\ell(z-y)\,dy
\longrightarrow
\Pi_\ell(z)
\]
uniformly for $z\in[0,L]$.
For the second term, uniform continuity of $\pi_\ell$ on
$[0,L+1]$ gives
$\frac{1}{\epsilon}
\int_{-\epsilon}^{0}
\pi_\ell(z-y)\,dy
\longrightarrow
\pi_\ell(z)$
uniformly for $z\in[0,L]$.
Finally, we have that
\[
\int_0^z
\pi_\ell(z-y)\bigl(W(y)-a\bigr)dy
=
\int_0^z\pi_\ell(z-y)W(y)dy
-
a\int_0^z\pi_\ell(z-y)dy.
\]
Since it holds
\[
\int_0^z\pi_\ell(z-y)dy
=
\int_0^z\pi_\ell(u)du
=
1-\Pi_\ell(z),
\]
the limit of \eqref{eq:approx_convolution_expansion} is
\[
(C-a)\Pi_\ell(z)
-\rho\pi_\ell(z)
+\int_0^z\pi_\ell(z-y)W(y)\,dy
-a\bigl(1-\Pi_\ell(z)\bigr)
=
Q(z)-a.
\]
Thus, we proved that
$\pi_\ell*v_{T,\epsilon,a}
\longrightarrow
Q-a$
uniformly on $[0,L]$ as $T\to\infty$ and
$\epsilon\downarrow0$.

If $Q-a$ had more than $2r+2$ strict sign changes, then there would
exist finitely many ordered points at which its values are nonzero and
alternate in sign more than $2r+2$ times. Uniform convergence would
preserve all these signs for sufficiently large $T$ and sufficiently
small $\epsilon$, contradicting
\eqref{eq:erlang_variation_diminishing}. Therefore we show that
\begin{equation}
\label{eq:output_sign_change_bound}
\operatorname{sc}(Q-a)
\le
2r+2
\qquad
\text{for every }a>\zeta.
\end{equation}
If $\{Q>\zeta\}$ had more than $r+2$ components, choose one point in
$r+3$ of them and one separating point between each consecutive pair.
Choose $a>\zeta$ sufficiently close to $\zeta$ that $Q>a$ at all the
selected positive points.  At the separating points $Q-a<0$.  This gives
at least $2r+4$ strict sign variations of $Q-a$, a contradiction. Our proof is thus completed.

\subsubsection{Proof of \Cref{claim:no_repeated_birth}}
If $D_{i+1}(t)=0$, then we have
\[
W_i(t)=-\eta_i c_i(t)\le0\le a,
\]
so $W_i(t)>a$ is impossible. If $D_{i+1}(t)=1$, then the constraint
$D_{i+1}(t)\le c_i(t)\le D_i(t)$
forces $c_i(t)=1$, and therefore, it holds that
\[
W_i(t)=\mu_{s_{i+1}}(t)-\eta_i.
\]
This proves \eqref{eq:Wi_superlevel_identity}.
Fix a connected component $I$ of the set
$\{\mu_{s_{i+1}}>\eta_i+a\}$.
On $I$, we have that
\[
\Gamma_i(t)
=
\mu_{s_{i+1}}(t)
-
\eta_i\bI_{\{\chi_i(t)>0\}}
>0.
\]
By \eqref{eq:switch_safe}, we know that $\widehat\Phi_i$ is nondecreasing on $I$. By \Cref{claim:no_active_singular_arc}, it holds that
\[
D_{i+1}=1
\quad\Longleftrightarrow\quad
\widehat\Phi_i<0
\quad\text{a.e.}
\]
The strict negative set of a nondecreasing function has at most one connected component. Combining this fact with \eqref{eq:Wi_superlevel_identity} finishes the proof of the claim.

\subsection{Proof of \Cref{lem:gen_D2_connected}}
\label{app:proof_D2_connected}
We begin with an optimal $\beta$ representative satisfying the safe island bound from
\Cref{lem:gen_safe_component}. We then select, within the optimal set, a representative whose
$D_2$-support is connected.

Similar to the proof of \Cref{lem:gen_safe_component}, we start with an optimal cutoff-valued solution of the compatible
regularized problem.  We first rule out active ties in the outer switching
gap. Following the same procedure as the proof of \Cref{claim:no_active_singular_arc}, we can show that
%\begin{claim}[No active outer switching tie]
%For every regularized optimum,
%\[
%\operatorname{Leb}
%\left(
%\{D_2=1\}\cap\{\widehat\Phi_1=0\}
%\right)=0.
%\]
%Consequently, after changing the control on a null set,
\begin{equation}\label{claim:no_active_outer_tie}
D_2=1
\quad\Longleftrightarrow\quad
\widehat\Phi_1<0
\quad\text{a.e.}
\end{equation}
We next show that the strict negative set of the outer gap is connected.
Consider a positive-length interval $I$ on which mode $1$ is used.  On
$I$, since mode $1$ is used, we have that
\[
D_2=D_3=\cdots=D_m=0.
\]
The costate equation in \eqref{eq:corrected_costate} gives that
$p_{s_2}'(t)=p_{s_2-1}'(t)=0$,
so, from \eqref{eq:def_mu}, we know that the term
$\mu_{s_2}(t)=C_2+p_{s_2-1}(t)-p_{s_2}(t)$
is constant on $I$.
We now use a coordinate $u$ to denote the base-level active time, i.e., $u=U(t)$ where $U(t)$ is defined in \eqref{eq:def_U_selector}.

For every branch $r\ge2$, summing
\eqref{eq:adjacent_gap_derivative} while the actual mode is $1$ gives that (with the expression shows up in \Cref{claim:singular_slope})
\begin{equation}
\label{eq:outer_branch_u_derivative}
\frac{d}{du}
\left(
\widehat{\mathcal H}_r-\widehat{\mathcal H}_1
\right)
=
\mu_{s_2}
-
\eta_1\bI_{\{\chi_1>0\}}
-
\zeta
+
\tau r_\star(u).
\end{equation}
On the interval $I$, the first term of the right hand side of \eqref{eq:outer_branch_u_derivative} is constant.  The indicator can switch only from
one to zero and therefore creates only an upward jump.  Since the function $r_\star$
is strictly increasing, the right-hand side of
\eqref{eq:outer_branch_u_derivative} is strictly increasing in $u$.
All deeper branches have the same derivative on a mode-$1$ interval,
so their pairwise differences are constant there.  Their minimum, the function
$u\mapsto\widehat\Phi_1(t(u))$, following definition in \eqref{eq:full_switching_gap}, is consequently strictly convex.

Suppose $\{\widehat\Phi_1<0\}$ had two separated components.  If the
separating set were a single point, filling that point changes the control
only on a null set and merges the two components.  Otherwise, the closed
interval between the two components has positive length and satisfies
$\widehat\Phi_1\ge0$.  
On that interval, we have $D_2=0$ almost everywhere. 
%by the strict sign rule, so the preceding derivative calculation applies on the whole interval.  
In the $u$-coordinate, we know that $\widehat\Phi_1$ is strictly
convex.  A sublevel set of a convex function is an interval; equivalently,
if its two endpoint values are zero, strict convexity makes every interior
value negative.  This contradicts $\widehat\Phi_1\ge0$ on the separating
interval.  Hence, after a null-set modification, the strict negative set
has at most one component.  By \Cref{claim:no_active_outer_tie}, we have that
\begin{equation}
\label{eq:J2_regularized}
J_2\le1
\end{equation}
for every compatible regularized optimum.

We now let all perturbations vanish.  Choose a sequence
\[
\delta_n\downarrow0,
\qquad
\varepsilon_{i,n}\downarrow0\quad(i=2,\dots,m),
\qquad
\zeta_n\downarrow0,
\qquad
0<\tau_n<\frac{\zeta_n}{2},
\qquad
\tau_n\downarrow0,
\]
and choose an optimal cutoff-valued, late-filled solution
$(\bD^n,\bm c^n,\bx^n)$
of the corresponding regularized problem. For each $n$, the
component bounds hold
\[
J_2(\bD^n)\le1,
~~~
J_m(\bD^n)\le2,
\text{~~~and~~~}
J_i(\bD^n)
\le
4(m-i+1)+2,
~~~ i=3,\dots,m-1.
\]
Pad every level by zero-length intervals up to its stated maximum number
of components.  The active horizon and every collapsed-time endpoint
then belong to a fixed compact finite-dimensional set.  Passing to a
subsequence, all endpoints converge.
The convergence of the $\beta$ controls follows recursively.  Since
$D_2^n$ is the indicator of one interval, endpoint convergence gives
$D_2^n\to D_2^0$ almost everywhere and in $L^1$.  If
$D_{i-1}^n\to D_{i-1}^0$ in $L^1$, then we have
$L_{i-1}^n(t)=\int_0^tD_{i-1}^n(s)ds$
converges uniformly. Convergence of the collapsed endpoints of level
$i$ therefore gives $D_i^n\to D_i^0$ almost everywhere and in $L^1$.
Induction yields convergence of the entire $\beta$ schedule.
Because the $\alpha$ controls are chosen late-filled and the masses $\bA$
are fixed, the optional part of every $c_i^n$ is described by one cutoff
in the cumulative clock of
$\{D_i^n=1,\ D_{i+1}^n=0\}$.
The total lengths of these optional sets and the required optional masses
converge with the $\beta$ endpoints. Passing to a further subsequence, the
cutoff locations converge, and hence
$\bm c^n\to\bm c^0$ almost everywhere and in $L^1$.  The linear state
equations then give uniform convergence $\bx^n\to\bx^0$, and every
terminal and running expression converges.
Thus $(\bD^0,\bm c^0)$ is an exact optimum of the original fixed-$\bA$
problem, and the closed endpoint bounds give
\[
J_2(\bD^0)\le1,
~~~
J_m(\bD^0)\le2,
\text{~~~and~~~}
J_i(\bD^0)
\le
4(m-i+1)+2,
~~~ i=3,\dots,m-1.
\]
Our proof is completed.

\subsection{Proof of \Cref{prop:m3_finite_M}}
For a fixed feasible pair $(\bm A,\bm z)$, the active-time controls,
state trajectory, and all terminal and running expressions are
independent of $M$. The only difference between the finite-$M$
objective and its limiting counterpart is that the base covering term
$\frac{G_0(\bm D,\bm c,\bm x)}{H_K(M)}$
is replaced by
$\frac{G_0(\bm D,\bm c,\bm x)}{K}$.
Since $G_0(\bm D,\bm c,\bm x)$ is the expected number of allocated
units, we have that
\[
0\le G_0(\bm D,\bm c,\bm x)\le K.
\]
Therefore, for every feasible $(\bm A,\bm z)$,
\[
0
\le
\Theta_{\bm A}(\bm z)
-
\lim_{M\rightarrow\infty}\Theta_{\bm A}(\bm z)
\le
\frac{K-H_K(M)}{H_K(M)}.
\]
Minimizing over $\bm A$ and
$\bm z\in\mathcal Z(\bm A)$ does not change the bound. Our proof is completed from combining \Cref{lem:TightOblivious}, \Cref{thm:PoissonOPT}, \Cref{thm:gen_threshold_structure}, and using the fact that the value of $\PoisOPTRe_K(\bm{s}, M)$ decreases with $M$.

\end{APPENDICES}

\end{document}